\documentclass[1p]{elsarticle}
\usepackage{epsf}  
\usepackage{amsfonts, amssymb, amsthm, amsmath}  
\usepackage{graphicx}
\usepackage{psfrag}
\usepackage{hyperref}
\newtheorem{thm}{Theorem}  
\newtheorem*{lemmaa}{Lemma A3}  
\newtheorem*{defna1}{Definition A1}
\newtheorem*{defna2}{Definition A2}
\newtheorem*{claima}{Claim A4}

\newtheorem{cor}[thm]{Corollary}  
\newtheorem{lemma}[thm]{Lemma}  
\newtheorem{remark}[thm]{Remark}  
\newtheorem{defn}[thm]{Definition}  
\newtheorem{prop}[thm]{Proposition}  
\newtheorem{claim}[thm]{Claim}  
\newtheorem{example}[thm]{Example}  
\numberwithin{thm}{section}  
\def\pf{\noindent\emph{Proof: }}  
\def\stop{\hfill$\square$}  
\def\A{A_{\text{loc}}}
\providecommand{\C}[2]{\ensuremath {C^{#1,\underline{#2}}}}

\def\nT{T^{nice}}

\providecommand{\totl}[1]{\ensuremath{\lceil #1\rceil }}
\providecommand{\totb}[1]{\ensuremath{\underline {#1}}}

\newcommand{\Q}{y}
\newcommand{\w}{w}
\newcommand{\ex}{\bold}
\providecommand {\e}[1]{\mathfrak t^{#1}}

\newcommand{\T}{\mathbf{T}}

\DeclareMathOperator{\Dom}{Dom}
\DeclareMathOperator{\dist}{dist}
\DeclareMathOperator{\id}{id}
\DeclareMathOperator{\expl}{Expl}
\DeclareMathOperator{\coker}{Coker}
\newcommand{\dbar}{\bar{\partial}}

\newcommand{\exploded}{exploded }

\providecommand{\et}[2]{\ensuremath{\ex T^{#1}_{#2}}}
\providecommand{\lrb}[1]{\ensuremath{\left(#1\right)}}
\providecommand{\abs}[1]{\left\lvert #1\right\rvert}  
\providecommand{\norm}[1]{\left\lVert #1\right\rVert}
  
\author{Brett Parker\\  \texttt{brettdparker@gmail.com} }  
  
\title{Holomorphic curves in exploded manifolds: compactness}

\begin{document}
\begin{abstract}

 This paper establishes compactness results for the moduli stack of holomorphic curves in suitable exploded manifolds. This together with \cite{reg} and \cite{egw} allows the definition of Gromov-Witten invariants of these exploded manifolds. 

\end{abstract}

\maketitle
\tableofcontents

\section{Introduction}

This paper establishes compactness results for the moduli stack of holomorphic curves in exploded manifolds. An introduction to exploded manifolds may be found in \cite{iec}.  

The topology in which we shall establish compactness is the $\C\infty1$ topology, which is both a topology and a level of regularity. A $\C\infty1$ function can be viewed as a generalization of a function  on a manifold with cylindrical ends which is smooth and has exponential convergence at cylindrical ends.  In \cite{evc} it is shown that the moduli stack of holomorphic  curves has regularity $\C\infty1$ in the sense that in a neighborhood of any holomorphic curve for which the $\dbar$ equation is transversal, the moduli stack of holomorphic curves can be represented by $\C\infty1$ family of holomorphic curves. More generally, a virtual moduli space is constructed in \cite{egw} which is locally constructed using $\C\infty1$ families of curves. The compactness established in this paper allows us to know that components of this virtual moduli space are compact, and define Gromov-Witten invariants for suitable exploded manifolds.  

\

The topology in which the moduli space of curves is compact seems unnatural from the perspective of differential topology because of bubbling phenomena and node formation. Similar phenomena occur in the category of exploded manifolds, but unlike the smooth category, all the bubbling and node formation behavior can happen within connected smooth (or $\C\infty1$)  families of maps. Still, in talking about a converging sequence of holomorphic curves, we must deal with the fact that the domain of the curves may change.

We shall say that a sequence of $\C\infty1$ curves 
\[ f_{i}:(\ex C_{i},j_{i})\longrightarrow \ex B\]
 converges to a $\C\infty1$ curve $f$ if the following holds:

  There exists a
  sequence of  families of $\C\infty1$ curves,
  
  \[\begin{array}{ccc} (\ex {\hat C},j_i) & \xrightarrow{\hat f^i} &{\ex B}
  \\ \downarrow & &
  \\ \ex F&  & \end{array}\]
  and a sequence of points $p_{i}\in \ex F$ so that 
  \begin{enumerate}
\item The fiber of $(\ex {\hat C},j_{i})$ over $p^{i}$ is $(\ex C_{i},j_{i})$, and  the restriction of $\hat f^{i}$ to the fiber over $p^{i}$ is $f^{i}$.
 \item
  This sequence of families converges in $\C\infty1$ to the family
  
  \[\begin{array}{ccc} (\ex {\hat C},j) & \xrightarrow{\hat f} &{\ex B}
  \\ \downarrow & &
  \\ \ex F&  & \end{array}\]
  in the sense that the fiberwise complex structures $j_{i}$ converge in $\C\infty1$ to $j$ and the maps $\hat f^{i}$ converge in $\C\infty 1$ to $\hat f$
 \item
   The sequence of points $p^i$ in $\ex F$  converge to some point $p\in\ex F$ so that  $f$ is given by the restriction of $\hat f$ to the fiber over $p$.  
  
 \end{enumerate}

One scary aspect  of this definition is that the corresponding topology on the moduli stack of holomorphic curves is non-Hausdorff in the same way that the topology on an exploded manifold is non-Hausdorff. In particular all the points $p\in F$ for which $p_{i}\rightarrow p$ correspond to limits of our sequence of curves.  This is to be expected as we want to view the moduli stack of holomorphic curves as analogous to an exploded manifold rather than simply as a topological space.   
In  \cite{evc}, it is shown that the above topology on the moduli stack of $\C\infty1$ curves is the natural topology in which a substack is open if  every family of curves intersects it in a open subfamily.

\

To understand  technical assumptions on an almost complex structure and taming form used in this paper, it is useful to think of the smooth part $\totl{\ex B}$ of an exploded manifold $\ex B$ as a $C^{\infty}$ version of a locally ringed space or a `smooth manifold with singularities', with a natural notion of tangent spaces and differential forms defined in section \ref{strict taming}. In section \ref{dbar compatible section}, we define the notion of a $\dbar\log$ compatible almost complex structure $J$ on an exploded manifold. Lemma \ref{J embedding} gives that for any such $J$, $\totl{\ex B}$ has an almost complex structure, and may locally be regarded as a holomorphic subset of some almost complex manifold. This result allows standard analysis of (pseudo)-holomorphic curves in almost complex manifolds to be applied to holomorphic curves in $\totl{\ex B}$. Section \ref{strict taming} contains a definition of a taming form, which may be thought of as a symplectic form $\omega$ on $\totl{\ex B}$ with an extra positivity condition that ensures that there exists a $\dbar\log$ compatible almost complex structure $J$ on $\ex B$ so that $\omega$ is positive on holomorphic planes. It is proved in section \ref{dbar compatible section} that for such a taming form $\omega$, the space of $\dbar\log$ compatible almost complex structures tamed by $\omega$ is nonempty and smoothly $n$-connected for all $n$.

 In the case that $\ex B$ is a smooth manifold, a taming form is simply a symplectic form and any smooth almost complex structure is $\dbar\log$ compatible, so the technical conditions of sections \ref{strict taming} and \ref{dbar compatible section} are automatically satisfied.

\

Section \ref{local area bound section} deals with an important extra requirement for compactness of the moduli space of curves. This requirement is that the local area of a holomorphic curve must be bounded by topological information about that curve. Of course, if $\ex B$ is a compact manifold, then the area of a holomorphic curve is bounded by the integral of a taming form. Lemma \ref{immersion bound} proves that if there is an integral affine map $\totb{\ex B}\longrightarrow [0,\infty)^{N}$ which is injective restricted to each stratum of $\totl{\ex B}$, then the local area of any holomorphic curve is bounded by the integral of a taming form. Lemma \ref{immersion bound} also proves that if there is an integral affine immersion $\totb{\ex B}\longrightarrow \mathbb R^{N}$, then the local area of any holomorphic curve $f$ is bounded by the integral of a taming form plus some combinatorial information about the ends of the tropical curve $\totb{f}$.

\

Section \ref{estimates section} provides estimates for the local behavior of holomorphic curves. It ends with Proposition \ref{decomposition proposition} on page \pageref{decomposition proposition} which provides a decomposition of a holomorphic curve into a bounded number of regions with bounded derivative and conformal geometry and a bounded number of low energy annuli with behavior constrained by Lemma \ref{strong cylinder convergence} on page \pageref{strong cylinder convergence}.  

Section \ref{compactness section} contains the proof of compactness of the moduli space of curves with bounded genus, number of punctures, energy, and local area. As proved in section \ref{local area bound section}, the local area bound follows from topological conditions if there is an immersion of $\totb{\ex B}$ in some $\mathbb R^{N}$. A more general version of compactness involving the case of curves in a family of exploded manifolds is stated on page \pageref{completeness theorem}.

The technical assumptions required on the target $(\ex B,J)$ are that $\ex B$ is basic and complete, and the almost complex structure  $J$ is $\dbar\log$ compatible, and tamed by a taming form $\omega$. As explained in section 4 of \cite{iec},  the tropical structure of $\ex B$ associates a polytope to each point in $\totl{\ex B}$, and these polytopes are glued together to form the tropical part $\totb{\ex B}$ of $\ex B$. The assumption that $\ex B$ is basic means that any two of these polytopes are glued together in at most one way. The assumption that $\ex B$ is complete means that $\ex B$ is compact and each of these polytopes is complete. For further details, see definitions 4.6 and 3.15 of \cite{iec}.

\

 This compactness result can be viewed as a generalization of the compactness results for holomorphic curves found in \cite{gromov}, \cite{ruan}, \cite{Li}, and \cite{IP}. In the algebraic case, similar compactness results have been obtained by Abramovich and Chen in \cite{Chen} and \cite{acgw} and Gross and Siebert in \cite{GSlogGW}. A different approach to compactifying the moduli space of curves in the symplectic case is given by Ionel in \cite{IonelGW}.
 
 \
 
The following are some examples of almost complex exploded manifolds satisfying the technical requirements used in this paper.

\begin{itemize}
\item Any compact manifold with an almost complex structure tamed by a symplectic form will satisfy the above assumptions. 
\item The explosion $\expl M$ of any complex manifold $M$ with normal crossing divisors so that each divisor is an embedded submanifold will be basic, and the tropical part of $\expl M$ will admit an immersion into $[0,\infty)^{n}$, so local area bounds on holomorphic curves will follow from energy bounds. If $M$ is compact and there is a symplectic form on $M$ which tames the complex structure, then $\expl M$ will satisfy the above assumptions. This case is useful for defining and computing Gromov-Witten invariants relative to normal crossing divisors.
\item Given any compact  symplectic manifold $M$ with orthogonally intersecting codimension 2 symplectic submanifolds, we can construct an exploded manifold $\ex M$ using a similar construction to the explosion functor. (For further details of this construction see section 14 of \cite{elc}.)
On such an exploded manifold $\ex M$, the set of  almost complex structures $J$  tamed by the symplectic form and satisfying the above assumptions is nonempty and $n$-connected for all $n$. This case is useful to define Gromov-Witten invariants relative to these symplectic submanifolds.
\end{itemize}

\section{Taming forms} \label{strict taming}

Taming forms play the role of  symplectic forms in taming  almost complex structures on exploded manifolds. The reason that a separate definition is needed is as follows: The most natural definition of a symplectic form on $\ex B$ is a closed 2-form $\omega$ which is a symplectic form on each tangent space. Such symplectic forms are inappropriate for taming holomorphic curves, because the integral of $\omega$ will often be infinite on curves of interest. 

We shall use taming forms which may be regarded as symplectic forms on $\totl{ \ex B}$ with an extra positivity condition. For this to make sense, we shall define the tangent and cotangent spaces of $\totl{\ex B}$ within this section. Recall that every smooth map $\ex B\longrightarrow \mathbb R$ factors through $\totl{\ex B}$, so the smooth real valued functions on $\ex B$ may be regarded as pulled back from a sheaf of smooth functions on $\totl{\ex B}$. In fact, $\totl{\ex B}$ with this sheaf of smooth functions may be regarded as a $ C^{\infty}$ version of a locally ringed space, so we may define the tangent and cotangent spaces as usual.

\begin{defn}The tangent space $T_{p}\totl{\ex B}$ to the smooth part of $\ex B$ at a point $p$ is the vector space of derivations $v$ so that for any two smooth locally defined real valued functions $f$ and $g$ on $\totl{\ex B}$, 
\[v(fg)=v(f)g(p)+f(p)v(g)\]
\[v(f+g)=v(f)+v(g)\]
and for any constant real function $c$,
\[v(cf)=cv(f)\] 
\end{defn}

We shall sometimes also use the notation $T_{p}\totl{\ex B}$ when $p$ is a point in $\ex B$ to refer to the tangent space to $\totl{\ex B}$ at the image of $p$ in $\totl{\ex B}$.

\
 
The tangent space $T_{p}\totl{\ex B}$ may be higher or lower dimensional than $\ex B$, and $T\totl{\ex B}$ may not be a vector bundle. 
 If $\ex B$ is a single coordinate chart,  then $\totl{\ex B}$ is equal to some closed subset  of $\mathbb C^{n}$. Tangent vectors can then be identified with vectors on $\mathbb C^{n}$ which vanish on any smooth function on $\mathbb C^{n}$ which vanishes on $\totl{\ex B}$. Therefore, in this case, we may regard  $T_{p}\totl{\ex B}$ as a linear subspace of $T_{p}\mathbb C^{n}$ and  $T\totl{\ex B}$ as a closed subset of $T\mathbb C^{n}$. For example, the smooth part of $\et 1{[0,1]}$ is equal to the subset of $\mathbb C^{2}$ where $\zeta_{1}\zeta_{2}=0$. The tangent space at $(0,0)$ is $4$ (real) dimensional, but the tangent space elsewhere on $\totl{\et 1{[0,1]}}$ is $2$ (real) dimensional.

\

\begin{defn} Define the nice tangent space, $\nT_{p}\totl{\ex B}\subset T_{p}\totl{\ex B}$ to be the subspace of $T_{p}\totl{\ex B}$ consisting of vectors $v$ satisfying the following property: Given any finite collection of smooth functions $h_{i}$ and $f_{i}$ on $\ex B$, if the one form $\sum h_{i}df_{i}$ on $\ex B$ is $0$, then 
\[\sum_{i}h_{i}v(f_{i})=0\]
\end{defn}

Clearly, the nice tangent space is a linear subspace of $T_{p}\totl{\ex B}$. It is conceivable that $\nT_{p}\totl{\ex B}\neq T_{p}\totl{\ex B}$ for some particularly singular $\totl{\ex B}$, but I know of no such examples. In the case that $\totl{\ex B}$ is embedded as a closed subset of $\mathbb C^{n}$, $\nT \totl{\ex B}\subset T\mathbb C^{n}$ consists of vectors which vanish on all smooth one forms which pullback to be $0$ on $\ex B$. It shall be important for us that $\nT\totl{\ex B}\subset T\mathbb C^{n}$ is closed.   Note that given any smooth map $\ex A\longrightarrow \ex B$, there is an induced tangent map $T\totl{\ex A}\longrightarrow T\totl{\ex B}$ which sends the nice tangent space to the nice tangent space. The following lemma implies that smooth maps from manifolds to $\totl{\ex B}$ only see the nice tangent space of $\totl{\ex B}$.

\begin{lemma}\label{nice image}The derivative of any smooth map $f:\mathbb R\longrightarrow \totl{\ex B}$ has image contained in the nice tangent space to $\totl{\ex B}$.\end{lemma}

\pf

Note that a smooth map $f:\mathbb R\longrightarrow \totl{\ex B}$ means a map $f$ from $\mathbb R$ to $\totl{\ex B}$ so that $g\circ f$ is smooth for all smooth functions $g$ on $\totl{\ex B}$.
By restricting to an open subset, we may assume that $\totl{\ex B}$ is embedded as a closed subset of $\mathbb R^{n}$, so a smooth map to $\totl{\ex B}$ is just a smooth map to $\mathbb R^{n}$ contained in $\totl{\ex B}\subset\mathbb R^{n}$. The nice tangent space of $\totl{\ex B}$ is then a closed subset of $T\mathbb R^{n}$, so the inverse image under $df$ of the nice tangent space is a closed subset of $T\mathbb R$. 

At a point $p$ contained in a stratum $S$ of  $\totl{\ex B}$, any vector $v$ tangent to $S$ comes from a vector in $T_{p}\ex B$, therefore $T_{p}S\subset T_{p}\totl{\ex B}$ is  contained in the nice tangent space.

Suppose that $df$ restricted to $x\in\mathbb R$ has image outside of the nice tangent space of $\totl{\ex B}$. Then $df$ is not contained in the tangent space to the stratum of $\totl{\ex B}$ containing $f(x)$. It follows that for some neighborhood $U$ of $x$, $f$ is contained in a higher dimensional stratum of $\totl{\ex B}$ on $U-\{x\}$.  If on some neighborhood of $x$, $df$ is not contained in the nice tangent space, the argument may be continued to obtain points arbitrarily close to $x$ which are sent by $f$ to maximal dimensional strata of $\totl{\ex B}$. At such points $df$ must be contained in the nice tangent space. Because the set of points where $df$ is contained in the nice tangent space is closed, this contradicts our assumption that $df$  restricted to $x$  was not in the nice tangent space.

\stop

\begin{remark}In the case of a family $\pi_{\ex G}:\hat{\ex B}\longrightarrow \ex G$, we shall need to talk about the vertical tangent space $T_{vert}\hat{\ex B}$ of $\hat{\ex B}$, which is the sub bundle of $T\hat{\ex B}$ which is the kernel of $d\pi_{\ex G}$. Similarly, we shall need the vertical tangent space $\nT_{vert,p}\totl{\hat{\ex B}}\subset \nT_{p}\totl{\hat{\ex B}}$ which is the kernel of $d\totl{\pi_{\ex G}}$ intersected with the nice tangent space. Similarly, we shall need the vertical cotangent spaces $T^{*}_{vert}\hat{\ex B}$ and $T^{*}_{vert}\totl{\hat {\ex B}}$, and the vertical exterior differential operator $d_{vert}$.
\end{remark}

There is an obvious map 
\[\iota:T_{p}\ex B\longrightarrow \nT_{p}\totl{\ex B}\subset T_{p}\totl{\ex B}\] 
given by considering a derivation $v\in T_{p}\ex B$ on exploded functions as giving a derivation $\iota( v)\in T_{p}\totl{\ex B}$  on smooth functions. 

\begin{remark} In section \ref{dbar compatible section}, we shall construct $\dbar\log$ compatible almost complex structures $J$ which induce an almost complex structure $J'$ on  $\nT\totl{\ex B}$. These $\dbar\log$ compatible almost complex structures will be canonically determined on the kernel and cokernel of $\iota$ by the following two properties:   

\begin{enumerate}\item
If $v$ is in the kernel of $\iota$, and $\tilde z$ is any exploded function, 
\[(Jv)\tilde z=i(v\tilde z)\]

\item If $\tilde z$ is any exploded function so that $\totl{\tilde z}$ is defined close to $p$ and $\totl{\tilde z}$ vanishes at $p$, then for any $w\in \nT_{p}\ex B$,

\[(J'w)\totl{\tilde z}=i(w\totl{\tilde z})\]
\end{enumerate}

\end{remark}

To understand how the above could possibly hold, observe that if $h$ is any smooth $\mathbb C^{*}$-valued function, then $v(h\tilde z)=h(p)v\tilde z$ so long as $v$ is in the kernel of $\iota$, and $w\totl{h\tilde z}=h(p)w\totl{\tilde z}$ so long as $\totl{\tilde z}$ makes sense close to $p$ and vanishes at $p$.

 Note also that if $v$ is any vector in $T_{p}\ex B$, then it is proved in \cite{iec} that $v\tilde z$ is equal to a complex number times $\tilde z(p)$, so if $\totl{\tilde z(p)}=0$, then $\iota (v)\totl{\tilde z}=\totl{\iota(v)\tilde z}=0$. The second condition above may therefore be regarded as a condition on $J'$ acting on the cokernel of $\iota$.  The following lemma proves that there exists a unique complex structure on the cokernel of $\iota$ satisfying this condition.

\begin{lemma}\label{cokernel complex} There exists a unique complex structure $J$ on the cokernel of $\iota:T_{p}\ex B\longrightarrow \nT_{p}\totl{\ex B}$ so that given any exploded function $\tilde z$ so that $\totl{\tilde z}$ is defined close to $p$ and vanishes at $p$,  
then
\[J(d\totl{\tilde z}_{p})=i d\totl{\tilde z}_{p}\] 
where $d\totl{\tilde z}_{p}$ is the  element of $\mathbb C\otimes_{\mathbb R}(\coker\iota)^{*}$ which sends $w\in \nT_{p}\totl{\ex B}$ to $w\totl{\tilde z}$.

\end{lemma}

\pf

Assume with out losing any generality that $p$ is in a coordinate chart $\mathbb R^{n}\times \et mP$ and has tropical part contained in the interior of $P$. Any smooth function defined in a neighborhood of $p$ may be written as a smooth function of $\mathbb R^{n}$ and smooth monomials $\zeta_{k}$ on $\et mP$. Therefore the linear forms on $\nT_{p}\totl{\ex B}$ are generated by $dx_{i}$ and the real and imaginary parts of $d\zeta_{k}$ for some basis $\{\zeta_{k}\}$ of smooth monomials on $\et mP$. Each of the forms $d\zeta_{k}$ vanish on the image of $\iota$, and $dx_{i}$ are linearly independent on the image of $\iota$, therefore, the dual of the cokernel of $\iota$ is spanned by the real and imaginary parts of $d\zeta_{k}$. Therefore, there exists at most one complex structure on the cokernel of $\iota$ that sends the imaginary part of $d\zeta_{k}$ to the real part of $d\zeta_{k}$. 

In coordinates on  $\et mP$, there is a canonical complex structure $J_{0}$ so that  $\tilde z_{j}^{-1}dz_{j}\circ J_{0}= i\tilde z_{j}^{-1}d\tilde z_{j}$.  The one-form $d\zeta_{i}$ is holomorphic using this complex structure on $\et mP$. Therefore, 
\[\sum_{k}f_{k}d\zeta_{k}+h_{k}d\bar \zeta_{k}\]
 vanishes as a one form on $\et mP$ if and only if both $\sum_{k}f_{k}d\zeta_{k}$ and $\sum_{k}h_{k}d\bar\zeta_{k}$ vanish. Note that a linear combination $\sum_{k} a_{k}d\zeta_{k}+b_{k}d\bar\zeta_{k}$ vanishes on $\nT_{p}\totl{\ex B}$ if and only if there exists smooth functions $f_{k}$ and $h_{k}$ so that $f_{k}(p)=a_{k}$ and $h_{k}(p)=b_{k}$ and the above one-form vanishes.  If this happens, then the one-form $\sum_{k}if_{k}d\zeta_{k}-ih_{k}d\bar\zeta_{k}$ also vanishes, so $\sum_{k}a_{k}id\zeta_{k}-b_{k}id\bar \zeta_{k}$ also vanishes. It follows that there exists a  well defined complex structure $J$ on $ \coker\iota$ with the property that $Jd\zeta_{k}=id\zeta_{k}$ and $Jd\bar\zeta_{k}=-id\bar\zeta_{k}$.

If $\tilde z$ is any exploded function defined near $p$ so that $\totl{\tilde z}$ is well defined near $p$ and vanishes at $p$, then there exists a smooth monomial $\zeta_{k}$ and a smooth complex valued function $f$ so that $\totl{\tilde z}=f\zeta_{k}$ near $p$. Then $d\totl{\tilde z}_{p}$  is equal to $f(p)d\zeta_{k}$, so    
$J(d\totl{\tilde z}_{p})=id\totl{\tilde z}_{p}$ as required.

\stop

\begin{defn} A taming form on $\nT_{p}\totl{\ex B}$ is a symplectic form on the vector space $\nT_{p}\totl{\ex B}$ so that 
\begin{itemize}
\item The image of 
\[\iota:T_{p}\ex B\longrightarrow \nT_{p}\totl{\ex B}\] is symplectic.
\item Identifying the cokernel of $\iota$ with the symplectic orthogonal of $\iota(T_{p}\ex B)\subset \nT_{p}\totl{\ex B}$ gives a symplectic form on the cokernel of $\iota$. This symplectic form must be positive on the  holomorphic planes of the canonical complex structure on the cokernel of $\iota$.
\end{itemize}
\end{defn}

\

We may define the cotangent space of $\totl{\ex B}$ at a point $p$ to be the dual to $T_{p}\totl{\ex B}$.  Given any smooth function $h$ on $\totl{\ex B}$, there is a section $dh$ of this cotangent bundle of $\totl{\ex B}$ defined as usual by 
\[dh(v):=v(h)\]
Below we shall define smooth differential forms on $\totl{\ex B}$. Be warned that although $d$ of a function $h$ is well defined, the exterior derivative of such differential forms on $\totl{\ex B}$ is not necessarily well defined at points where $\totl{\ex B}$ is not smooth.

\begin{defn} A smooth differential $n$-form $\theta$ on $\totl{\ex B}$ is a choice of  $n$-form $\theta_{p}$ on $T_{p}\totl{\ex B}$ for each $p\in\totl{\ex B}$ so that each point has a neighborhood $U$ with a finite number of smooth functions $h_{i,j}$ so that 
\[\sum_{i}h_{i,0} dh_{i,1}\wedge \dotsb \wedge dh_{i,n}\]
restricted to each point $p$ in $U$ is equal to $\theta_{p}$. 

Given a submersion $\pi:\hat {\ex B}\longrightarrow \ex G$, a smooth family $\theta$ of $n$-forms on $\hat{\totl{\ex B}}$ is  for each $p$ a choice of $n$-form $\theta_{p}$ on the kernel of $d\totl\pi: T_{p}\totl{\ex B}\longrightarrow T_{\pi(p)}\totl{\ex G}$, so that each point has a neighborhood $U$ and a finite number of smooth functions $h_{i,j}$ so that 
\[\sum_{i}h_{i,0} dh_{i,1}\wedge \dotsb \wedge dh_{i,n}\]
restricted the vertical tangent space of $\totl{\ex B}$ at each point $p$ in $U$ is equal to $\theta_{p}$.
\end{defn}

The local description of $T\totl{\ex B}$ as a subset of $T\mathbb C^{n}$ implies that we may also think of smooth $n$-forms on $\totl{\ex B}$ as the restriction of smooth $n$-forms on $\mathbb C^{n}$. This implies that any $n$-form on  $T_{p}\totl{\ex B}$ may be extended to a smooth $n$-form on $\totl{\ex B}$. 

Given any smooth $n$-form $\theta$ on $\totl{\ex B}$, we may pull back $\theta$ via the map $\iota:T\ex B\longrightarrow T\totl{\ex B}$ to get a section $\iota^{*}\theta$ of $\bigwedge^{n}T^{*}\ex B$. As $\iota^{*}$ commutes with exterior differentiation of functions, sums and wedge products, $\iota^{*}\theta$ is a smooth differential $n$-form on $\ex B$. We shall say that such differential forms are generated by functions.

\begin{defn}\label{genbf}Say that a (real or complex valued) differential form on $\ex B$ is generated by functions if it is locally equal to a differential form constructed from smooth (real or complex valued) functions on $\ex B$ by the operations of exterior differentiation, wedge products and sums.

In the case of a family $\hat{\ex B}\longrightarrow \ex G$, say that a section of $\bigwedge T_{vert}^{*}\hat{\ex B}$ is generated by functions if it is locally constructible from smooth functions using the operations of $d_{vert}$, wedge products and sums.
\end{defn}

Note that any smooth  $n$-form on $\totl{\ex B}$ defines a $n$-form on $\nT_{p}\totl{\ex B}\subset T_{p}\totl{\ex B}$. The defining feature of $\nT_{p}\totl{\ex B}$ is that any smooth one-form which is nonzero on $\nT_{p}\totl{\ex B}$ pulls back to be a one-form generated by functions on $\ex B$ which is nonzero in a neighborhood of $p$. It follows that any one-form generated by functions lifts uniquely to the nice cotangent space of $\totl{\ex B}$.

\begin{defn} A taming form on $\ex B$ is a smooth two-form $\omega$ on $\totl{\ex B}$ so that the two-form $\iota^{*}\omega$ on $\ex B$ is closed, and so that $\omega$  restricts to give  a taming form on $\nT_{p}\totl{\ex B}$ for all $p\in \totl{\ex B}$.

A family of taming forms on $\hat{\ex B}\longrightarrow \ex G$ is a smooth family $\omega$ of two-forms on $\totl{\hat{\ex B}}$  which is a taming form restricted to each fiber of $\hat{\ex B}\longrightarrow \ex G$.
\end{defn}

For example, if $\{\zeta_{j}\}$ is a basis for the smooth monomials on $\et mP$, then 
\[\sum_{j}id\zeta_{j}\wedge d\bar \zeta_{j}\]
is a taming form on $\et mP$.

\

Where there is no possibility for confusion, we shall refer to the two-form $\iota^{*}\omega$ on $\ex B$ simply as $\omega$.

\section{$\dbar\log$ compatible almost complex structures}\label{dbar compatible section}

We shall make a technical assumption on our almost complex structure which will be  weak enough that every taming form will tame a nonempty connected space of such almost complex structures in the sense of Definition \ref{tame definition} below.

 The condition on our almost complex structure shall be that  $f^{-1}\dbar f$ is generated by functions for any exploded function $f:\ex B\longrightarrow \ex T$.  We may consider $f^{-1}df$ as a $\mathbb C$ valued $1$-form on $\ex B$ with real and imaginary parts the pullback from $\ex T$ of the real and imaginary parts of $\tilde z^{-1}d\tilde z$. We can therefore define

\[\dbar\log f:=f^{-1}\dbar f:=\frac 12(f^{-1}df+if^{-1}df\circ J)\]

\begin{defn}
Say that an almost complex structure $J$ on $\ex B$ is $\dbar\log$ compatible if every locally defined exploded function $f$ satisfies the condition that $f^{-1}\dbar f$ is generated by functions.

Say that a family of almost complex structures $J$ on a family $\hat {\ex B}$ is $\dbar\log$ compatible if for every locally defined exploded function $f$, 
\[f^{-1}\dbar f:=\frac 12(f^{-1}d_{vert}f+if^{-1}d_{vert} f\circ J)\]
is generated by functions.
\end{defn}

Note in particular that on a smooth manifold, every smooth almost complex structure is $\dbar\log$ compatible. Also, if $(\ex B,J)$ is a holomorphic exploded manifold, every smooth exploded function is locally in the form $f=g\e a\tilde z^{\alpha}$ where $\tilde z^{\alpha}$ is holomorphic and $g$ is a $\mathbb C^{*}$ valued smooth function. In this case, $\dbar \log f=g^{-1}\dbar g$ which is generated by functions, so in this case $J$ is $\dbar\log$ compatible.

\

 The following is a useful criteria for checking that an almost complex structure is $\dbar\log$ compatible.

\begin{lemma} \label{J construction}

Let $J_{0}$ be a $\dbar\log$ compatible almost complex structure. Then another almost complex structure $J$ is $\dbar\log$ compatible if and only if $J-J_{0}$ sends smooth one-forms to one-forms generated by functions.

\end{lemma}

\pf

The difference between $\dbar\log f$ using $J$ and $J_{0}$ is 
\[i f^{-1}df\circ(J-J_{0})\]
If $(J-J_{0})$ sends smooth one-forms to one-forms generated by functions,  it follows that $\dbar \log f$ using $J$ is generated by functions because $\dbar \log f$ using $J_{0}$ is generated by functions. If $J$ and $J_{0}$ are both $\dbar\log $ compatible, then the above expression implies that $(J-J_{0})$ sends the real and imaginary parts of $f^{-1}df$ to one-forms generated by functions. As any smooth one-form may be expressed as a sum of smooth functions times real or imaginary parts of some $f^{-1}df$, it follows that $(J-J_{0})$ sends smooth one-forms to one-forms generated by functions.

\stop

\begin{remark}\label{preserves gbf} The standard complex structure $J_{0}$ on $\mathbb C^{n}\times \et mP$ is $\dbar\log$ compatible, and sends one forms generated by functions to one forms generated by functions, so Lemma \ref{J construction} implies that any $\dbar\log$ compatible almost complex structure sends one-forms generated by functions to one-forms generated by functions. 
\end{remark}

Recall that any one-form generated by functions lifts uniquely to the nice cotangent space of $\totl{\ex B}$. The following lemma shows that the action of $J$ on the sheaf of one-forms generated by functions induces a unique almost complex structure on the nice  tangent space of $\totl{\ex B}$.

\begin{lemma} \label{J lift}Given any $\dbar\log$ compatible almost complex structure $J$ on $\ex B$, for all $p\in\totl{\ex B}$, there exists a unique complex structure $J'$ on $\nT_{p}\totl{\ex B}$ satisfying the following property:

If $\theta$ is any one-form generated by functions, and $\theta_{p}$ is its lift to a  linear form on $\nT_{p}\totl{\ex B}$, then 
\[(\theta\circ J)_{p}=\theta_{p} \circ J'\]

\

For such a $J'$, the map \[\iota:T_{p}\ex B\longrightarrow \nT_{p}\totl{\ex B}\] is complex, and the induced complex structure on the cokernel of $\iota$ is the canonical complex structure from Lemma \ref{cokernel complex}.
\end{lemma}

\pf

Suppose that $p$ is contained in the stratum of a coordinate chart $\mathbb R^{n}\times \et mP$ corresponding to the interior of $P$. Let $\{\zeta_{k}\}$ be some basis of nonconstant smooth monomials on $\et mP$. The sheaf of $\mathbb C$-valued one-forms generated by functions is generated by $dx_{1},\dotsc,dx_{n}$,  $d\zeta_{k}$ and $d\bar\zeta_{k}$, so the space of complex valued linear forms on $\nT_{p}\totl{\ex B}$ is generated as a complex vector space by $(dx_{j})_{p}$, $(d\zeta_{k})_{p}$ and $(d\bar\zeta_{k})_{p}$.

 Because $\zeta_{k}$ is the smooth part of some exploded function $\tilde z$, 
 \[d\zeta_{k}\circ J=id\zeta_{k}- 2i\zeta_{k}\dbar\log\tilde z\]
The fact that $\zeta_{k}(p)=0$ then implies that
\[(d\zeta_{k}\circ J)_{p}=i(d\zeta_{k})_{p}\]
Therefore $J$ exchanges the  imaginary and real parts of $(d\zeta_{k})_{p}$, and
\[(d\bar\zeta_{k}\circ J)_{p}=-i(d\bar\zeta_{k})_{p}\]

As argued in the proof of Lemma \ref{cokernel complex}, the complex linear relations between our generators are generated by linear relations involving only the $(d\zeta_{k})_{p}$, and linear relations involving only the $(d\bar\zeta)_{p}$.

Therefore, we may define $J'$ so that
\[(d\zeta_{k})_{p}\circ J'=i(d\zeta_{k})_{p}=(d\zeta_{k}\circ J)_{p}\]
\[(d\bar\zeta_{k})_{p}\circ J'=-i(d\bar\zeta_{k})_{p}=(d\bar\zeta_{k}\circ J)_{p}\]
\[(dx)_{p}\circ J'=(dx\circ J)_{p}\]

Such a $J'$ has the required property that for any one form  $\theta$ generated by functions $\theta_{p}\circ J'=(\theta\circ J)_{p}$. Note that this property uniquely defines $J'$. This property also implies that  $\iota:T_{p}\ex B\longrightarrow \nT_{p}\totl{\ex B}$ is complex because $\iota^{*}(\theta_{p})$ is equal to $\theta$ restricted to $p$, so $\iota^{*}(\theta_{p}\circ J')=(\iota^{*}\theta_{p})\circ J$.  Note also that $J'$ induces the canonical complex structure on the cokernel of $\iota$ because $(d\zeta_{k})_{p}\circ J'=i(d\zeta_{k})$.

\stop

\

Our taming form $\omega$ is defined as a smooth two form on $\totl{\ex B}$ which is a symplectic form on $\nT_{p}\totl{\ex B}$ for all $p\in\totl{\ex B}$. Our complex structure $J$ is defined as an endomorphism of $T\ex B$, so the definition below of $J$ being tamed by $\omega$ uses the canonical  lift $J'$ of $J$ to $\nT_{p}\totl{\ex B}$.

\begin{defn}\label{tame definition}A $\dbar\log$ compatible almost complex structure $J$ is tamed by a taming form $\omega$ if for all $p\in\totl{\ex B}$, the lift $J'$ of $J$ to  $\nT_{p}\totl{\ex B}$ is tamed by $\omega$ in the sense that  $\omega$ is positive on every $J'$-holomorphic plane in $\nT_{p}\totl{\ex B}$. 

A family of $\dbar\log$ compatible almost complex structures is tamed by a family of taming forms if restricted to each fiber, the corresponding almost complex structure is tamed by the corresponding taming form.
 \end{defn}

\begin{example}\end{example}
Recall from \cite{iec} that the explosion functor constructs an exploded manifold $\expl M$ from a complex manifold with normal crossing divisors $M$. If $\omega$ is any symplectic form on $M$ which tames the complex structure, then the pullback of $\omega$ to $\expl M$ via the smooth part map $\expl M\longrightarrow M$ is a taming form on $\expl M$ that tames the complex structure.  
 
\

\begin{remark}There is a natural topology on the nice tangent space of $\totl{\ex B}$ which is the coarsest  topology in which the projection to $\totl{\ex B}$ and all one-forms generated by functions are continuous. With this topology, if $\totl{\ex B}$ locally embeds as a closed subset of $\mathbb R^{n}$, then the nice tangent space of $\totl{\ex B}$ locally embeds as a closed subset of $T\mathbb R^{n}$. In this topology, $J'$ is continuous (and we show in Lemma \ref{J embedding} that $\mathbb R^{n}$ may be chosen so that $J'$ extends to an almost complex structure on $\mathbb R^{n}$). Similarly, any  two-form $\omega$ on $\totl{\ex B}$  locally extends to a two-form on $\mathbb R^{n}$. It follows that if $\omega$ is positive on $J'$-holomorphic planes in $\nT_{p}\totl{\ex B}$, then $\omega$ is also positive on $J'$-holomorphic planes in $\nT_{p'}\totl{\ex B}$ for $p'$ in a neighborhood of $p$.\end{remark}

\

\

We may use the following method to average two  $\dbar\log$ compatible almost complex structures tamed by a given symplectic form.

\begin{lemma}\label{average J}

Suppose that $J$ and $J_{0}$ are two families of $\dbar\log$ compatible almost complex structures tamed by $\omega$ on $\hat{\ex B}$, and $\rho:\hat{\ex B}\longrightarrow [0,1]$ is smooth. The following construction is well defined and produces a family of $\dbar\log$ compatible almost complex structures $J_{\rho}$.
 
Let $Z$ be the endomorphism of $T_{vert}\hat{\ex B}$ defined by 
\[Z:=-J\circ J_{0}\]
then define 
\[W:=(1+Z)^{-1}(1-Z)\]
and 
\begin{equation}\label{Jdef}J_{\rho}:=(1+\rho W)^{-1}(1-\rho W)J_{0}\end{equation}

\end{lemma}

\pf

 Restricted to $T_{vert,p}\hat{\ex B}$, the endomorphism $1+Z$ is invertible because $1+Z$ is twice the identity on the kernel of $\iota:T_{vert,p}\hat{\ex B}\longrightarrow \nT_{vert,p}\totl{\hat{\ex B}}$, and if $v$ is any vector not in this kernel, then $\omega(v,(1+Z)J_{0}v)>\omega(v,J_{0}v)>0$.
Therefore, the following defines a smooth endomorphism
\[W:=(1+Z)^{-1}(1-Z)\]
We may recover $Z$ from $W$ by the same formula
\[Z=(1+W)^{-1}(1-W)\]
Note that \[(1+\rho W)=(1+Z)^{-1}\lrb{(1+\rho)+(1-\rho)Z}\]
$\lrb{(1+\rho)+(1-\rho)Z}$ is equal to twice the identity on the kernel of $\iota$, and on any vector $v$ not in this kernel, $\omega(v,\lrb{(1+\rho)+(1-\rho)Z}J_{0}v)>0$. It follows that $(1+\rho W)$ is invertible, and the following formula defines a smooth endomorphism $J_{\rho}$ of $T_{vert}\hat{\ex B}$.
\[J_{\rho}:=(1+\rho W)^{-1}(1-\rho W)J_{0}\]

To check that $J_{\rho}\circ J_{\rho}=-1$, note that $J^{2}=-1$ is equivalent to $Z^{-1}=J_{0}ZJ^{-1}_{0}$ which is equivalent to $-W=J_{0}W J^{-1}_{0}$, which also holds for $\rho W$. It follows that $J_{\rho}\circ J_{\rho}=-1$, so $J_{\rho}$ defines a complex structure on $T_{vert,p}\hat{\ex B}$ for each $p$. As $J_{\rho}$ is equal to $J_{0}$ when restricted to vectors in the kernel of $\iota$, $J_{\rho}$ is a family of almost complex structures in the sense of \cite{iec}. 

We shall now show that $J_{\rho}-J_{0}$ sends smooth one-forms to one-forms generated by functions, so Lemma \ref{J construction} will imply that $J_{\rho}$ is $\dbar\log$ compatible.
\[\begin{split}J_{\rho}-J_{0}&=(1+\rho W)^{-1}(-2\rho W)J_{0}
\\ &=-2 \rho (1+\rho W)^{-1}(1+Z)^{-1}(1-Z)J_{0}
\\ &=2 \rho (1+\rho W)^{-1}(1+Z)^{-1}(J-J_{0})\end{split}\]

The dual action of $J_{\rho}-J_{0}$ on $T^{*}_{vert}\ex B$ is the composition of the dual of $2 \rho (1+\rho W)^{-1}(1+Z)^{-1}$ with the dual of $(J-J_{0})$. The first action sends smooth one-forms to smooth one-forms, the second action sends these smooth one-forms to one-forms generated by functions. Therefore, $(J_{\rho}-J_{0})$ sends smooth one-forms to one-forms generated by functions, so Lemma \ref{J construction} implies that $J_{\rho}$ is a family of $\dbar\log$ compatible almost complex structures.

\stop

\begin{lemma}
\label{average J'}Let $J$, $J_{0}$ and $J_{\rho}$ be as in Lemma \ref{average J}, and let $J'$, $J_{0}'$ and $J'_{\rho}$ be the canonical lifts of $J$, $J_{0}$ and $J_{\rho}$ to $\nT_{vert, p}\totl{\hat{\ex B}}$.  Then  $J'_{\rho}$ satisfies the following equation

\[Z':=-J'\circ J'_{0}\] 
\[W':=(1+Z')^{-1}(1-Z')\]
\begin{equation}\label{J'def}J'_{\rho}=(1+\rho(p) W')^{-1}(1-\rho(p) W')J'_{0}\end{equation}

\end{lemma}

\pf

\

 As this construction may be restricted to a fiber of $\hat{\ex B}\longrightarrow \ex G$, we may simplify notation by  restricting to the case of a single target $\ex B$. Given any one-form  $\theta$ generated by functions, use the notation $\theta_{p}$ for the corresponding linear form on $\nT_{p}\totl{\ex B}$. Recall that the complex structure $J'_{\rho}$ induced on $\nT_{p}\totl{\ex B}$ by $J_{\rho}$ satisfies the following defining property: 
\[\theta_{p}\circ J'_{\rho}=(\theta\circ J_{\rho})_{p}\]
 The induced complex structures $J'$ and $J_{0}'$ also satisfy similar equations. 
Note that
\[\theta_{p}\circ Z'=-\theta_{p}\circ J'\circ J_{0}'=-(\theta\circ J)_{p}\circ J_{0}'=-(\theta\circ J\circ J_{0})_{p}=(\theta\circ Z)_{p}\]
 Similarly, 
\[\theta_{p}\circ W'=(\theta\circ W)_{p}\]
As all terms in the formula for $J_{\rho}$ obey similar equations, it follows that
\[J_{\rho}'=(1+\rho(p) W')^{-1}(1-\rho(p) W')J'_{0}\]

\stop

\

The following lemma shows that the set of $\dbar\log$ compatible almost complex structures tamed by a given taming form is nonempty and smoothly $n$-connected for all $n$. 

\begin{lemma}
The set of $\dbar\log$ compatible  almost complex structures tamed by a given taming form is nonempty and connected. 

More generally, suppose that $\omega$ is a family of taming forms on  $\hat {\ex B}$, and $J$ is a family of  $\dbar\log$ compatible  almost complex structures  tamed by $\omega$,  defined on a neighborhood of some compact subset $K$ of $\hat{ \ex B}$. Then there exists a family of $\omega$-tame  $\dbar\log$ compatible almost complex structures $J_{1}$ on $\hat{\ex B}$ so that $J_{1}$ and $J$ coincide on a neighborhood of $K$. 
\end{lemma}

\pf

Note that we only need to prove  the general case  of  a family of taming forms on $\hat{\ex B}\longrightarrow \ex G$, as the special case of a single target $\ex B$ follows from the general case applied to $\hat {\ex B}=\ex B$, $K=\emptyset$  for existence of a tamed $J$, and $\hat {\ex B}=\ex B\times \mathbb R$, $K=\ex B\times\{0,1\}$ for connectedness of the space of tamed $J$.

\

 Suppose that $K$ contains all strata of $\hat {\ex B}$ with tropical dimension bigger than $m$. We shall show that $J$  may be extended to a   family  of $\omega$-tamed,  $\dbar\log$ compatible almost complex structures $J_{1}$  defined on a neighborhood of the union of $K$ with all strata with tropical dimension $m$. The lemma follows from this by induction.

The proof shall proceed as follows: First, we shall construct an appropriate almost complex structure $J_{0}$ on a neighborhood of a stratum with tropical dimension $m$, then we shall patch together $J$ and $J_{0}$ using Lemma \ref{average J}, and verify using Lemma  \ref{average J'} that the resulting almost complex structure $J_{\rho}$ is tamed by $\omega$.

\

Consider a stratum of $\hat{\ex B}$ with tropical dimension $m$. The existence of equivariant coordinate charts from  \ref{equivariant section} compatible with the submersion $\hat{\ex B}\longrightarrow \ex G$  implies that a neighborhood of this stratum is isomorphic to some $\et mP$ bundle $\et mP\rtimes M$ over a manifold $M$, and that the map $\et mP\rtimes M\longrightarrow \ex G$ factors as an equivariant bundle map followed by an isomorphism onto an open subset of $\ex G$
\[\et mP\rtimes M\longrightarrow \et {m'}Q\rtimes M'\hookrightarrow \ex G\]

The map $M\longrightarrow M'$ is a submersion, so it makes sense to talk of the vertical tangent space $T_{vert}M$ of $M$. It also makes sense to talk about the vertical tangent spaces $T_{vert}\et mP$ and $\nT_{vert}\totl{\et mP}$  of $\et mP$, because the intersection of $T_{vert}\hat{\ex B}$ or $\nT_{vert}\totl{\hat{\ex B}}$ with the tangent space to each $\et mP$ fiber of $\et mP\rtimes M\longrightarrow M$ is independent of which fiber is chosen, and equivariant under the $(\mathbb C^{*})^{m}$ action on the fibers. Because on each fiber the map $\et mP\longrightarrow \et {m'}Q$ is given by monomials in coordinate functions,  $T_{vert}\et mP$ is a complex subbundle of $T\et mP$ when $\et mP$ is given the canonical complex structure in which each smooth monomial is complex. Similarly,  $\nT_{vert,p}\totl{\et mP}$ is a complex subspace of $\nT_{p}\totl{\et mP}$.

For any $p$ in our stratum, $\nT_{vert,p}\totl{\hat{\ex B}}$  splits into $T_{vert,p}M\times \nT_{vert,p}\totl{\et mP}$, with the inclusion of $T_{vert,p}M$ canonically given by the inclusion of $M$ as our stratum of $\totl{\hat{ \ex B}}$, and the projection to $T_{vert,p}M$ given by the derivative of the bundle map $\et mP\rtimes M\longrightarrow M$. We must modify this splitting to be compatible with a splitting given by the symplectic structure.

\begin{claim}\label{orthogonal split}After a coordinate change, we may assume that the splitting of $\nT_{vert,p}\totl{\hat {\ex B}}$ will  be into  the image of $T_{vert,p}\hat{\ex B}$ and its symplectic orthogonal. \end{claim}

To prove Claim \ref{orthogonal split},  consider a coordinate change locally in the form of $(\tilde z,x)\mapsto (\tilde z,f(\tilde z,x))$, where $f:\et mP\rtimes M\longrightarrow M$ is smooth, preserves fibers of $\totl{\hat {\ex B}}\longrightarrow \totl{\ex G}$, and is given by  the identity restricted to our stratum. We may choose $f$ so that if $p$ is in our stratum, and $v\in \nT_{vert,  p}\totl{\hat{\ex B}}$ is $\omega$-orthogonal to the image of  $T_{vert,p}\hat{\ex B}$, then $df(v)=0$. After this change of coordinates, the splitting of  $\nT_{vert,p}\totl{\hat {\ex B}}$ will  be into  the image of $T_{vert,p}\hat{\ex B}$ and its symplectic orthogonal as required. This completes the proof of Claim \ref{orthogonal split}.
 
 \
 
We may extend our splitting of $\nT_{vert}\totl{\hat{\ex B}}$ on our stratum to a splitting of $T_{vert}(\et mP\rtimes M)$ by choosing a `horizontal' sub-bundle of $T_{vert}(\et mP\rtimes M)$ which projects isomorphically to $TM$. Choose this horizontal sub-bundle to be invariant under the action of $(\mathbb C^{*})^{m}$ on the fibers of $\et mP\rtimes M$. This horizontal sub-bundle gives us compatible splittings of $T_{vert}(\et mP\rtimes M)$ into $T_{vert}\et mP\oplus T_{vert}M$, and $\nT_{vert}\totl{\et mP\rtimes M}$ into $\nT_{vert}\totl{\et mP}\oplus T_{vert}M$.

 \
 
Note that the restriction of $\omega$ to our stratum gives a family of symplectic forms on $M$. Let $J_{M}$ be any family of almost complex structures on $M$ tamed by this family of symplectic forms. $J_{M}$ lifts uniquely to a split almost complex  structure $J_{0}$ on $T_{vert}\et mP\rtimes M$ which is the canonical complex structure from $T_{vert}\et mP$ on fibers. 

To see that $J_{0}$ is $\dbar\log$ compatible, consider $\dbar\log \tilde z_{i}$ where $\tilde z_{i}$ is a coordinate function from $\et mP$ on a locally defined $\et mP\times \mathbb R^{s}$ chart on $\et mP\rtimes M$. As $J_{0}$ is the canonical complex structure on fibers, $\dbar\log \tilde z_{i}$ vanishes on $T_{vert}\et mP$. As $J_{0}$ and $d_{vert}\log \tilde z_{i}$ are invariant under the action of $(\mathbb C^{*})^{m}$ on $\et mP$, $\dbar \log \tilde z_{i}$ is invariant under this action of $(\mathbb C^{*})^{m}$, and is therefore the pullback of some family of forms from $\mathbb R^{s}$. Therefore $\dbar\log \tilde z_{i}$ is generated by functions. As any exploded function on $\et mP\times \mathbb R^{s}$ is a finite product of such $\tilde z_{i}$ with a smooth function, it follows that $\dbar \log$ of any exploded function is generated by functions, so $J_{0}$ is $\dbar\log$ compatible.

\
\begin{claim}\label{J0} The $\dbar\log$ compatible almost complex structure $J_{0}$ is tamed by $\omega$ in a neighborhood of our stratum. 
\end{claim}

For any $p$ in our stratum, consider the canonical lift $J_{0}'$ of $J_{0}$ to $\nT_{vert,p}\totl{\hat{\ex B}}$. As $J_{0}$ is a split almost complex structure, it preserves the space of one-forms on $\et mP\rtimes M$ which are pulled back from $M$. It follows that $J_{0}'$ must preserve the splitting of $\nT_{vert,p}\hat{\ex B}$ into $\nT_{vert,p}\totl{\et mP}\times T_{vert,p}M$. Note that as the $T_{vert,p}M$ factor is the image of the  complex map $\iota$, $J_{0}'$ must coincide with $J_{M}$ on this factor and therefore be positive on holomorphic planes. On the complimentary factor, $\nT_{vert,p}\totl{\et mP}$, $J_{0}'$ is determined by the canonical complex structure on the cokernel of $\iota$. Claim \ref{orthogonal split} gives that $\nT_{vert,p}\totl{\et mP}$ is the symplectic orthogonal to the image of $\iota$. By the definition of a taming form, $\omega$ must be positive on holomorphic planes within the symplectic orthogonal of the image of $\iota$ when the symplectic orthogonal is given the canonical complex structure of the cokernel of $\iota$. We therefore get that  $\omega$ is positive on $J'_{0}$-holomorphic planes within two symplectic orthogonal subspaces of $\nT_{vert,p}\totl{\hat{\ex B}}$. 
  Therefore   $\omega$ is positive on all $J'_{0}$-holomorphic planes within $\nT_{vert,p}\totl{\hat{\ex B}}$. It follows that $\omega$ is positive on $J_{0}'$-holomorphic planes in some neighborhood of our stratum, so $J_{0}$ is tamed by $\omega$ on a neighborhood of our stratum. This completes the proof of Claim \ref{J0}

\
 
 We shall now patch together $J_{0}$ and $J$ using Lemma \ref{average J}.
Choose a smooth function $\rho:\hat{\ex B}\longrightarrow [0,1]$ supported in the region where $J$ is defined and identically equal to $1$ in a neighborhood of $K$. On the subset of our stratum where $\rho>0$, we may choose $J_{M}$ to coincide with the almost complex structure that $J$ defines on $M$. Define $J_{\rho}$ on a neighborhood of $K$ and our stratum  by averaging $J_{0}$ and $J$ using $\rho$ as in Lemma \ref{average J}.

\begin{claim}\label{J tame}$J_{\rho}$ is tamed by $\omega$ in a neighborhood of our stratum.\end{claim}
 
 Let $p$ be any point in our stratum where $\rho(p)\in(0,1)$. Denote by $J'$, $J_{0}'$ and $J_{\rho}'$ the canonical lifts of $J$, $J_{0}$ and $J_{\rho}$ to $\nT_{vert,p}\totl{\hat{\ex B}}$. 
 
 Lemma \ref{average J'} states the following: let  $Z':=-J'\circ J'_{0}$, then define
\[W':=(1+Z')^{-1}(1-Z')\]
Then 
\[J'_{\rho}:=(1+\rho W')^{-1}(1-\rho W')J'_{0}\]

Because we chose $J_{M}$ to be equal to $J$ restricted to the image  of $\iota:T_{vert, p}\hat{\ex B}\longrightarrow \nT_{vert, p}\totl{\hat{\ex B}}$, $J'$ and $J'_{0}$ coincide on the image of $\iota$. The cokernel of $\iota$ has its canonical complex structure. Therefore  $(Z'-1)$ is a transformation $N$ which squares to $0$ because its image is contained in the image of $\iota$, and its  kernel contains the image of $\iota$. We may then write
\[W'=(2+N)^{-1}(-N)=-N/2\]
and 
\[\begin{split}J'_{\rho}&=(1-\rho N/2)^{-1}(1+\rho N/2)J'_{0}
\\&=(1+\rho N/2)^{2}J'_{0}
\\&=(1+\rho N )J'_{0}
\\&=(1-\rho)J'_{0}+\rho(1+N)J'_{0}
\\&=(1-\rho)J'_{0}+\rho J'\end{split}\]
 As $\omega$ is positive on $J'_{0}$ and $J'$ holomorphic planes in $\nT_{vert,p}\totl{\hat{\ex B}}$, it follows that $\omega$ is positive on $J'_{\rho}$ holomorphic planes in $\nT_{vert,p}\totl{\hat {\ex B}}$, and $J'_{\rho}$ is therefore tamed by $\omega$ in some neighborhood of $p$. As $J_{\rho}$ is equal to either $J$ or $J_{0}$ at other points within our stratum, it follows that $J_{\rho}$ is tamed by $\omega$ on some neighborhood of our stratum. This completes the proof of Claim \ref{J tame}.
 
 \
 
 As $\rho$ is $1$ in a neighborhood of $K$, $J_{\rho}$ coincides with $J$ in a neighborhood of $K$. It follows that $J_{\rho}$ gives  a $\omega$-tame,  $\dbar\log$ compatible  extension of $J$ to a neighborhood of $K$ and our stratum.

We may similarly extend $J$ to a neighborhood of $K$ and all strata of tropical dimension $m$, then continue with strata of smaller tropical dimension until $J$ has been extended to all of $\hat{\ex B}$  to be a $\dbar\log$ compatible almost complex structure tamed by $\omega$.

\stop

\

In the lemma below, we shall give an explicit method for constructing $\dbar\log$ compatible almost complex structures. In the lemma that follows it, we shall use this to show that if $J$ is a $\dbar\log$ compatible almost complex structure on $\ex B$, there locally exists a holomorphic embedding of the smooth part of $\ex B$ into some smooth almost complex manifold. This allows us to use standard results about $J$ holomorphic curves in smooth manifolds.

\

\begin{lemma}\label{local J construction} Suppose that 
$J_{0}$ is a $\dbar\log$ compatible almost complex structure on a $2n$-dimensional  coordinate chart.
Let \[\alpha_{1},\dotsc, \alpha_{n}\] be pointwise linearly independent smooth one-forms on this coordinate chart so that
 \[\{\alpha_{i},J_{0}\alpha_{i}\}\] form a basis for the cotangent space.
  Let \[\{v_{i},-J_{0}v_{i}\}\] be the dual basis for the tangent space. Then given any smooth one-forms \[\theta_{1},\dotsc,\theta_{n}\] generated by functions, let $X$ be the $n\times n$ matrix with $(i,j)$ entry \[X_{i,j}:=(J_{0}v_{i})(\theta_{j}+J_{0}\alpha_{j})\] Then on the region where $X$ is invertible, there exists a unique $\dbar\log$ compatible  almost complex structure $J$ so that 
\[J\alpha_{i}=\theta_{i}+J_{0}\alpha_{i}\] In particular, this almost complex structure $J$ exists where $\theta_{i}$ is small enough, and  $\log \dbar$ compatible almost complex structures $J$ which are close by to $J_{0}$ are determined by the one-forms $(J-J_{0})\alpha_{i}$. 
\end{lemma}

\

\pf

Given any choice of smooth one-forms $\beta_{i}$ so that $\{\alpha_{i},\beta_{i}\}$ form a basis for the cotangent space, there exists a unique smooth almost complex structure $J$ so that $J\alpha_{i}=\beta_{i}$. 
As $\{v_{i},-J_{0}v_{i}\}$ is a dual basis to $\{\alpha_{i},J_{0}\alpha_{i}\}$,  $\{\alpha_{i},\beta_{i}\}$ is a basis if and only if the matrix with $(i,j)$ entries $\beta_{j}(J_{0}v_{i})$ is invertible. Setting $\beta_{i}=\theta_{i}+J_{0}
\alpha_{i}$ implies that wherever the matrix $X$ with the above entries is invertible, there exists a unique smooth almost complex structure $J$ so that $J\alpha_{i}=\theta_{i}+J_{0}\alpha_{i}$. We must check that $J$ is $\dbar\log$ compatible.

 Let $\ex B$ indicate the open subset of our coordinate chart where $X$ is invertible. The one-forms $\alpha$ define a vector bundle inclusion 
\[\iota:\mathbb R^{n}\times \ex B\longrightarrow T^{*}\ex B\]
\[\iota(x,p):=\sum x_{i}\alpha_{i}(p)\]
and the vectors $v_{i}$ define a vector bundle projection 
\[\pi:T^{*}\ex B\longrightarrow \mathbb R^{n}\times \ex B\]
\[\pi(\beta,p):=(\beta(v_{1}),\dotsc,\beta(v_{n}),p)\]
so that:\\$\pi\circ\iota=\id$,\\ $\pi\circ J_{0}\circ\iota$ is the zero section,\\ $\iota\circ \pi$ is a projection, \\ $-J_{0}\circ\iota\circ\pi\circ J_{0}:=\id-\iota\circ\pi$ is also a projection,\\ and our matrix $X$ we are assuming to be invertible may be regarded as \[X:=\pi J_{0}J\iota\]
In what follows, we shall write  $J\circ J_{0} \circ \iota$ in terms of $J\circ \iota$, allowing us to write $J-J_{0}$ in terms of $(J-J_{0})\circ \iota$ which is determined by $\theta_{i}$.
\[\begin{split}JJ_{0}\iota&=JJ_{0}\iota X X^{-1}
\\&=J(J_{0}\iota\pi J_{0})J\iota X^{-1}
\\&=\lrb{J\iota\pi J\iota+\iota}X^{-1} 
\\&= \lrb{(J-J_{0})\iota\pi J\iota +J_{0}\iota\pi J\iota +\iota}X^{-1}
\\&= \lrb{(J-J_{0})\iota\pi J\iota +J_{0}\iota\pi (J-J_{0})\iota +\iota}X^{-1}
\\&= \lrb{(J-J_{0})\iota\pi J\iota -(J_{0}\iota\pi J_{0})J_{0}(J-J_{0})\iota +\iota}X^{-1}
\\&= \lrb{(J-J_{0})\iota\pi J\iota +J_{0}(J-J_{0})\iota -\iota\pi J_{0}(J-J_{0})\iota  +\iota}X^{-1}
\\&= \lrb{(J-J_{0})\iota\pi J\iota +J_{0}(J-J_{0})\iota -\iota\pi J_{0}J\iota -\iota +\iota}X^{-1}
\\&= \lrb{(J-J_{0})\iota\pi J\iota +J_{0}(J-J_{0})\iota}X^{-1}-\iota
\end{split}\]
Therefore, 
\[\begin{split}J-J_{0}&=(J-J_{0})\iota\pi -(J-J_{0})J_{0}\iota\pi J_{0}=(J-J_{0})\iota\pi-(JJ_{0}\iota+\iota)\pi J_{0} 
\\&=(J-J_{0})\iota\pi -\lrb{(J-J_{0})\iota\pi J\iota +J_{0}(J-J_{0})\iota}X^{-1}\pi J_{0}\end{split}\]
As $(J-J_{0})\iota$ applied to any smooth section of $\mathbb R^{n}\times \ex B$ is a smooth one-form on $\ex B$ which is generated by functions, and as noted in Remark \ref{preserves gbf}, $J_{0}$ of any smooth one-form generated by functions is another smooth one-form generated by functions, the above formula implies that $(J-J_{0})$ applied to any smooth one-form is a one-form which is generated by functions. Lemma \ref{J construction} then implies that $J$ is $\dbar\log$ compatible as required.

\stop

%
%
%
%

\begin{lemma}\label{J embedding} Let $J$ be family of $\dbar\log$ compatible almost complex structures on $\hat{\ex B}$. Then around every point $p$ in $\hat{\ex B}$, there exists some neighborhood  $U$ so that  there exists some smooth almost complex structure $J'$ on $\mathbb R^{2n}$ and some holomorphic map
\[f:(U,J)\longrightarrow(\mathbb R^{2n},J')\]
so that 
\[\totl f:\totl U\longrightarrow \mathbb R^{2n}\]
is an embedding in the sense that any smooth function on $U$ is $\totl f$ composed with  some smooth function defined on an open neighborhood of $ f(U)$  in  $\mathbb R^{2n}$.

\end{lemma}

\pf
By multiplying our family by a trivial factor of $\mathbb R$ if necessary, we may assume that $\hat{\ex B}$ is even dimensional.
The above lemma is local, so we may restrict to the case of a neighborhood of a point $p$ in a single coordinate chart $U=\mathbb C^{k}\times \et mP$. (As usual we shall assume that $p$ is in the stratum of this coordinate chart corresponding to the interior of $P$). By using a change of coordinates, we may also assume that at the point $p$, our almost complex structure is equal to the standard complex structure on $\mathbb C^{k}\times \et mP$ restricted to the vertical tangent space. 

We may consider $U=\mathbb C^{k}\times \et mP$ as being a subset of $V:=\mathbb C^{k}\times \ex T^{n'}\times \lrb{\et 11}^{n} $ cut out by some monomial equations and $\totl U$ as being a subset of $\totl V=\mathbb  C^{k+n}$ cut out by related monomial equations. 
Extend $J$ to a $\dbar\log$ compatible almost complex structure on neighborhood of $p$ in $V$ as follows: 

Lemma \ref{local J construction} implies that we may construct a $\dbar\log$ compatible $J''$ on some neighborhood of $p$ in $V $ by specifying some one-forms generated by functions to be equal to $(J''-J_{0})$ composed with the imaginary part of $dz_{j}$ and  $\tilde z_{i}^{-1}d\tilde z_{i}$. (Here $J_{0}$ indicates the standard almost complex structure, the $z_{j}$ are coordinates on $\mathbb C^{k}$ and $\tilde z_{i}$ are coordinates on $\ex T^{n'}\times (\et11)^{n}$.)

 For each coordinate $\tilde z_{i}$, choose  the corresponding one-form to be $0$ restricted to  $T_{p}V $, and to restrict to the vertical tangent space of  $U\subset V$ to be $(J-J_{0})$ applied to the restriction of the imaginary part of $\tilde z_{i}^{-1}d\tilde z_{i}$ to the vertical tangent space of $U$. This is possible because of the following observations:
 \begin{itemize} \item
 $(J-J_{0})$ applied to the restriction of the imaginary part of  $\tilde z_{i}^{-1}d\tilde z_{i}$ to the vertical tangent space of $U$ is a family of one-forms generated by functions, which vanishes on $T_{vert,p}U$ because $J=J_{0}$ on $T_{vert,p}U$.
 \item Any family of one-forms generated by functions on $U$ extends by definition to a one-form on $U$ which is generated by functions. As our family of one-forms vanishes on $T_{vert,p}U$, we may extend it to a one-form which is generated by functions and vanishes on $T_{p}U$. 
 
 \item Any smooth function on $U$ extends by definition to a smooth function on $V$. As one-forms generated by functions are constructed from smooth functions by the operations of exterior differentiation and multiplication, it follows that any one-form on $U$ which is generated by functions extends to a one-form which is generated by functions on $V$. As our one-form vanishes on $T_{p}U$, the extension may be chosen to vanish on $T_{p}V$. 
 \end{itemize}
  Similarly, for each coordinate $ z_{j}$, choose  the corresponding one-form to be $0$ restricted to  $T_{p}V$, and to restrict to the vertical tangent space of  $U$ to be $(J-J_{0})$ applied to the restriction of the imaginary part of $ dz_{j}$ to the vertical tangent space of $U$.

The matrix $X$ from Lemma \ref{local J construction} is the identity at $p$, therefore the almost complex structure $J''$ constructed this way is defined in a neighborhood of $p$ in $V $ and is $\dbar\log$ compatible. By restricting to a neighborhood and rescaling, we can assume that $J''$ is defined on all of $V$.

To summarize, we now have for any $\dbar\log$ compatible family of almost complex structures $J$ on $\hat{\ex B}$ and a point $p\in \hat{\ex B}$,  there exists an open neighborhood $U$ of $p$, a $\dbar\log$ compatible almost complex structure $J''$ on some $V:=\mathbb C^{k}\times \ex T^{n'}\times \lrb{\et 11}^{n} $ and a fiberwise holomorphic embedding
\[g: (U,J)\longrightarrow (\mathbb C^{k}\times \ex T^{n'}\times \lrb{\et 11}^{n},J'') \]

The smooth part of $\mathbb C^{k}\times \ex T^{n'}\times \lrb{\et 11}^{n} $ is $\mathbb C^{k+n}$. A one-form generated by functions on $\mathbb C^{k}\times \ex T^{n'}\times \lrb{\et 11}^{n} $ is equivalent to a smooth one-form on $\mathbb C^{k+n}$. The $\dbar\log$ compatible almost complex structure $J''$ sends one-forms generated by functions to one-forms generated by functions, so it gives an almost complex structure $J'$ on $\mathbb C^{k+n}$. The composition of $g$ with the smooth part map is the required
holomorphic map 
\[f:(U,J)\longrightarrow (\mathbb C^{k+n},J')\]
As the smooth part of $f$ is equal to the smooth part of $g$, it is an embedding.

\stop

\section{Local area bounds from energy bounds}\label{local area bound section}

In this section, we prove that in some cases, a bound on the energy of holomorphic curves in $\ex B$ 
implies a local area bound for those holomorphic curves, and in other cases, we prove that an energy bound and a specification of the number and type of punctures for a holomorphic curve implies a local area bound.

\

To understand the definition of a local area bound below, note that
there is a functor $\mathcal F$ which sends each exploded manifold $\ex B$ to a disjoint union of smooth manifolds 
\[ \mathcal F(\ex B):=\coprod_{p\in \totb {\ex B}}\ex B_{p} \]
so that every map of a connected smooth manifold $M$ to $\ex B$ is equivalent to a smooth map to $\ex B_{p}$, the subset of  $\ex B$ with tropical part $p\in\totb{\ex B}$. For example, $\mathcal F(\et mP)=\coprod_{p\in P}(\mathbb C^{*})^{n}$. Any smooth metric on $\ex B$ gives a smooth metric on $\ex B_{p}$ for all $p$, so we can define all the usual concepts from differential geometry.

\begin{defn}\label{lab}Given an exploded manifold $\ex B$ with a metric $g$, a curve $f:\ex C\longrightarrow \ex B$ has local area bounded by $c$ if for every point $p\in \ex B$, the area of $\ex C$ contained within a distance $1$ of $p$ is  bounded by $c$. 

\

More formally, regard $\ex C$ as a disjoint union of manifolds by applying the functor $\mathcal F$. The metric $g$ defines an area form on $\ex B$ which pulls back to give a density on $\mathcal F(\ex C)$. The set of points  $q\in\mathcal F(\ex C)$ so that $f(q)$ is within distance $1$ from $p$ using the metric $g$ is a subset of $\mathcal F(\ex C)$. The area of $\ex C$ contained within distance $1$ of $p$ is defined as the integral of this density over this subset of $\mathcal F(\ex C)$.  This integral is well defined if it is bounded because the contribution of any subset is non negative.

\end{defn}

The following lemma gives a sufficient condition for the local area of holomorphic curves to be bounded in terms of topological data. Although it is not the sharpest condition, it is presented here because it is easy to state and covers most cases of interest. The bound on local area from the lemma below is in terms of the $\omega$-energy of the curve $f$ and some combinatorial data involving the external edges of the tropical curve $\totb{f}$.  To understand what is meant by external edges, recall that the domain $\totb{\ex C}$ of $\totb{f}$ is a connected complete integral affine graph. If $\totb{\ex C}$ is not equal to $\mathbb R$, then the external edges of $\totb{f}$ are those edges isomorphic to $[0,\infty)$; these edges have a canonical outgoing unit vector. In  the exceptional case that $\totb{\ex C}=\mathbb R$, we shall say that $\mathbb R$ has two external edges with outgoing unit vectors the positive and negative unit vectors in $\mathbb R$.

\begin{lemma}\label{immersion bound}
Suppose that a complete exploded manifold $\ex B$ has a metric $g$ and a $\dbar\log$ compatible almost complex structure $J$ tamed by a taming form $\omega$. Suppose further that $\ex B$ is basic, and that there exists an  integral affine map $\totb{\ex B}\longrightarrow \mathbb R^{n}$ with injective derivative (restricted to each stratum). 

Given any holomorphic curve $\ex C\longrightarrow \ex B$, composition of the tropical part of the curve with $\totb{\ex B}\longrightarrow \mathbb R^{n}$ gives a map $\totb{\ex C}\longrightarrow \mathbb R^{n}$. The derivative of this map applied to an outgoing unit vector on an external edge is  $(m_{1},\dotsc,m_{n})$. Let $E'$ be the sum of all positive $m_{i}$ from each external edge. Let $E$ denote the $\omega$-energy of the curve.
 Then  there exists a constant $c$ so that the local area of any holomorphic curve is bounded by  $c(E+E')$.

\

A similar statement holds for families of exploded manifolds. In particular, 
 suppose that a family of basic exploded manifolds $\hat {\ex B}$ with a complete metric $g$ has a family of $\dbar\log$ compatible almost complex structures $J$  tamed by a family of taming forms $\omega$, and that there is an integral affine map $\totb{\hat{\ex B}}\longrightarrow \mathbb R^{n}$ which has injective derivative when restricted to the strata of each exploded manifold $\ex B$ in the family.
 
  Then there exists a continuous  function $c:\hat{\ex B}\longrightarrow (0,\infty)$ so that given any holomorphic curve in $\hat{\ex B}$,  the local area of that curve within a ball of radius $1$ around $p$  is bounded by $c(p)(E+E')$.
   
\end{lemma}

\begin{remark}
Note that in the special case that $\totb{\ex B}$ may be immersed in the negative quadrant of $\mathbb R^{n}$, $E'$ is always $0$, so the local area of any holomorphic curve is bounded by $cE$. In particular, this holds if $\ex B$ is basic and all strata of $\totb{\ex B}$ are diffeomorphic to some $[0,\infty)^{k}$, as we may embed each stratum equal to $[0,\infty)$ in a different coordinate direction. If $\ex B$ is obtained by applying the explosion functor to a log scheme $M$, then this condition of an immersion $\totb{\ex B}\longrightarrow\mathbb [0,\infty)^{n}$ is equivalent to $M$ being almost generated in the language of Gross and Siebert from \cite{GSlogGW}, and if there is an immersion $\totb{\ex B}\longrightarrow \mathbb R^{N}$, then $M$ is quasi generated. 

 \end{remark}

\pf

We shall first prove the case of a single exploded manifold $\ex B$.  To obtain a local area bound around a point $p\in\ex B$, we shall make a modification of $\omega$ to a  two-form $\omega_{p}$ on a refinement of $\ex B$. This two-form $\omega_{p}$ shall be nonnegative on holomorphic planes so that in the ball of radius $1$ around $p$, the area of holomorphic curves is bounded by some (uniform) constant times the integral of $\omega_{p}$, and so that the integral of $\omega_{p}$ is bounded by $E+E'$. 

The structure of this proof is as follows: We shall first construct $\omega_{p}$, leaving some variables in the construction to be determined later. Then we shall verify that the integral of $\omega_{p}$ is bounded by $E+E'$. Finally, we shall ensure that $\omega_{p}$ is non negative on holomorphic planes, and is sufficiently positive on holomorphic planes in the ball of radius $1$ around $p$ to give a bound on the local area of holomorphic curves.

\

Recall that we have an integral affine map $\totb{\ex B}\longrightarrow \mathbb R^{n}$ which is injective restricted to each stratum of $\totb{\ex B}$. For $j=1,\dotsc, n$, construct a map 
\[r_{j}:\ex B\longrightarrow (0,\infty)\e {\mathbb R}\subset \mathbb R\e{\mathbb R}\] with the property that in a neighborhood of any point on $\ex B$,  \[r_{j}={\abs f}\] for some  exploded function $f$ locally defined on $\ex B$, with tropical part the pullback of  the $j$th coordinate function on $\mathbb R^{n}$ under the map $\totb{\ex B}\longrightarrow \mathbb R^{n}$. (The absolute value of $c\e{x}$ is $\abs c\e{x}$.) Such a map locally exists, and if $r_{j}$ and $r_{j}'$ are two such maps, they may be averaged using any smooth $\mathbb R$-valued function $\rho$ to produce another such map $r_{j}^{1-\rho}r_{j}'^{\rho}=r_{j}(r_{j}^{-1}r_{j}')^{\rho}$. (Note that this averaged map is in the required form because $(r_{j}^{-1}r_{j}')^{\rho}$ is a smooth $(0,\infty)$-valued function which may also be regarded as an exploded function.) We may therefore patch together locally defined $r_{j}$'s using a partition of unity. 

 There is a smooth one-form $r_{j}^{-1}dr_{j}$ on $\ex B$ which is  locally equal to the real part of $f^{-1}df$.

\

Define a one-form $\alpha_{j}$ on $\ex B$ by
\[\alpha_{j}:=r_{j}^{-1}dr_{j}\circ J\]

The assumption that $J$ is $\dbar\log$ compatible implies that $d\alpha_{j}$ is a two-form  which is generated by functions. In particular, locally, $\alpha_{j}=d\log\abs f\circ J$, so $d\alpha$ is locally equal to $d$ of the imaginary part of $\dbar\log f$, which is generated by functions because $J$ is $\dbar\log$ compatible.

\

There is a refinement $\ex B'$ of $\ex B$ given by subdividing each polytope in $\totb{\ex B}$ by its intersection with planes in the coordinate directions of $\mathbb R^{n}$ passing through the image of $p$. Our modified form $\omega_{p}$ shall be a smooth differential form on $\ex B'$. We shall modify $\omega$ by $\sum_{j}d(h_{j}\alpha_{j})$, where the functions $h_{j}$ are smooth functions on $\ex B'$ (and obey further conditions to be enumerated later).

Let $\omega_{p}$ be the differential form on $\ex B$  given by
\[\omega_{p}:=\omega+\sum_{j}d(h_{j}\alpha_{j})=\omega+\sum_{j}h_{j}d\alpha_{j}+\sum_{j}dh_{j}\wedge (r^{-1}dr_{j}\circ J)\]

 To ensure  that $\omega_{p}$ is still nonnegative on holomorphic planes, we shall  require that $\abs {h_j}$ is small, and also choose $h_{j}$ so that $dh_{j}$ is some nonpositive function times $r_{j}^{-1}dr_{j}$; so $h_{j}$ is a nonincreasing function of $r_{j}$. Choose $h_{j}$ so that it is some constant less than $(2\pi)^{-1}$ when $r_{j}$ is small enough, and $0$ when $r_{j}$ is large enough (recall that $x\e {a}<y\e b$ if $a>b$ or $a=b$ and $x<y$, as $\e1$ is thought of as being infinitesimal). 
 
 \
 
 We shall now verify that the integral of $\omega_{p}$ is bounded by $E+E'$. Theorem 3.4 from \cite{dre}, which is a version of Stokes' theorem, allows us to compute the integral of $f^{*}d(h_{j}\alpha_{j})$ over $\ex C$. This theorem states that the integral of $d$ of a form in $\Omega^{*}_{c}$ over an exploded manifold with boundary is equal to the integral of that form over the boundary of the exploded manifold.  To be in $\Omega^{*}_{c}(\ex C)$, a form must have complete support, vanish on all integral vectors, and also vanish on all external edges of $\ex C$. Our form $h_{j}\alpha_{j}$ is a form on some refinement of $\ex C$ which vanishes on all integral vectors, but it does not vanish on those external edges on which $h_{j}$ is nonzero. To remedy this, we apply Stokes' theorem to the exploded manifold with boundary $\{f^{*}r_{j}\geq 1\e {a_{j}}\}\subset \ex C$ for some sufficiently large constant $a_{j}$ so that  when $r_{j}\leq 1\e a_{j}$, $h_{j}$ is constant, and $\ex C$ consists only of external edges. Then our version of Stokes' theorem tells us 
\[\int_{\{f^{*}r_{j}\geq 1\e {a_{j}}\}\subset \ex C}f^{*}d(h_{j}\alpha_{j})=\int_{\{f^{*}r_{j}=1\e{a_{j}}\}}h_{j}\alpha_{j}\]
  As the right hand side of the above equation is constant for $a_{j}$ large enough, it must compute the integral of $f^{*}d(h_{j}\alpha_{j})$ over all of $\ex C$. As $h_{j}$ is chosen to be smaller than $1/2\pi$, this right hand side is less than the integral of $\alpha_{j}/2\pi$ around the boundary of $\{f^{*}r_{j}\geq 1\e {a_{j}}\}\subset \ex C$. Each external edge in the direction $(m_{1},\dotsc,m_{n})$ contributes $m_{j}$ to this integral if $m_{j}$ is positive, and contributes nothing if $m_{j}$ is not positive, because for $a_{j}$ chosen large enough, such an edge does not intersect the region where $f^{*}r_{j}=1\e a_{j}$. Therefore, 
  \begin{equation}\label{omegap}\int_{\ex C}f^{*}\omega_{p}\leq E+E'\end{equation}
 where $E$ is the integral of $\omega$ and $E'$ is equal to the sum of all positive $m_{j}$ which appear in the derivatives  $(m_{1},\dotsc,m_{n})$ of external edges of $\totb f$ with respect to outgoing unit vectors.

\

 We shall now specify conditions on $h_{j}$ so that $\omega_{p}$ remains non negative on holomorphic planes.   Recall that $\omega$ is defined as a two-form on $\totl{\ex B}$.  The two-form $d\alpha_{j}$ is generated by functions, so it lifts (possibly non-uniquely) to a two-form on $\totl{\ex B}$, which we shall again simply call $d\alpha_{j}$. Then $d\alpha_{j}(v,Jv)$ is a continuous function of $v\in \nT\totl{\ex B}$. The set of nice tangent vectors $v$ on $\totl{\ex B}$ for which $\omega(v,Jv)=1$ is compact, so there exists some $M$ so that $\sum_{j}\abs{d\alpha_{j}(v,Jv)}\leq M$ for any $v$ so that $\omega(v,Jv)=1$. Then choose
\[\abs{h_{j}}<\frac 1{2Mn}\] 
so
\[\abs{\sum_{j}h_{j}d\alpha_{j}(v,Jv)}\leq \frac 12\omega(v,Jv)\]
The above inequality holds both for any $v$ in the nice tangent space to $\totl{\ex B}$ and any $v$ in $T\ex B$. In particular, for all $v\in T\ex B$, 
\begin{equation}\label{omegap pre}\begin{split}\omega_{p}(v,Jv)\geq &\frac 12 \omega(v,Jv)+\sum_{j}(dh_{j}\wedge \alpha_{j})(v,Jv)\\&=\frac 12\omega(v,Jv)+\sum_{j}-h'_{j}( (r_{j}^{-1}dr_{j}(v))^{2}+(r_{j}^{-1}dr_{j}(Jv))^{2})\end{split}\end{equation}
where $h'_{j}$ is the smooth function on $\ex B'$ defined by
\[dh_{j}=h'_{j}r_{j}^{-1}dr_{j}\]
Note that as we have chosen $h_{j}$ to be non-increasing as a function of $r_{j}$, inequality (\ref{omegap pre}) above implies that $\omega_{p}$ is non-negative on holomorphic planes.

\

Now, we must use inequality (\ref{omegap pre}) 
to get a uniform lower bound for $\omega_{p}(v,Jv)$ in a neighborhood of $p$. To achieve this, we need to construct $h_{j}$ so that $-h_{j}'$ is uniformly large in a neighborhood of $p$.  We may construct our $h_{j}$ uniformly for any $p$ as a function of $r_{j}(\cdot)r_{j}^{-1}(p)$. Our only constraints on $h_{j}$ are that it is monotone decreasing, positive and bounded by $1/2Mn$ and $1/2\pi$.  How large we can get $-h_{j}'$ in a region is inversely proportional to the amount that $\log\lrb{r_{j}(\cdot)r^{-1}_{j}(p)}$ varies in that region, because a change of variables which is linear in $\log\lrb{r_{j}(\cdot)r^{-1}_{j}(p)}$ scales these two quantities inversely. 

In particular, choose some non-increasing function $\rho:\mathbb R\longrightarrow [0,1/2Mn)$ so that $\rho(x)=0$ for $x>2$,  $\rho(x)$ is a constant $c$ for $x<-2$, and  $\rho'(x)>\epsilon>0$ for all $x\in[-1,1]$. Then for any  $p$ in $\ex B$, we may define
\[h_{j}(p'):=\begin{cases}c\hspace{4 cm} \text{if }\totl{r_{j}(p')/r_{j}(p)}=0\\
 \rho(L^{-1}\log(r_{j}(p')/r_{j}(p)))\ \ \ \   \text{ if }r_{j}(p')/r_{j}(p)\in(0,\infty)
 \\ 0 \hspace{4 cm}\text{if }\totl{r_{j}(p)/r_{j}(p')}=0
  \end{cases}\]
If  $\omega_{p}$ is defined using such  functions $h_{j}$, then
\begin{equation}\label{opest}\omega_{p}(v,Jv)\geq \frac 12\omega(v,Jv)+ \frac \epsilon{L}\sum_{j}((r^{-1}_{j}dr_{j}(v))^{2}+(r^{-1}_{j}dr_{j}(Jv))^{2})\end{equation}
at all points  $p'$ so that
\[\abs{\log\lrb{r_{j}(p')r^{-1}_{j}(p)}}<L\text{ for all }j\]

As $\ex B$ is compact and $r_{j}^{-1}dr_{j}$ is a smooth one-form, the size of $r_{j}^{-1}dr_{j}$ is bounded in whatever metric we have chosen. Therefore, by choosing $L$ large enough (independent of $p$), the inequality (\ref{opest}) holds in a ball of radius $1$ around $p$.

It remains to prove that there exists some constant $c$ so that for any vector $v$ in $T\ex B$,  
\begin{equation}\label{rts}\abs v^{2}\leq c\lrb{\frac 12\omega(v,Jv)+ \frac \epsilon{L}\sum_{j}((r^{-1}_{j}dr_{j}(v))^{2}+(r_{j}^{-1}dr_{j}(Jv))^{2})}\end{equation}

The right hand side of the above equation is equal to $c\hat\omega(v,Jv)$, where 
\[\hat \omega:=\frac 12\omega-\frac \epsilon L\sum_{j}r^{-1}_{j}dr_{j}\wedge\alpha_{j}\]
 As $\hat \omega$ is a smooth two-form on the  compact exploded manifold $\ex B$, to prove the inequality (\ref{rts}), it suffices to prove that this two-form is positive on holomorphic planes. For any $v\in T_{p}\ex B$ which is nonzero as a vector in $\nT_{p}\totl{\ex B}$, $\omega(v,Jv)>0$, so $\hat \omega(v,Jv)>0$. It  therefore suffices to check positivity for nonzero vectors $v$ in $T\ex B$ for which $dg(v)=0$ for all smooth functions $g$.

Restrict attention to a coordinate chart. Within our coordinate chart, $r_{j}=\abs{g\tilde z^{\alpha}}$ so $r_{j}^{-1}dr_{j}$ is the real part of $\tilde z^{-\alpha}d\tilde z^{\alpha}$ plus a one-form which is generated by functions, and as $J$ is $\dbar\log$ compatible, $r_{j}^{-1}dr_{j}\circ J$ is the imaginary part of $\tilde z^{-\alpha}d\tilde z^{\alpha}$ plus a one-form which is generated by functions. The choice of $r_{j}$ as coming from a map $\totb{\ex B}\longrightarrow \mathbb R^{n}$ which is injective restricted to strata of $\totb{\ex B}$ implies that the real  and imaginary parts of the corresponding $\tilde z^{\alpha}\frac\partial {\partial \tilde z^{\alpha}}$ generate the space of vectors which project to zero on $\totl{\ex B}$, so on any of these vectors $v\neq0$, $\sum_{j}(r_{j}^{-1}dr_{j}(v)^{2}+r_{j}^{-1}dr_{j}(Jv)^{2})>0$. Therefore, $\hat \omega$ is positive on all holomorphic planes. The compactness of $\ex B$ then implies that there exists a constant $c$ so that the area of any piece of a holomorphic curve is bounded by $c$ times the $\hat \omega$ area. The required inequality, (\ref{rts}) follows.

 As the area of a parallelogram spanned by $v$ and $Jv$ is bounded by the maximum of $\abs v^{2}$ and $\abs {Jv}^{2}$, the inequalities (\ref{rts}) and (\ref{opest}) imply that the area form on the unit ball around $p$ is bounded on holomorphic planes by $c\omega_{p}$. As $\omega_{p}$ is non-negative on holomorphic planes elsewhere, it follows that the local area of any holomorphic curve within the unit ball around $p$ is bounded by the integral of $c\omega_{p}$, which, as noted in inequality (\ref{omegap}), is bounded by $c(E+E')$.

\

With the case of a single exploded manifold understood, the proof of the family case is an easy generalization. 
We may now define $r_{j}$ and $r_{j}^{-1}dr_{j}$ exactly as before. When defining $\alpha_{j}:=dr_{j}\circ J$, note that $\alpha_{j}$ is only defined as a smooth family of differential forms, because $J$ acts on $T^{*}_{vert}\hat {\ex B}$, and not $T^{*}\hat {\ex B}$. Nevertheless, as before, $d\alpha_{j}$ is a smooth family of two-forms generated by functions. Similarly, $d (h_{j}\alpha)$ should be interpreted as taking the differential in every fiber of the family, obtaining a smooth family of two-forms.  

Follow the rest of the argument after restricting to a compactly contained subfamily so that the required  bounds on $h_{j}$ and $L$ are possible to achieve. This implies that given any compact subset of our family, there exists a constant $c$ so that the area within a unit ball of a holomorphic curve is bounded by $c(E+E')$. Noting that we have chosen our metric on $\hat {\ex B}$ complete, and allowing our constant $c$ to depend continuously on position in the base of the family gives the required result.

\stop  

\

The following observation gives a kind of conservation of momentum condition for the tropical part of curves in exploded manifolds. One use of this conservation of momentum condition is to prove local area bounds from topological bounds on holomorphic curves in some cases when $\totb{\ex B}$ does not immerse into some $\mathbb R^{N}$. 

\begin{lemma}\label{cm} Let $ \ex B$ be a basic exploded manifold  with an integral affine map \[A:\totb{\ex B}\longrightarrow \mathbb R\] There exists a smooth differential form $\alpha$ on $\ex B$ locally equal to a one-form generated by functions plus the imaginary part of $g^{-1}dg$ for some exploded function $g$ with tropical part $\totb g=A$.

Given any such one-form $\alpha$ and any  map $f:\ex C\longrightarrow \ex B$ from a compact two dimensional exploded manifold $\ex C$, define an external edge of $\totb{\ex C}$ to be any edge of $\totb{\ex C}$ with closure equal to $[0,a)$ for some $0<a\leq \infty$. Define the momentum of an external edge to be the derivative of $A\circ\totb f$.

Then the sum of the momentum of all external edges of $\ex C$ is equal to 
\[-\frac 1{2\pi}\int_{\ex C}f^{*}d\alpha\] 

\end{lemma}

\pf

A one-form satisfying the conditions enumerated for $\alpha$ locally exists. Given any two such one-forms $\alpha_{1}$ and $\alpha_{2}$,  the one-form $\alpha_{1}-\alpha_{2}$ is generated by functions because if $h$ and $g$ are locally defined exploded functions with tropical part equal to $A$, then $h/g$ is a smooth function, so the imaginary part of $h^{-1}dh-g^{-1}dg$ is a smooth one-form generated by functions. It follows that if $\alpha_{1}$ and $\alpha_{2}$ are two one-forms satisfying the requirements of $\alpha$, then so is $\rho\alpha_{1}+(1-\rho)\alpha_{2}$. We may therefore use a partition of unity to patch together locally defined candidates for $\alpha$ into a globally defined $\alpha$.

Note that $d\alpha$ is a two-form which is generated by functions on $\ex B$, so $f^{*}d\alpha$ is generated by functions on $\ex C$ and therefore the integral of $f^{*}d\alpha$ is well defined and equal to the integral of $f^{*}d\alpha$ over the smooth components of $\ex C$. For computing this integral, we may therefore replace $\ex C$ with a complete exploded manifold with boundary by adding a boundary component to each external edge. In particular, if an external edge is covered by a coordinate $\tilde z\in \et 1{[0,a)}$, then we may replace that coordinate chart by the subset of $\et 1{[0,b]}$ where $\tilde z\geq 1\e b$ for some $0<b<a$. On the region where $\totl{\tilde z}=0$ in such a coordinate chart, $f^{*}\alpha$ is equal to $m$ times the imaginary part of $\tilde z^{-1}d\tilde z$, where $m$ is the momentum of that external edge, so the integral of $f^{*}\alpha$ over that boundary component  is equal to $-2\pi m$. As our new surface with boundary is complete, the version of Stokes' theorem proved in \cite{dre} applies, and the integral of $f^{*}d\alpha$ over $\ex C$ is equal to the integral of $\alpha$ over the boundary, which is $-2\pi$ times the sum of the momentum of the external edges of $\ex C$. 

\stop

\

 \section{Estimates}\label{estimates section}
 
 In this section, we shall prove the analytic estimates required for our compactness theorem. These estimates do not differ much from the standard estimates used to study (pseudo)holomorphic curves in smooth symplectic manifolds. They are also proved using roughly the same methods as in the smooth case. We shall make the following assumptions:
 
  \begin{enumerate}
  \item We'll assume that our target $\ex B$ is complete. This is required to give us the kind of bounded geometry found in the smooth case when the target is compact or has a cylindrical end. We shall also assume that $\ex B$ is basic. The assumption that $\ex B$ is basic  is mainly for convenience in the arguments within this section. It is most useful in establishing  coordinates and notation for the construction of families of smooth curves in the proof of Theorem \ref{completeness theorem}.
 \item We'll assume our almost complex structure $J$ on $\ex B$ is $\dbar\log$ compatible. This assumption allows us to use Lemma \ref{J embedding} to apply standard analysis of $J$-holomorphic curves in manifolds to study the smooth part of our curves in $\ex B$.  In section \ref{local area bound section}, this assumption was also used to get  local area bounds for holomorphic curves.
 \item We'll assume there is a taming form $\omega$ which tames $J$.
\item We'll assume a local area bound for our curves in the sense of definition \ref{lab} on page \pageref{lab}.
 \end{enumerate}

 The strategy of proof for most of our estimates is roughly as follows:  first prove a weak estimate in a metric on $\totl{\ex B}$, and then improve this to get an estimate in an actual metric on $\ex B$. We shall use this strategy to prove estimates for the the derivative at the center of  holomorphic disks with bounded local area and small $\omega$-energy, and to get strong estimates on the behavior of holomorphic cylinders with bounded local area and small $\omega$-energy. 

 \
 
 This section ends with Proposition \ref{decomposition proposition}, which is a careful statement of the fact that any holomorphic curve with bounded local area, $\omega$-energy, genus, and number of punctures can be decomposed into a bounded number of components, which either have  `bounded conformal structure  and bounded derivative', or are annuli with small $\omega$-energy, (and so we have strong estimates on their behavior.) 
 
 \
 
 Throughout this section, we shall use the notation $(\ex B, J,\omega,g,\A)$ to indicate a smooth, basic, complete exploded manifold $\ex B$, a $\dbar\log$ compatible almost complex structure $J$ tamed by a taming form $\omega$, a metric $g$, and a local area bound $\A$.
 
  If we say that a constant depends continuously on $(\ex B,J,\omega,g,\A)$, we mean that if we have a family $\hat{\ex B}\longrightarrow\ex G$ with such a structure on each fiber, the constant can be chosen a continuous function on $\ex G$. 
  
 \

 For a lot of our analysis, we shall be thinking of $\ex B$ as a disjoint union $\coprod_{p\in \totb{\ex B}}\ex B_{p}$ of smooth manifolds, where $\ex B_{p}$ is the smooth manifold with the property that any map of a smooth manifold to $\ex B$ with image $p\in\totb{\ex B}$ is equivalent to a smooth map to $\ex B_{p}$. Thinking this way of $\ex B$ as a disjoint union of smooth manifolds allows us to use methods and intuitions familiar from differential geometry, however as $\coprod_{p\in \totb{\ex B}}\ex B_{p}$ is usually not compact, we must bound its geometry. 
 
 Recall that the topology on $\ex B$ is the pullback of the topology from $\totl{\ex B}$--- this is usually a much coarser topology than the topology on $\coprod_{p\in\totb{\ex B}}\ex B_{p}$. Given any metric $g$ on some complete $\ex B$, the injectivity radius  is a continuous function in this coarser standard topology on $\ex B$, and is therefore bounded below.

\begin{lemma}\label{bounded geometry}
 Given a complete exploded manifold $\ex B$ with a metric $g$, there exists some radius $R>0$ smaller than the injectivity radius of $(\ex B,g)$ so that the following holds:

 Given a sequence of points $p_{n}\in\ex B$ converging to $p$, consider the ball of radius $R$ around $p_{n}$ as a smooth manifold $B_{R}(p_{n})$. \begin{description}\item[]There exists a sequence of diffeomorphisms
 
 \[f_{n}:B_{R}(p_{n})\longrightarrow B_{R}(p)\]
 
 \item[] so that\begin{description}\item[] $(f_{n})_{*}g$ converges to $g$ in $C^{\infty}$,\item[] and given any smooth tensor field $\theta$ on $\ex B$ (such as an almost complex structure), $(f_{n})_{*}(\theta)$ converges to $\theta$ in $C^\infty$.
 \end{description}
 \
 
 Similarly, given a family $\hat{\ex B}\longrightarrow \ex G$ of complete exploded manifolds with fiberwise metrics $g$, we may choose $R>0$ to depend continuously on $\ex G$ so that the following holds: 

 Given a sequence of points $p_{n}\in\hat {\ex B}$ converging to $p$, consider the ball of radius $R$ around $p_{n}$ in the fiber containing $p_{n}$ as a smooth manifold $B_{R}(p_{n})$. Then the above conclusion holds: There is a sequence of diffeomorphisms $f_{n}:B_{R}(p_{n})\longrightarrow B_{R}(p)$
 so that $(f_{n})_{*}g$ converges to $g$ in $C^{\infty}$, and given any smooth vertical tensor field $\theta$ on $\hat {\ex B}$, $(f_{n})_{*}(\theta)$ converges to $\theta$ in $C^\infty$.
 
  \end{description}
\end{lemma}

\pf
By choosing $R$ small enough, we may assume that $B_{R}(p)$ is contained in a single coordinate chart, and $B_{R}(p_{n})$ is eventually contained in the same coordinate chart.
 The lemma then follows from the fact that all smooth tensor fields (and their derivatives) can be locally written as some sum of smooth functions on $\ex B$ times some basis tensor fields.
  
\stop

\

The above lemma tells us that the local geometry of a complete \exploded manifold is bounded in the same way that the geometry of a compact manifold is bounded. 
 
\

\
 
 Note that $\totl{\ex B}$ has a smooth structure which consists of the sheaf of smooth functions on the topological space $\totl{\ex B}$. As noted in Section \ref{dbar compatible section}, a $\dbar\log$ almost complex structure on $\ex B$ lifts to an almost complex structure $J'$ on the nice tangent space to $\totl{\ex B}$ so it makes sense to talk of holomorphic curves in $\totl{\ex B}$ as $J'$-holomorphic smooth maps of nodal Riemann surfaces to $\totl{\ex B}$. (Note that Lemma \ref{nice image} implies that the derivative of such a smooth map must have image contained in the nice tangent space to $\totl{\ex B}$, so being $J'$-holomorphic is well defined.)
 
 Lemma \ref{J embedding}  states that if $J$ is a $\dbar\log$ compatible almost complex structure on $\ex B$, then $\ex B$ is covered by open subsets $U$ with holomorphic maps $\iota:(U,J)\longrightarrow (\mathbb R^{N},J')$ so that $\totl\iota:\totl U\longrightarrow \mathbb R^{N}$ is an embedding. The removable singularity theorem for $J'$-holomorphic curves (first proved in example 1.4.B of \cite{gromov}) implies that the smooth part of any holomorphic curve in $\ex B$ is a smooth holomorphic curve in $\totl{\ex B}$, which is locally equal to a $J'$-holomorphic curve in $\mathbb R^{N}$ which lies in the image of $\iota$. 

\

\begin{defn}Given $(\ex B,J,\omega)$, let $J'$ denote the induced almost complex structure on the nice tangent space to $\totl{\ex B}$, and define  
\[\langle v,w\rangle_\omega:=\frac 12(\omega(v,J'w)+\omega(w,J'v))\]
Note that $\langle\cdot,\cdot\rangle_{\omega}$ defines a metric on $\nT_{p}\totl{\ex B}$, but pulls back to a pseudometric on $T_{p}\ex B$ as it vanishes  on any vector $v$ which acts as the zero derivation on smooth functions. 

Given any smooth map from a manifold $f:M\longrightarrow \totl{\ex B}$, Lemma \ref{nice image} implies that $Tf$ has image contained in the nice tangent space to $\totl{\ex B}$, so \[\langle df(\cdot),df(\cdot)\rangle_{\omega}\] defines a smooth, symmetric, positive semidefinite form on $M$, so we can define the length of smooth paths in $\totl {\ex B}$ or the area of holomorphic curves in $\totl{\ex B}$ as usual.

Define the distance function $d_{\omega}$ on $\totl{\ex B}$ as follows: Let $d_{\omega}(p,q)$ be the infimum of the length of all smooth paths in $\totl{\ex B}$ joining $p$ with $q$. 
\end{defn}

\begin{remark}\label{metric embedding}As pointed out in Lemma \ref{J embedding}, $\totl{\ex B}$ locally embeds as a closed holomorphic subset of some $(\mathbb R^{N},J')$. The metric $\langle\cdot,\cdot\rangle_{\omega}$ is locally equal to the pullback of some smooth metric on $\mathbb R^{N}$ to $\totl{\ex B}$.\end{remark}  

 To see that $\langle\cdot,\cdot\rangle_{\omega}$ is locally the pullback of some metric on $\mathbb R^{N}$, note that $\omega$ and $J'$  are locally equal to the restriction of some smooth two-form  and almost complex structure on $\mathbb R^{N}$. Therefore $\langle\cdot,\cdot\rangle_{\omega}$ is the restriction of some smooth symmetric form on $\mathbb R^{N}$ which is positive definite on the image of $\nT\totl{\ex B}$. As the image of $\nT\totl{\ex B}$ is closed, we may choose the extension of $\langle\cdot,\cdot\rangle_{\omega}$ to $\mathbb R^{N}$ to be smooth and positive definite, (ie a metric).
 
 The fact that $d_{\omega}$ is a nondegenerate metric on $\totl{\ex B}$ follows from the Remark \ref{metric embedding} above.
 
 \

 The next lemma is a version of the monotonicity lemma for holomorphic curves in $\totl{\ex B}$. This uses the pseudo-metric $d_{w}$ instead of a metric on $\ex B$. 
 
\begin{lemma}\label{omega monotonicity}
Given $(\ex B,J,\omega)$ and a choice of $\epsilon>0$, \begin{description}\item[]there exists some $E$ depending continuously on $(\ex B,J,\omega)$\item[] so that any non constant holomorphic curve in $\totl{\ex B}$ \begin{description}\item[] which passes through a point $p$\item[] and which is a proper map when restricted to the $d_{\omega}$-ball of size $\epsilon$ around $p$\end{description}  has $\omega$-energy greater than $E$.\end{description}
\end{lemma}
 
 \pf

 Using that $J$ is $\dbar\log$ compatible, Lemma \ref{J embedding}, and the fact that $\ex B$ is  complete gives  that it can be covered by a finite number of open subsets $U$ with holomorphic maps 
 \[f: (U,J)\longrightarrow(\mathbb R^{N},J')\]
 so that $\totl{f}$ is a holomorphic embedding.
 
\ 

  Using either Remark \ref{metric embedding} or the fact that $d_{\omega}$ induces the topology on $\ex B$ and $\totl{\ex B}$, 
  any $d_{\omega}$-ball of size $\epsilon$ centered around a point $p$ contained in our coordinate chart contains the inverse image of an $\epsilon'$-ball in $\mathbb R^{N}$ centered at $p$, and we may choose $\epsilon'$ to depend continuously on $p$ within our coordinate chart. For each point $p$, there is a positive constant $c$ so that \[\omega(v,Jv)\geq c(v\cdot v)\] where $v$ is any vector in $\nT_{p}\totl{\ex B}$ considered as a vector in $T_{p}\mathbb R^{N}$. The $c$ in the above equation may be chosen to depend continuously on $p$. It follows that the $\omega$-area of any holomorphic curve in $\totl{\ex B}$ is bounded below by a constant times the amount of its area  contained within an $\epsilon'$-ball  around $p$ in $\mathbb R^{N}$.
  
     For a given point $q\in\totl{U}$ in the smooth part of a coordinate chart $U$, the usual monotonicity lemma from \cite{Hummel} implies that there exists some lower bound  $E_{U}(q)>0$ for the $\omega$-energy of holomorphic curves which are  non constant in $\totl{U}$, pass through $q\in \totl U$ and are proper when restricted to the $d_{\omega}$-ball of size $\epsilon$ around $q\in\totl{U}$.  This lower bound $E_{U}(q)$ can be chosen to depend continuously on $q\in\totl{U}$.
  
   When we cover  $\ex B$ by a finite number of such open sets $U_{i}$, the function $E(q):=\max_{i}E_{U_{i}}(q)$ is positive and lower semicontinuous on $\totl {\ex B}$, which is compact. Therefore, there exists some lower bound $E>0$ for $E(q)$. In the case of a family $\hat {\ex B}\longrightarrow \ex G$, cover $\hat {\ex B}$ by coordinate charts so that each fiber is only contained in a finite number of coordinate charts, then $E(q)$ has a lower bound $E$ which can be chosen to depend continuously on $\ex G$. 

%
 
 \stop

%

 \begin{lemma}\label{cylinder bound}
  Given $(\ex B,J,\omega)$, for any  $\epsilon>0$, \begin{description}\item[]there exists some $E>0$ depending continuously on $(\ex B,J,\omega)$ \item[] so  that any $J$ holomorphic map  $f$ of $\{e^{-(R+1)}<\abs z<e^{(R+1)}\}\subset \mathbb C$ with $\omega$-energy less than $E$\begin{description}\item[] is contained inside a  $d_{\omega}$-ball of radius $\epsilon$  on the smaller annulus \[\{e^{-R}\leq\abs z\leq e^R\}\]
  \end{description}\end{description}
 
 \end{lemma}

 \pf
 
 Suppose to the contrary that this lemma is false. All distances in the following shall use $d_{\omega}$, the  pseudo metric defined by $\omega$.

 Choose $E$ small enough so that  Lemma \ref{omega monotonicity} holds for  $\frac \epsilon 8$ balls. There must be a path in $\{e^{-(R+1)}<\abs z< e^{R+1}\}$ joining $z=1$ with an end of this annulus whose image is contained entirely inside the ball $B_{\frac \epsilon 8}(f(1))$. (Otherwise, $f$ restricted to some subset would be a proper map to $B_{\frac \epsilon 8}(f(1))$, and we could use Lemma \ref{omega monotonicity} to get a contradiction.) Suppose without loss of generality that this path connects $1$ with the circle $\abs z=e^{-(R+1)}$.

Suppose for this paragraph that there exists some point $z$ so that $e^{-(R+\frac 34)}<\abs z<e^{-(R+\frac 14)}$ and $f(z)\notin B_{\frac{3\epsilon} 8}(f(1))$.
As above, there must exist some path connecting $z$ with an end which is contained in the ball $B_{\frac \epsilon 8}(f(z))$. Depending on which direction this path goes, the annulus $A$ of length $\frac 14$ to the right or left of $z$ has the property that every circle slice of the annulus contains at least one point with image in $B_{\frac \epsilon 8}(f(z))$, and one point with image in $B_{\frac \epsilon 8}(f(1))$, so each circle slice has $d_{\omega}$ length at least $\frac \epsilon 4$.  
  Then the Cauchy-Schwarz inequality tells us that 
\[\int_{A}\abs {df}^{2}_{\omega}:=\int_0^{2\pi}\int_0^{\frac 14}\abs{ df}^2_\omega
\geq \frac 2\pi \lrb{\frac \epsilon {16}}^{2} \]
This in turn gives a uniform lower bound for the energy of our curve. This case can therefore be discarded, and we can assume that all $z$ satisfying 
$e^{-(R+\frac 34)}<\abs z<e^{-(R+\frac 14)}$ are contained in $B_{\frac{3\epsilon} 8}(f(1))$.

As we are assuming this lemma is false, there must be some point $z_2$ so that $e^{-R}\leq\abs {z_2}\leq e^R$ and $f(z_2)\notin B_{\epsilon}(f(1))$. Repeating the above argument (with $z_2$ in place of the point $z=1$) gives us that $f(z)$ must be contained in $B_{\frac {3\epsilon}8}(f(z_2))$ when $e^{R+\frac 14}\leq \abs z\leq e^{R+\frac 34}$. 

We then have that  one end of our annulus is contained inside $B_{\frac{3\epsilon} 8}(f(1))$ and the other is contained inside $B_{\frac {3\epsilon}8}(f(z_2))$ where $f(z_{2})\notin B_{\epsilon}(f(1))$. Therefore we can apply  Lemma \ref{omega monotonicity} to some point in the image of $f$ with distance at least $\frac \epsilon 8$ from both balls. This gives  a lower bound for the energy and a contradiction to the assumption that the lemma was false.

\stop

\

The following gives useful coordinates for the analysis of holomorphic curves:
\begin{lemma}\label{coordinate bound}
Given  $(\ex B,J,\omega)$, \begin{description}\item[]there exist constants $c_i$ and $\epsilon>0$ depending continuously on $(\ex B,J,\omega)$\item[] so that for all points $p\in \ex B$,\item[] there exists some coordinate neighborhood $U$ of $p$ \begin{description}\item[] so that the smooth part, $\totl{U}$ can be identified with a  relatively closed subset of the open $\epsilon$ ball in $\mathbb C^n$ with almost complex structure $\hat J$ and some flat $\hat J$ preserving connection $\nabla$  so that
\begin{enumerate}
\item $p$ is sent to $0$.
\item The map $(U,J)\longrightarrow (\mathbb C^{n},\hat J)$ is holomorphic.
\item  The metric $\langle\cdot,\cdot\rangle_\omega$ on $\totl U$ is close to the standard flat metric on $\mathbb C^n$ so that 
\[\frac 1{c_{0}} \langle v,v\rangle_\omega\leq \abs v^2\leq c_{0} \langle v,v\rangle_\omega\]
\item\label{cb3} $\hat J$ at $0\in\mathbb C^{n}$ is the standard complex structure, and the derivative of $\hat J$ using the standard connection on $\mathbb C^{n}$ is bounded by $c_{0}$.

 \item The torsion tensor 
 \[\T_\nabla (v,w):=\nabla_vw-\nabla_wv-[v,w]\]
 is bounded by $c_0$, and has its $k$th derivatives bounded by $c_k$.
 \end{enumerate}
 \end{description}
\end{description}
\end{lemma}

Note that such a flat $J$ connection always locally exists, but may not exist globally. The only point in this lemma which does not follow from Lemma \ref{J embedding} and calculation in local coordinates is the fact that the constants involved do not have to depend on the point $p$. This follows from the fact that $\ex B$ is complete.

\begin{lemma}\label{dbar of derivative}
 Given a holomorphic map $f$ of the complex unit disk to an exploded manifold with a flat $J$ preserving connection $\nabla$, if $z=x+iy$ are standard coordinates on the disk, the partial derivative of $f$ in the $x$ direction $f_x$ defines a map to $\mathbb C^n$ defined by parallel transporting $f_x$ back to the tangent space at $f(0)$ which we then identify with $\mathbb C^n$. Such an $f_x$ satisfies the following equation involving the torsion tensor of $\nabla$, $\T_\nabla$.
 
 \[\dbar f_x=\frac 12 J\T_\nabla(f_y,f_x)\]
\end{lemma}

\pf

\[\begin{split}\dbar f_x&=\frac 12(\nabla_{f_x}f_x+J\nabla_{f_y}f_x)
 \\&=\frac 12 \nabla_{f_x}(f_x+J f_y) +\frac 12J(\nabla_{f_y}f_x-\nabla_{f_x}f_y)
 \\&=\frac 12 J\T_\nabla(f_y,f_x)  
  \end{split}\]
\stop

As the $\dbar$ above is the standard linear $\dbar$ operator, this expression is good for applying the following
standard elliptic regularity lemma for the linear $\dbar$ equation. (An exposition of the proof can be found in \cite{MS}). This allows us to get bounds on higher derivatives from bounds on the first derivative of holomorphic functions.

\begin{lemma}\label{elliptic regularity}
For for a given number $k$ and $1<p<\infty$, there exists a constant $c$ so that given any complex function $f$ on the unit disk $D(1)$,

\[\norm f_{L_{k+1}^p(D(\frac 12))}\leq c\left(\norm{\dbar f}_{L_k^p(D(1))}+\norm f _{L_k^p(D(1))}\right)\]

\end{lemma}

\begin{lemma}\label{derivative bound}
Given $(\ex B,J,\omega)$ and an energy bound $E_{0}$, \begin{description}\item[] there  exists\begin{description}\item[] some energy $E$,\item[] some distance $r>0$,\item[] and a constant $c$,\end{description} each depending continuously on $(\ex B,J,\omega,E_{0})$\item[] so that any holomorphic map $f$ of a disk $\{\abs z\leq 1\}$ into $\ex B$ which either\begin{description}\item[] has $\omega$-energy less than $E$\item[] or is contained inside a $d_{\omega}$-ball  of radius $r$ and has $\omega$-energy less than $E_{0}$,\end{description} satisfies 
\[\abs {df}_\omega <c\text{ at }z=0\] 
\end{description}
\end{lemma}
We shall omit the proof of the above Lemma which is a standard bubbling argument, similar but easier than the proof of the following:

\begin{lemma}\label{strict derivative bound}
Given $(\ex B,J,\omega,g)$, for any $\A>0$ and $E_{0}>0$, \begin{description}\item[] there exists\begin{description} \item[] some $\epsilon>0$,\item[] distance $r>0$, \item[] and constant $c$,\end{description} each depending continuously on $(\ex B,J,\omega,g,\A,E_{0})$ \item[] so that any  holomorphic map $f$ of a disk $\abs z\leq 1$ into $\ex B$ with local area bounded by $\A$, and either \begin{description}\item[]contained in a $d_\omega$-ball of radius $r$ and having $\omega$-energy less than $E_{0}$ \item[] or having $\omega$-energy less than $\epsilon$\end{description} satisfies
\[\abs {df}_g <c\text{ at }z=0\] 
\end{description}
\end{lemma} 

\pf

Suppose that this lemma was false. Then there would exist some sequence of maps $f_i$ satisfying the above conditions with $\abs {df_i(0)}\rightarrow\infty$. 
From this sequence, we could obtain a sequence of rescaled $J$ holomorphic maps $\tilde f_i
:D(R_i)\longrightarrow \ex B$ from the standard complex disk of radius $R_i$ 
so that
\[\abs{d \tilde f_i}\leq 2\]
\[\abs{d\tilde f_i(0)}=1\]
\[\lim_{i\rightarrow\infty}R_i=\infty\]

We achieve this in the same way as in any standard bubbling off argument such as the proof of lemma 5.11 in \cite{compact}. We can then use Lemma \ref{coordinate bound}, Lemma \ref{dbar of derivative}, and Lemma \ref{elliptic regularity} to get a bound on the higher derivatives of $\tilde f_i$ on $D(R_i-1)$.

As $\ex B$ is complete, we can choose a subsequence so that $f_i(0)$ converges to some $p\in \ex B$. Lemma \ref{bounded geometry} tells us that the geometry of $(g,J,\omega)$ around $f_i(0)$ converges to that around $p$.  By using the diffeomorphisms from Lemma \ref{bounded geometry}, consider all our maps as maps sending $0$ to $p$.  Then choose a subsequence that converges on compact subsets to a non constant holomorphic map  $f:\mathbb C\longrightarrow \ex B$. Note that $f$ either has $\omega$-energy less than $\epsilon$ or is contained in a ball of radius  $r$ in the $\omega$ pseudo-metric and has $\omega$-energy less than $E_{0}$. Because in either case, $f$ has finite $\omega$-energy, we can use Lemma \ref{cylinder bound} to tell us that $f$ must converge in the $\omega$ pseudo-metric as $\abs z\rightarrow\infty$. The standard removable singularity theorem for pseudo holomorphic curves proved in \cite{gromov} implies that $\totl{f}$ extends to a holomorphic map of $\mathbb CP^{1}$ to $\totl{\ex B}$. Then by choosing $\epsilon$ or $r$ small enough,  Lemma \ref{omega monotonicity} and the standard removable singularity theorem for holomorphic curves imply that $\totl f$ must be constant, because $f$ either has $\omega$-energy less than $\epsilon$ or is contained in a $d_\omega$-ball of radius $r$. 

Now that we have that the image of $f$ is a point in $\totl{\ex B}$, we know that $f$ must be contained entirely within some $(\mathbb C^*)^n$ worth of points over a single  point in $\totl{\ex B}$ and a single point in $\totb{\ex B}$. This $(\mathbb C^*)^n$ has the standard complex structure, so $f$ gives us $n$ entire holomorphic maps from $\mathbb C\longrightarrow \mathbb C^*$. As $f$ is non constant, at least one of these maps must be non constant. This map must have infinite degree, as it must have dense image in the universal cover of $\mathbb C^*$ which is just the usual complex plane. It therefore has infinite  local area, contradicting the fact that it must have local area less than $\A$.

To prove the family case, it is important to note that (as stated at the start of this section), when we say our constants depend continuously on  $(\ex B,J,\omega,g,\A,E_{0})$, we mean that given any finite dimensional family, $(\hat{\ex B}\longrightarrow\ex G,J,\omega,g,\A,E_{0})$, the constants can be chosen to depend continuously on $\ex G$. With this understood, the proof in the family case is analogous to the above proof.

\stop

 \begin{lemma}\label{omega cylinder bound}
 
  Given $(\ex B,J,\omega)$ and an energy bound $E_{0}$, 
  \begin{description}\item[] there exists\begin{description}\item[] an energy bound $E>0$,\item[] a distance $r>0$,\item[] and for any $0<\delta<1$ a constant $c$,
  \end{description} each depending continuously on $(\ex B,J,\omega,E_{0})$ \item[] so that any holomorphic map $f$ from the annulus $e^{-(R+1)}<\abs z<e^{(R+1)}$ to $\ex B$ with $\omega$-energy less than $E$, or with $\omega$-energy less than $E_{0}$ and contained inside a ball of $d_{\omega}$ radius $r$, satisfies the following inequality
  \[d_{\omega} (f(z),f(1))\leq ce^{-\delta R}(\abs z^\delta+\abs z^{-\delta})
  \ \ \ \ \ \ \ \ \ \text{ for } e^{-R}\leq\abs z\leq e^R\]  
  \end{description}
 \end{lemma}

 \pf

For this proof we shall use coordinates $z=e^{t+i\theta}$, and the cylindrical metric in which $\{\frac \partial{\partial t},\frac\partial{\partial \theta}\}$ are an orthonormal frame.

Choose our energy bound $E>0$ small enough that we can use the derivative bound from Lemma \ref{derivative bound}, and Lemma \ref{cylinder bound} implies that the smaller annulus is contained a small enough ball that we can use the coordinates of Lemma \ref{coordinate bound} (we also choose $r$ small enough that this is true).  

Item \ref{cb3} of Lemma \ref{coordinate bound} states that the almost complex structure $J$ used in these coordinates is equal to the standard complex structure $J_{0}$ at $0$, and that the derivative of $J$ is bounded by a constant $c_{1}$. Therefore, we may estimate the standard $\dbar$ operator applied to $f$ in these coordinates as follows:
  \[\abs{\dbar f}:=\frac 12\abs{df+J_{0}\circ df\circ j}=\frac 12 \abs{(J_{0}-J)\circ df\circ j}\leq c_1\abs {df}\abs f\]

Let $A_{p}$ be an affine coordinate change on our coordinates from Lemma \ref{coordinate bound} which sends $p$ to $0$ so that $J$ in these new coordinates is equal to the standard complex structure at $0$. Then, the derivative of $J$ in these new coordinates is bounded by $c_{1}\norm {dA_{p}}\norm{dA_{p}^{-1}}$. We may choose such coordinate transformations $A_{p}$ for all $p$ in a ball with size depending on $c_{1}$ so that $\norm {dA_{p}}\norm{dA_{p}^{-1}}$ is bounded by $2$. Choosing our energy bound $E>0$  or our distance $r>0$ small enough means that we can ensure that the image of $\{e^{-R}\leq\abs z\leq e^R\}$ is contained inside this ball. Then
\[\abs {\dbar A_{p} f}\leq 2c_{1}\abs{dA_{p}f}\abs{A_{p}f}\]

We can run the above inequality through the inequality from Lemma \ref{elliptic regularity} to obtain the following estimate on $\dbar f$ restricted to the interior of a disk (which holds on the interior of the cylinder, where we have the derivative bound from Lemma \ref{derivative bound}.) 

\begin{equation}\label{dbar bound}\begin{split}\norm{\dbar A_{p} f}_{L^q(D(\frac 12))}&\leq 2c_1\norm {A_{p}f}_{L^q_1(D(\frac 12))}\norm {A_{p}f}_{L^\infty(D(\frac 12))}
\\ &\leq c\lrb{\norm {A_{p}f}_{L^{q}(D(1))}+\norm{\dbar A_{p}f}_{L^q (D(1))}}\norm {A_{p}f}_{L^\infty(D(\frac 12))}
\\&\leq c\lrb{\norm {A_{p}f}_{L^{q}(D(1))}+\norm{2c_{1}\abs{d{A_{p}f}}\abs {A_{p}f}}_{L^q (D(1))}}\norm {A_{p}f}_{L^\infty(D(\frac 12))}
\\& \leq c_2\norm {A_{p}f}_{L^\infty (D(1))}\norm {A_{p}f}_{L^\infty(D(\frac 12))}\end{split}\end{equation}  
The second inequality uses Lemma \ref{elliptic regularity} and the last inequality uses the derivative bound from  Lemma \ref{derivative bound}.
 Here $c_2$ is a constant depending only on the $q$ from $L^{q}$ and $(\ex B, J, \omega,E_{0})$. We can fix $q$ to be something bigger than $2$. By choosing $E$ or $r$ small, we can force $\abs{A_{p} f}$ to be as small as we like on the smaller cylinder using Lemma \ref{cylinder bound}. 
 
 Now we can use Cauchy's integral formula 
 
 \[2\pi i A_{p}f(z_0)=-\int_{\abs z=1}\frac {A_{p}f(z)}{z-z_0}dz+\int_{\abs z=e^{2l}}\frac {A_{p}f(z)}{z-z_0}dz+\int_{1\leq\abs z\leq e^{2l}}\frac{\dbar A_{p}f(z)}{z-z_0}\wedge dz\]
 or in our coordinates, 
 \[\begin{split} A_{p}f(t_0,\theta_0)&=-\frac 1{2\pi}\int_0^{2\pi} \frac{A_{p}f(0,\theta)}{1-e^{t_0}e^{i(\theta_0-\theta)}}d\theta+
 \frac 1{2\pi}\int_0^{2\pi} \frac{A_{p}f(2l,\theta)}{1-e^{t_0-2l}e^{i(\theta_0-\theta)}}d\theta
 \\&+ \frac 1{2\pi}\int_0^{2\pi}\int_0^{2l}\frac {\dbar A_{p}f(\theta,t)}{1-e^{t_0-t}e^{i(\theta_0-\theta)}}dtd\theta \end{split}\]
 
 Let us now consider each term of this expression for $A_{p}f(l,\theta_0)$ in the middle of a cylinder under the assumption that the average of $f(2l,\theta)$ is $p$, so the average of $A_{p}f(2l,\theta)$ is $0$.
 
 The first term:
 \[\abs{\frac 1{2\pi}\int_0^{2\pi} \frac{A_{p}f(0,\theta)}{1-e^{l}e^{i(\theta_0-\theta)}}d\theta}\leq\frac 1{e^l-1}\max \abs{A_{p}f(0,\theta)}\]
 
 The second term using that the average of $A_{p}f(2l,\theta)$ is $0$,
 \[\begin{split}\abs{\frac 1{2\pi}\int_0^{2\pi} \frac{A_{p}f(2l,\theta)}{1-e^{l-2l}e^{i(\theta_0-\theta)}}d\theta}
  &=\abs{\frac 1{2\pi}\int_0^{2\pi} A_{p}f(2l,\theta)\left(\frac 1{1-e^{-l}e^{i(\theta_0-\theta)}}-1\right)d\theta}
   \\&=\abs{\frac 1{2\pi}\int_0^{2\pi} A_{p}f(2l,\theta)\left(\frac {e^{-l}e^{i(\theta_0-\theta)}}{1-e^{-l}e^{i(\theta_0-\theta)}}\right)d\theta}
   \\&\leq \frac 1{e^l-1}\max\abs{A_{p}f(2l,\theta)}
   \end{split}\]

   The third term:
 \[\begin{split}\abs{\frac 1{2\pi}\int_0^{2\pi}\int_0^{2l}\frac {\dbar A_{p}f(\theta,t)}{1-e^{t_0-t}e^{i(\theta_0-\theta)}}dtd\theta}
    &
\leq\norm{\dbar A_{p} f}_{L^3}\frac 1{2\pi}\norm{\frac 1{1-e^{l-t}e^{i(\theta_0-\theta)}}}_{L^{\frac 32}}
\\& \leq c(l+1) \norm{\dbar A_{p}f}_{L^3}
\\&  \leq 
c_3 (l+1)^{2}\left(\max_{t\in[-1,2l+ 1]}\abs {A_{p}f(t,\theta)}\right)
\left(\max_{t\in[0,2l]}\abs {A_{p}f(t,\theta)}\right)
\end{split}\]  
The last inequality above uses the inequality (\ref{dbar bound}).  
 The constant $c_3$ depends only on $(\ex B,J,\omega)$. It is zero if $J$ is integrable.
 
 To summarize the above, we have the following expression which holds if the average of $f(\theta,2l)$ is $p$.  
 \begin{equation}\label{c estimate}\begin{split}
  \abs {A_{p}f(l,\theta)}\leq & \frac 1{e^l-1}\left(\max\abs{A_{p}f(0,\theta)}+\max\abs {A_{p}f(2l,\theta)}\right)
  \\&+ c_3(l+1)^{2}\left(\max_{t\in[-1,2l+ 1]}\abs {A_{p}f(t,\theta)}\right)
\left(\max_{t\in[0,2l]}\abs {A_{p}f(t,\theta)}\right)
\end{split}
 \end{equation}

Stokes' theorem applied to $A_{p}f(dt+id\theta)$ implies that the amount that the average of $A_{p}f(t_0,\theta)$  changes with $t_0$ inside $[0,2l]$ is determined by the integral of $\dbar A_{p} f$, which is dominated as above by a term of the form 
\begin{equation}\begin{split}\label{average estimate}\max_{t_{0},t_{1}\in[0,2l]}&\abs{\frac 1{2\pi}\int_{0}^{2\pi} A_{p}f(t_{0},\theta)d\theta  -\frac 1{2\pi}\int_{0}^{2\pi} A_{p}f(t_{1},\theta)d\theta }
\\ &\leq c_4(l+1)^{2} \left(\max_{t\in[-1,2l+ 1]}\abs {A_{p}f(t,\theta)}\right)
\left(\max_{t\in[0,2l]}\abs {A_{p}f(t,\theta)}\right)\end{split}\end{equation}

As we have chosen $A_{p}$ so that $\norm{dA_{p}}\norm{dA_{p}^{-1}}\leq 2$, we may remove the dependence on $A_{p}$ of the above inequalities.  Define the variation of $f(t,\theta)$ on the cylinder $[a,b]\times S^{1}$ as follows:
 \[Vf([a,b]):=\max_{(t_{i},\theta_{i})\in[a,b]\times S^{1}}\abs {f(t_{1},\theta_{1})-f(t_{2},\theta_{2})}\]
The inequality (\ref{c estimate}) implies that 
\begin{equation}\label{v est}Vf[x,x]\leq \frac {c_{5}}{e^l-1}Vf[x-l,x+l]  + c_{5}(l+1)^{2}Vf[x-l-1,x+l+1]Vf[x-l,x+l]\end{equation}
and the inequality (\ref{average estimate}) implies that the change in the average of $f$ over the cylinder $[a,a+l]$ may be estimated by 
\[c_{6}(l+1)^{2} Vf[a-1,a+l+1]Vf[a,a+l]\]
therefore we may estimate $Vf[a,a+l]$ by
\[Vf[a,a+l]\leq c_{6}(l+1)^{2} Vf[a-1,a+l+1]Vf[a,a+l]+ \sup_{x\in [a,a+l]}Vf[x,x]\]
Applying inequality (\ref{v est}) above then gives
\begin{equation}\begin{split}\label{Vf1}Vf[a,a+l]&\leq c_6(l+1)^{2} Vf[a-1,a+l+1]Vf[a,a+l]
\\ &\ \ + \sup_{x\in [a,a+l]}\Big(\frac {c_{5}}{e^l-1}Vf[x-l,x+l] 
\\ & \ \ \ \ \ \ \ \ \ \ \ \ \ \ \ \ + c_5(l+1)^{2}Vf[x-l-1,x+l+1]Vf[x-l,x+l]\Big)
 \\&\leq c_{7}\lrb{\frac 1{e^{l}-1}+(l+1)^{2}Vf[a-l-1,a+2l+1]}Vf[a-l,a+2l]\end{split}\end{equation}

We shall now complete the argument for one particular $\delta$; say $\delta=\frac 12$. 
If we choose $l$ large enough, and then make $Vf$ small enough by making our energy bound $E$ small (or directly making $r$ small), we can get the estimate
\begin{equation}\label{Vf2} Vf[a,a+l]\leq \frac {e^{-\delta l}}3Vf[a-l,a+2l]
\ \ \ \ \ \text{ for }[a-l-1,a+2l+1]\subset[-R,R]\end{equation}

How large $l$ needs to be, and how small $E$ or $r$ need to be depends continuously on $(\ex B,J,\omega,E)$. (We have no need for this inequality when $R$ is not much larger than  $l$, because our lemma's estimate for any bounded distance from the edges of the cylinder follows immediately from Lemma \ref{cylinder bound}.) Applying inequality (\ref{Vf2})  three times, we get that on the appropriate part of the cylinder,  
\[\begin{split}Vf[a-l,a+2l]&\leq Vf[a-l,a] +Vf[a,a+l]+Vf[a+l,a+2l]
\\&\leq \frac {e^{-\delta l}}3\lrb{Vf[a-2l,a+l]+Vf[a-l,a+2l]+Vf[a,a+3l]}
\\ &\leq e^{-\delta l}Vf[a-2l,a+3l]\end{split}\]
therefore
\[Vf[a,a+l]\leq\frac {e^{-\delta 2l}}3Vf[a-2l,a+3l]\]

Inductively continuing this argument gives that if $[a-nl,a+(n+1)l]\subset[-R+1,R-1]$,
\begin{equation}\label{vf3}Vf[a,a+l]\leq \frac {e^{-\delta n l}}3 Vf[a-nl,a+(n+1)l]\end{equation}
In particular, if the above inequality holds for $n$, then each of $Vf[a-l,a]$, $Vf[a,a+l]$ and $Vf[a+l,a+2l]$ is each bounded by $\frac{e^{-\delta nl}}3Vf[a-(n+1)l,a+(n+2)l]$. Then applying inequality (\ref{Vf2}) gives that inequality (\ref{vf3}) also holds when $n$ is replaced by $n+1$.

The required estimate follows from inequality (\ref{vf3}). In particular, because Lemma \ref{coordinate bound} stipulates that $d_{\omega}$ is comparable to the norm we have been using in our coordinates,   $d_{\omega}(f(z),f(1))$ can be bounded by the sum of the variation of adjacent cylinders of length $l$ starting with a cylinder containing $z$ and ending with a cylinder containing $1$. If $e^{R}>\abs z>1$, then  $z$ may be placed in a cylinder of length $l$ surrounded by at least $(R-\log \abs z)/l-2$ other cylinders, so inequality \ref{vf3} implies that the variation on this cylinder containing $z$ is bounded by $\abs z^{\delta}e^{(-R+2l)\delta}V[-R,R]$. The bound on the variation of adjacent cylinders improves exponentially until the cylinder containing $1$ is reached, therefore for $e^{R}>\abs z>1$,  $d_{\omega}(f(z),f(1))$ is bounded by $\abs z^{\delta}e^{-R}$ times a constant which is independent of $R$. Similarly, if $e^{-R}<\abs z<1$, then $d_{\omega}(f(z),f(1))$ is bounded by $\abs z^{-\delta}e^{-R}$ times a constant which is independent of $R$. Therefore, $d_{\omega}(f(z),f(1))$ is bounded by $(\abs z^{\delta}+\abs z^{-\delta})e^{-R}$ times a constant which is independent of $R$.

\

Note that once we have inequality (\ref{Vf2}) for any $\delta<1$, following the above argument verbatim gives the required result that $d_{w}(f(z),f(1))$ is bounded by $(\abs z^{\delta}+\abs z^{-\delta})e^{-R}$ times a constant which is independent of $R$. Above, we proved inequality (\ref{Vf2}) from (\ref{Vf1})  for any particular $\delta<1$ by choosing $l$ large enough, then forcing $Vf$ to be small by making our energy bound $E$ small or making $r$ small. As we don't want to use different $E$ or $r$ bounds for different $\delta$, we need some other way of making $Vf$ small. Our already proven estimate for $\delta=1/2$ gives that we may force $Vf$ to be as small as we like simply by restricting some distance $d$ into the interior of the cylinder.
 Therefore, for any $\delta<1$, we may obtain inequality (\ref{Vf2}) from (\ref{Vf1}) by choosing $l$ large enough, then restricting attention a distance $d$ into the interior of our cylinder. This constant $d$ depends only on $(\ex B,J,w,E_{0},\delta)$. Repeating the above argument  gives the required inequality (with a different constant). As $d_{w}(f(z),f(1))$ is already bounded by a constant (by assumption or using Lemma \ref{cylinder bound}), the  fact that our argument only applies after we have shortened our cylinder by $d$ only affects the bound by multiplying by some $d$-dependent constant.

%
%
%
%
%
%

 \stop

\

We may derive similar bounds for derivatives of $f$ measured by $d_{\omega}$ using standard elliptic bootstrapping as follows: Use coordinates from Lemma \ref{coordinate bound} centered on $f(1)$. Then  Lemma \ref{omega cylinder bound} and inequality (\ref{dbar bound}) on page \pageref{dbar bound} imply that
if $D_{z}(r)$ indicates the disk of radius $r$ around $z$, and all norms below indicate the norms in the standard metric from Lemma \ref{coordinate bound}, 
\[\norm{\dbar f}_{L^{p}(D_{z}(1))}\leq c'\norm{f}_{L^{\infty}(D_{z}(2))}^{2}\leq ce^{-2R\delta}(\abs z^{\delta}+\abs z^{-\delta})^{2}\]
 for some $c$ depending only on $E,r,\ex B,J,\omega$ and $p>2$.  
 
 Then Lemma \ref{elliptic regularity} implies that 
 \begin{equation}\label{eineq}\norm f_{L^{p}_{1}(D_{z}(1/2))}\leq c e^{-R\delta}(\abs z^{\delta}+\abs z^{-\delta})\end{equation}
 for some different constant $c$ depending on the same quantities as the previous constant. Then Lemma \ref{dbar of derivative} and the H\"older inequality implies that 
 \[\begin{split} \norm{\dbar f_{x}}_{L^{p}(D_{z}(1/2))} &\leq \norm {c\abs {f_{x}}\abs{f_{y}}}_{L^{p}(D_{z}(1/2))}
 \\ &\leq c\norm{f_{x}}_{L^{2p}(D_z(1/2))}\norm{f_{y}}_{L^{2p}(D_z(1/2))}
 \\& \leq c\norm{df}_{L^{2p}(D_z(1/2))}^{2}\end{split}\]
so applying our previous inequality (\ref{eineq}) with $2p$ in place of $p$ gives 
\[\norm{\dbar f_{x}}_{L^{p}(D_{z}(1/2))}\leq ce^{-2R\delta}(\abs z^{\delta}+\abs z^{-\delta})^{2}\]
with a different constant $c$ depending on the same quantities. (A similar inequality  applies to $f_{y}$). Then applying Lemma \ref{elliptic regularity}
 gives
 \[\norm{f}_{L^{p}_{2}(D_{z}(1/4))}\leq c e^{-R\delta}(\abs z^{\delta}+\abs z^{-\delta})\]
 The argument can be continued to produce decay of $\norm f_{L^{p}_{k}(D_{z}(2^{-k}))}$ by differentiating the equation from Lemma \ref{dbar of derivative}, and noting that the derivatives of the torsion tensor have universal bounds from Lemma \ref{coordinate bound}. We do not need to continue at this point, but merely to note that Sobolev embedding  gives that  our $L^{p}_{2}$ bound above implies a bound on the first derivative. As Lemma \ref{coordinate bound} tells us that the standard metric in our coordinate chart is equivalent to the $\omega$ pseudometric,   there exists some constant $c$ depending  on $(E,r,\ex B,J,\omega)$ so that
  \begin{equation}\label{df wbound}\abs {df(z)}_{\omega}\leq c e^{-R\delta}(\abs z^{\delta}+\abs z^{-\delta})\end{equation}
 Lemma \ref{omega cylinder bound} supplemented with equation (\ref{df wbound}) has the following immediate corollary:
\begin{cor}\label{cylinder energy decay}
Given $(\ex B,J,\omega)$,
  \begin{description}
  	 \item[] there exists some energy $E>0$ depending continuously on $(\ex B,J,\omega)$
   \begin{description}\item[] so that given any $E'>0$, 
   \begin{description}\item[]there exists some distance $R$ depending continuously on $(\ex B,J,\omega,E')$
   \begin{description}\item[]so that given any holomorphic map $f$ of a cylinder \[\{e^{-l-R}<\abs z<e^{l+R}\}\] to $\ex B$ with $\omega$-energy less than $E$,\item[] the $\omega$-energy of $f$ restricted to $e^{-l}<\abs z<e^{l}$ is less than $ E' $. 
   \end{description}
   \end{description}\end{description}
\end{description}
\end{cor}

 \begin{lemma}\label{strong cylinder convergence}
  Given  $(\ex B,J,\omega,g)$,  a local area bound  $\A<\infty$ and energy bound $E_{0}$, and a covering of $\ex B$ by  coordinate charts,

  \begin{description}\item[] there exists\begin{description}\item[] an open cover $\{U_{i}\}$ of $\ex B$\item[] and  for each $0<\delta<1$ a constant $c$  \end{description}
  \item[] so that the following is true:
 
  \begin{description}
  \item[]
  Each member of our open cover $U_{i}$ is contained in one of our chosen coordinate charts on $\ex B$.
  \item[]
  Given any holomorphic map $f$ with local area bounded by $\A$ and $\omega$-energy bounded by $E_{0}$ from a cylinder $e^{-(R+1)}\leq \abs z\leq e^{(R+1)}$ to $\ex B$ contained inside $U_{i}$,\begin{description}\item there exists a map $F$ given in coordinates as 
  \[F(z):=(c_1z^{\alpha_1}\e{a_1},\dotsc,c_kz^{\alpha_k}\e{a_k}, c_{k+1},\dotsc,c_n)\]
 \item[]so that 
 \[\dist\left( f(z), F(z)\right)\leq ce^{-\delta R}\left(\abs z^\delta+\abs{z}^{-\delta}\right)\text{ for }e^{-R}\leq\abs{z}\leq e^{R}\]
  where for $i
  \in [1,k]$, $c_i\in\mathbb C^*$, $a_i\in\mathbb R$,  $\alpha_i\in \mathbb Z$, and for $i\in [k+1,n]$, $c_{i}\in\mathbb R$. In the above, $\dist$ indicates distance in the metric $g$.
   \end{description} 
    \
    
 \end{description}

 \end{description}
  In the case that we a have a family $(\hat{\ex B}\longrightarrow\ex G,J,\omega,g,\A,E_{0})$ with the above structure, we can choose $c$ continuous on $\ex G$.   

 \end{lemma}

 \pf 
 
 Because $\ex B$ is complete, we can  choose a finite covering $\{U_{i}\}$ of  $\ex B$ by open subsets which are contained inside our coordinate charts and which are small enough in the $\omega$ pseudo metric so that Lemma \ref{omega cylinder bound} tells us that for any holomorphic map from a cylinder as above contained inside $U_{i}$, 
 \begin{equation}\label{omega distance estimate}d_\omega(f(z),f(1))\leq ce^{-\delta  R}(\abs z^{\delta }+\abs z^{-\delta })\end{equation}  In the family case, we choose our open cover $\{U_{i}\}$ to be fiberwise finite, and the above constant $c$ depends continuously on $\ex G$.

 In our coordinates,
 
 \[f(1)=(c_1\e{a_1},\dotsc,c_k\e{a_k}, c_{k+1},\dotsc,c_n)\]
 
 We can choose $\alpha_i$ so that the winding numbers of the first $k$ coordinates of $f(e^{i\theta})$ are the same as our model map $F$:
 \[F:=(c_1z^{\alpha_1}\e{a_1},\dotsc,c_kz^{\alpha_k}\e{a_k}, c_{k+1},\dotsc,c_n)\] 
 
 The metric $g$ can be compared to $d_{\omega}$ on the last coordinates $c_{k+1},\dotsc,c_n$, so we only need to prove convergence in the first $k$ coordinates. Let us do this for $f$ restricted to the first coordinate $f_1$.
 
Use coordinates $z=e^{t+i\theta}$ on the domain with the usual cylindrical metric and the similar  cylindrical metric on our target. (Choose members of our open cover to be compactly contained inside our coordinate charts so that this standard flat metric on the coordinate chart is uniformly comparable to $g$).
Use the notation 
\[\frac{f_1(z)}{c_1z^{\alpha_1}\e {a_1}}=e^{h(z)}\]
where $h$ is a $\mathbb C$ valued function so that $h(1)=0$. Note that this standard flat metric has as an orthonormal basis the real and imaginary parts of $\tilde z\frac\partial{\partial \tilde z}$, so $\dbar f_{1}$ measured using this cylindrical metric is equal to $\dbar h$ measured using the standard Euclidean metric on $\mathbb C$.

 As Lemma \ref{omega cylinder bound} holds, equation (\ref{df wbound}) on page \pageref{df wbound} also holds, so

 \begin{equation}\label{lkj}\abs{df(z)}_\omega\leq ce^{-\delta  R}(\abs z^{\delta }+\abs z^{-\delta })\text{ for }e^{-R}\leq \abs z\leq e^{R}\end{equation}

 (The $c$ in the above inequality is some new constant which is independent of $f$). So long as we have chosen members of our open cover $U_{i}$ compactly contained inside  our coordinate charts, the following inequality holds inside $U_{i}$:
 \begin{equation}\label{J bound}\abs{Jv-J_{0}v}_{g}\leq M{\abs v_{\omega}}\end{equation}
 where $J_{0}$ indicates the standard complex structure, and $M$ is a constant, which may have to depend on $\ex G$ in the case of a family $\hat{\ex B}\longrightarrow \ex G$. 
   To prove the inequality (\ref{J bound}) above, note that  $J_{0}$ is $\dbar \log$ compatible, so Lemma \ref{J construction} implies that $(J-J_{0})$ sends smooth one-forms to one-forms generated by functions. It follows that $g((J-J_{0})\cdot,(J-J_{0})\cdot)$ is equal to the pullback under $\iota:T\ex B\longrightarrow \nT\totl{\ex B}$ of some smooth  section $g'$ of $T^{*}\totl{\ex B}\otimes T^{*}\totl{\ex B}$. On the other hand, as $J$ is tamed by $\omega$, the $\omega$ pseudo metric is the pullback under $\iota$ of some smooth metric on $\totl{\ex B}$. The above inequality (\ref{J bound}) follows.

    The above inequality (\ref{J bound}) implies that  
  $\abs {\dbar f_1}_g$ (using $J_{0}$) in our coordinates is controlled by $\abs{df}_\omega$ because
  \[\abs {\dbar f_{1}}_{g}= \frac 12\abs{df_{1}-J\circ df_{1}\circ j+(J-J_{0})\circ df_{1}\circ j}_{g}\leq \frac M2\abs {df}_{\omega}\]
  To understand the above, $df_{1}$ should be considered as the first component of $df$ using our coordinates, so it is a one form with values in $f^{*}T\ex B$. 
  The above inequality and (\ref{lkj}) then give that for some new constant $c$ independent of $f$ we get the following inequality:
   
 \begin{equation}\label{dbar h estimate}\abs {\dbar h(z)}\leq c e^{-\delta  R}(\abs z^{\delta }+\abs z^{-\delta })\text{ for }e^{-R}\leq \abs z\leq e^{R}\end{equation}  
 where $\dbar h$ is now measured using the standard Euclidean metric on its target $\mathbb C$, and the cylindric metric on its domain.
If we choose our open cover small enough using $d_{\omega}$, Lemma \ref{strict derivative bound} gives us a bound for $\abs{df}_g$ on the interior of the cylinder, so we have that there exists some $c$ independent of $f$ so that 
 
 \begin{equation}\label{h derivative bound}\abs {dh(z)}< c \text{ for }e^{-R}\leq\abs z\leq e^R\end{equation}
 
 We can now proceed roughly as we did in the proof of Lemma \ref{omega cylinder bound}. In particular, Cauchy's integral theorem tells us that 
  
 \[\begin{split} h(t_0,\theta_0)&=-\frac 1{2\pi}\int_0^{2\pi} \frac{h(t_0-l,\theta)}{1-e^{l}e^{i(\theta_0-\theta)}}d\theta+
 \frac 1{2\pi}\int_0^{2\pi} \frac{h(t_0+l,\theta)}{1-e^{-l}e^{i(\theta_0-\theta)}}d\theta
 \\&+ \frac 1{2\pi}\int_0^{2\pi}\int_{t_0-l}^{t_0+l}\frac {\dbar h(t,\theta)}{1-e^{t_0-t}e^{i(\theta_0-\theta)}}dtd\theta \end{split}\]

 Cauchy's formula is also valid if we shift $h$ so that its average around the loop at $t_{0}+l$ is $0$. We may bound each term in the shifted formula as in the proof of inequality (\ref{c estimate}) in the proof of Lemma \ref{omega cylinder bound}, except we use the estimate (\ref{dbar h estimate}) to bound $\abs {\dbar h}$.

 Then we have the following estimate:
  \[\begin{split}
  \abs {h(t_0,\theta_0)-\frac 1{2\pi}\int_0^{2\pi}h(t_0+l,\theta)d\theta }\leq & \frac 1{e^l-1}\max\abs{h(t_0-l,\theta)-\frac 1{2\pi}\int_0^{2\pi}h(t_0+l,\theta)d\theta}
  \\&+\frac 1{e^l-1}\max\abs {h(t_{0}+l,\theta)-\frac 1{2\pi}\int_0^{2\pi}h(t_0+l,\theta)d\theta}
  \\&+ c(l+1)^{2}e^{-\delta R}e^{\delta l}\left( e^{\delta  t_{0}}+e^{-\delta t_{0}}\right)
\end{split}
 \]

(Of course, this is a new constant $c$, which is independent of $l$ or $h$.)

\

 Note that the change in the average of $h$ is determined by the integral of  $\dbar h$, which we can bound using estimate (\ref{dbar h estimate}). At the cost of increasing the constant $c$ in the inequality above, each term in absolute values above may therefore be replaced with the maximum difference between $h$ and its average at $t_{0}$, $t_{0}-l$ and $t_{0}+l$ respectively.   Define the variation of $h$ for a particular $t$ as follows:
 
 \[Vh(t):=\max_{\theta}\abs {h(t,\theta)-\frac 1{2\pi}\int_0^{2\pi}h(t,\theta)d\theta }\]
 We then have the following estimate (with a new constant $c$):
 \[Vh(t)\leq \frac 1{e^l-1}(Vh(t+l)+Vh(t-l))+ c(l+1)^{2}e^{\delta (l-R)}\left( e^{\delta  t}+e^{-\delta t}\right)\]
Fix $l$ large enough that $1/(e^{l}-1)<e^{-\delta l}/4$. Then 
\[Vh(t)\leq \frac{e^{-\delta l}}4 (Vh(t+l)+Vh(t-l))+X(t)\]
where
\[X(t):= ce^{-\delta R}(e^{\delta t}+e^{-\delta t})\] 
and $c$ is new constant which depends on $l$. 
 In order to apply the above inequality inductively, define $Ax(t):=\frac 12(x(t-l)+x(t+l))$.
 Applying the above inequality to estimate $Vh(t-l)$ and $Vh(t+l)$ then running these estimates through the above inequality again gives
\[Vh(t)\leq \frac {e^{-2l\delta}}4A^{2}Vh(t)+\frac {e^{-\delta l }}2 AX(t)+X(t)\]
Inductively continuing this estimation gives
\[Vh(t)\leq \frac {e^{-n\delta l}}{2^{n}}A^{n}Vh(t)+\sum_{i=0}^{n-1}\frac{e^{-i\delta l}}{2^{i}}A^{i}X(t)\]
Of course, the above inequality only holds if $t$ is a distance at least $nl$ from boundary of $[-R,R]$. Note that $A^{i}X(t)\leq e^{i\delta l}X(t)$ so the last term in the above inequality is bounded by $2X(t)$. The observation that $A^{n}Vh(t)$ is bounded independent of $t$ and $n$ then gives for some new constant $c$,  
 \[Vh(t)\leq ce^{-\delta R}(e^{\delta t}+e^{-\delta t})\text{ for }-R\leq t\leq R\] 
  The required estimate for $h$ then follows from the fact that the change in the average of $h$ is bounded by the estimate (\ref{dbar h estimate}). We therefore have
 \[h(z)\leq ce^{-\delta R}(\abs z^{\delta}+\abs z^{-\delta})\text{ for }e^{-R}\leq\abs z\leq e^{R}\]
 which is the required estimate.
 
 \stop

\

To prove compactness results, we shall divide our domain up into annuli with small energy, and other compact pieces with derivative bounds. For this we shall need some facts about annuli.
Recall the following standard definition for the conformal modulus of a Riemann surface which is an annulus:

\begin{defn}
The conformal modulus of an annulus $A$ is defined as follows. Let $S(A)$ denote the set of all continuous functions with $L^2$ integrable derivatives on $A$ which approach $1$ at one boundary of $A$ and $0$ at the other. Then the conformal modulus of $A$ is defined as 
\[R(A):=\sup_{f\in S(A)}\frac {2\pi}{\int_A(df\circ j)\wedge df}=\sup_{f\in S(A)}\frac {2\pi}{\int_A\abs{df}^2}\] 
\end{defn}

We can extend the definition of conformal modulus to include `long' annuli inside exploded curves as follows:

\begin{defn}
Call a (non complete) exploded curve $\ex A$ an exploded annulus of conformal modulus 
$\log x\e {-l}$ if it is connected, and there exists an injective holomorphic map $f:\ex A\longrightarrow \ex T$ with image $\{1<\abs {\tilde z}< x\e {-l}\}$. Call it an exploded annulus with semi infinite conformal modulus if it is equal to (a refinement of) $\{\abs{\tilde z}<1\}\subset\et 11$, and call it an exploded annulus with infinite conformal modulus if it is a refinement of $\ex T$.
\end{defn}

For interpreting the above, note that that the absolute value of $c\e a\in\mathbb C^{*}\e {\mathbb R}$ is defined as $\abs {c\e a}:=\abs c\e a\in(0,\infty)\e{\mathbb R}$.
(We use $x\e{-l}$ in the above definition because the tropical part of the resulting annulus will have length $l$.) 
The use of `$\log$' in the above is just notation, we do not attempt here to define a function log on $(0,\infty)\e{\mathbb R}$, however in the case that $l=0$, we may regard $x\e 0$ simply as $x\in (1,\infty)$, and then $\log x$ is the usual conformal modulus $R$.  
Two exploded annuli with the same conformal modulus may not be isomorphic, but they will have a common refinement.

\

We shall need the following lemma containing some useful properties of the (usual) conformal modulus.

\begin{lemma}
\label{annulus lemma}
\begin{enumerate}
\item\label{conformal modulus 1} An open annulus $A$ is conformally equivalent to \[\{1<\abs z< e^R\}\] if and only if  the conformal modulus of $A$ is $R$. If the conformal modulus of $A$ is infinite, then $A$ is conformally equivalent to either a punctured disk, or a twice punctured sphere.
\item\label{conformal modulus 2} If $\{A_i\}$ is a set of disjoint annuli $A_i\subset A$ none of which bound a disk in $A$, then 
\[R(A)\geq\sum R(A_i)\]
\item\label{conformal modulus 3}

 If $A_1$ and $A_2$ are annuli contained inside the same Riemann surface  so that
  \begin{enumerate}
  \item one connected component of the  boundary of $A_{1}$ is equal to a connected component of the boundary of $A_{2}$, and the orientations on these components agree,
  \item $R(A_2)<\infty$,\item and the other boundary of $A_1$ intersects the other boundary of $A_2$ and the circle at the center of $A_2$,\end{enumerate} then 
\[R(A_2)+\frac{16\pi^2}{R(A_2)}\geq R(A_1)\]
(This inequality is not sharp.)

\includegraphics{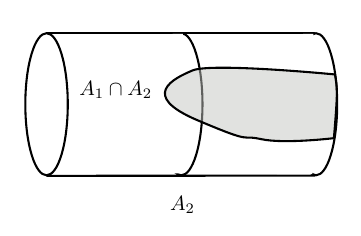}

\item \label{conformal modulus 4} If $A_1$ and $A_2$ are annuli contained inside a Riemann surface so that  every circle homotopic to the boundary inside $A_2$ contains a segment inside $A_1$ that intersects both boundaries of $A_1$, then
\[R(A_1)\leq\frac {4\pi^2}{R(A_2)}\]

\includegraphics{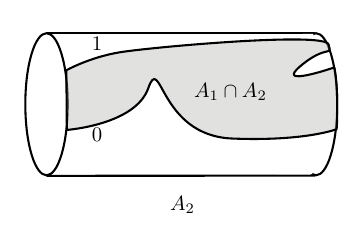}

\end{enumerate}

\end{lemma}

\pf

Item \ref{conformal modulus 1} is well known. Given an open cylinder $A$, there is a unique function $f$ with the required boundary conditions which minimizes $\int_{A}\abs {df}^{2}$, and this function is harmonic because of Dirichlet's principle. It follows that so long as $R<\infty$ there is a bijective holomorphic map $h: A\longrightarrow\{1<\abs z<e^{c}\}$ so that $\abs h=e^{cf}$. Therefore, any open annulus with finite conformal modulus is isomorphic to such an annulus. The conformal modulus of $\{1<\abs z<e^{R}\}$ is at least $R$ because of the function $f=\frac 1R\log\abs z$, and at most $R$ because of the Cauchy-Schwarz inequality. If $R(A)=\infty$, then the Uniformization Theorem implies that $A$ is conformally equivalent to a semi-infinite cylinder $S^{1}\times (0,\infty)$ or an infinite cylinder $S^{1}\times (-\infty,\infty)$. 

Item \ref{conformal modulus 2} follows from the case when $A$ contains two disjoint annuli.  Let $f$ be a  continuous function satisfying the boundary conditions which is constant outside these annuli, harmonic inside these annuli, and equal to $R(A_{1})/(R(A_{1})+R(A_{2}))$ between these annuli, where $A_{1}$ is the annuli with a boundary on which $f$ is $0$. This implies that  restricted to  $A_{i}$, $f$   is equal to a constant plus $R(A_{i})/(R(A_{1})+R(A_{2}))$ times the harmonic function which realizes the supremum in the definition of $R(A_{i})$. Then 
\[\begin{split}\int_{A}\abs{df}^{2}&=\int_{A_{1}}\abs {df}^{2}+\int_{A_{2}}\abs {df^{2}}
\\&= \lrb{\frac{R(A_{1})}{R(A_{1})+R(A_{2})}}^{2}\frac {2\pi}{R(A_{1})}+ \lrb{\frac{R(A_{2})}{R(A_{1})+R(A_{2})}}^{2}\frac {2\pi}{R(A_{2})}
\\&= \frac {2\pi}{R(A_{1})+R(A_{2})}\end{split}\] 
Therefore $R(A)\geq R(A_{1})+R(A_{2})$ as required.

 To prove item \ref {conformal modulus 3}, first set $R=R(A_2)$ and put coordinates  $(0,R)\times\mathbb R/2\pi\mathbb Z$ on $A_2$. Consider any function $f\in S(A_1)$. Without losing generality, we can assume that the shared boundary is $(0,\theta)$, and that some segment of the other boundary of $A_1$ is a curve in $A_2$ between $(\frac R2,0)$ and $(R,\theta_0)$ where $\theta_0\geq 0$. Then consider integrating $\abs {df}^2$ along diagonal lines $( Rt, -{4\pi}t+c)$. Because $\theta_{0}\geq 0$, each of these lines contains a segment inside $A_1$ on which the integral of $df$ is $1$. The length of each segment is bounded by $((4\pi)^2+R^2)^{\frac 12}$. The integral of $\abs {df}^2$ along this segment is therefore at least $(16\pi^2+R^2)^{-\frac 12}$. This tells us that the integral of $\abs {df}^2$ over $A_1$ is at least $\frac {2\pi R}{((4\pi)^2+R^2)}$, and therefore, 
\[R(A_1)\leq \frac{(16\pi^2+R^2)}{R}=R(A_2)+\frac{16\pi^2}{R(A_2)}\] 

To prove item \ref{conformal modulus 4}, consider a function $f\in S(A_1)$. We shall integrate $\abs{df}^2$ along segments  of the form $(c,t)$ inside $A_2\cap A_1$ traveling from one boundary of $A_1$ to the other. The integral of $df$ along such a segment is $1$, and the length of the segment is at most $2\pi$, so the integral of $\abs{df}^2$ along the segment is at least $\frac 1{2\pi}$, and the integral of $\abs{df}^2$ over $A_1$ is at least $\frac  {R(A_2)}{2\pi}$. Therefore, we have
\[R(A_1)\leq\frac {4\pi^2}{R(A_2)}\] 

\stop

\

The following proposition gives the existence of a decomposition of any holomorphic curve $f$ with bounded energy, local area,  and topology. This decomposition is into annular regions with small energy and complementary regions with bounded conformal modulus on which the derivative of $f$ is bounded. 
 
 \begin{prop}\label{decomposition proposition}
\begin{description}\item[]\item[] Given \begin{description}\item[]$(\ex B,J,\omega,g)$,\item an energy bound $E$,\item a local area bound $\A$, \item[] a number $N$ to bound genus and number of punctures,\item[] and small enough $\epsilon>0$ (which will bound the energy of some annular regions),\end{description}\item[] there exists\begin{description}\item a number bound $M$ depending lower semicontinuously on $(\ex B,J,\omega,g,E,\A,N,\epsilon)$ (M shall bound a number of annular regions)\item for any large enough distance $R$, a derivative bound $c$  and conformal bound $\hat R$ depending continuously on $(\ex B,J,\omega,g,E,\A,N,\epsilon, R)$
\item and an energy lower bound $\epsilon_{0}>0$ depending continuously on $(E,R,\epsilon)$\end{description}\item[] so that the following is true:
 
 Given any complete, stable holomorphic curve \[f:\ex C\longrightarrow \ex B\] with energy at most $E$, local area bounded by $\A$, and with genus and number of punctures at most $N$, there exists some collection of at most $M$ exploded annuli $\ex A_i\subset\ex C$, so that:
 \begin{enumerate}
 \item\label{first inductive condition} Each $\ex A_i$ has conformal modulus larger than $2R$ (so the conformal modulus of $A_{i}$ is $\log x \e {-l}$ where either $x>e^{2R}$ or $l>0$.)
 
 \item\label{non intersection} Put the standard cylindrical metric on $\ex A_i$, and use the notation $\ex A_{l,i}$ to denote the annulus consisting of all points in $\ex A_i$ with distance to the boundary at least $l$. Then  \[\ex A_{\frac R 2,i}\cap \ex A_{\frac R 2,j}=\emptyset\text{ if } i\neq j\]

\item \label{annulus energy bound} $f$ restricted to $ \ex A_i$ has $\omega$-energy less than $\epsilon$. 

\item\label{conformal bound} Each component of $\ex C-\bigcup \ex A_{R,i}$ is a smooth Riemann surface with bounded conformal geometry in the sense that any annulus inside one of these smooth components with conformal modulus greater than 
$\hat R$ must bound a smooth disk inside that component.

\item\label{proposition derivative bound} The following metrics can be put on $\ex A_i$ and each smooth component of $\ex C-\bigcup \ex A_{R,i}$:
\begin{enumerate}
\item On any component which is a smooth torus, use the unique flat metric in the conformal class of the complex structure so that the area of the torus is $1$. 

\item \label{disk condition}If the component is equal to a disk,  a conformal identification with the standard unit disk can be chosen so that if $\ex A_{R,i}$ is the bounding annulus, $0$ is in the complement of $\ex A_{i}$. Give components such as this the standard Euclidean metric.
\item If the component is equal to some annulus, give it the standard cylindrical metric on $\mathbb R/\mathbb Z\times(0,l)$. Give each $\ex A_i$ the analogous  standard metric.
\item Any component not equal to a torus, annulus, or  disk will admit a unique metric in the correct conformal class with curvature $-1$ so that boundary components are geodesic (there will be no components which are smooth spheres). Give these components this metric.

\end{enumerate}

On any component of $\ex C-\bigcup \ex A_{\frac {9R}{10},i}$, the derivative in the above metric is bounded by $c$ 

\[\abs {df}_g<c\]

Moreover, on $\ex A_{\frac {6R}{10},i}-\ex A_{\frac{9R}{10},i}$, the ratio between the flat metric from $\ex A_{i}$ and the metric from $\ex C-\bigcup \ex A_{R,i}$ is bounded between $c$ and $\frac 1c$.

\item \label{decomposition stability}

 $f$ restricted to any  component of $\ex C-\bigcup \ex A_{{R},i}$ which is a disk, annulus or torus has $\omega$-energy greater than $\epsilon_0>0$.
 \end{enumerate}
 \end{description}\end{prop}
 \includegraphics{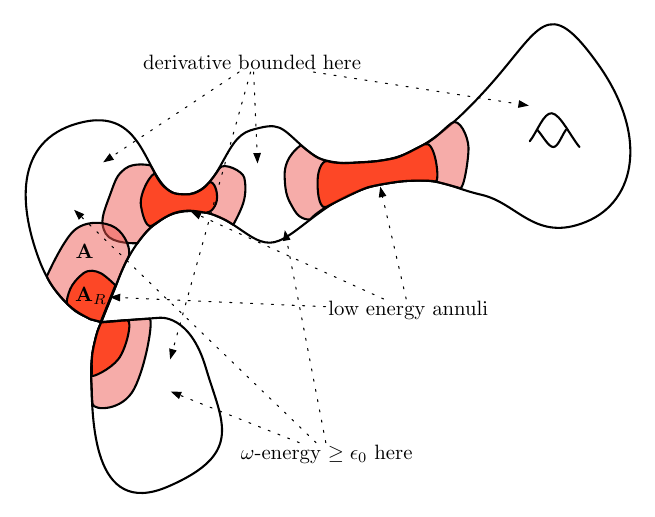}
 
  \pf

If $\ex C=\ex T$, then the proposition holds with the single annulus $\ex T$, as the curve must have no $\omega$-energy. We may therefore reduce to the case that $\ex C\neq \ex T$.
 
Lemma \ref{strict derivative bound} tells us that for $\epsilon$ small enough, any holomorphic map of the unit disk with $\omega$-energy less than $\epsilon$ and local area bounded by $\A$ must have derivative at $0$ bounded by $c_0$. We shall prove our theorem for $\epsilon$ small enough so that this is true, and small enough that Corollary \ref{cylinder energy decay} also holds with an energy bound of $\epsilon$. (This energy bound is referred to as $E$ in the notation of Corollary \ref{cylinder energy decay}). We shall also choose our distance $R$ greater than $4\pi$ (for use with Lemma \ref{annulus lemma}), and large enough that Corollary \ref{cylinder energy decay} tells us that if $f$ restricted to some smooth annulus $\ex A_i$ has $\omega$-energy less than $\epsilon$, then the $\omega$-energy of $f$ restricted to the smaller annulus $\ex A_{\frac R 2,i}$ is less than $\frac \epsilon 5$. 

We shall construct annuli in a number of stages, in particular, we shall
\begin{enumerate}
\item choose annuli around the edges of $\ex C$,
\item construct an annulus in the case that $\ex C$ is a smooth sphere,
\item construct annuli around points where the derivative of $f$ at the center of a disk is large--- this will allow us to prove the required derivative bound for $f$ later on,
\item construct an annulus inside any long cylinder in what remains--- this will allow us to prove our conformal bound on the complement of the annuli,
 \item and finally,  merge pairs of annuli which bound annular regions with low $\omega$-energy--- this will allow us to bound the number of constructed annuli.
\end{enumerate}

{\noindent\large{Stage 1}:}

Note that any edge of $\ex C$ has zero $\omega$-energy, so we can choose an exploded annulus $\ex A_i$ containing each edge with $\omega$-energy less than $ \frac \epsilon 5$.  Note, for use with item \ref{decomposition stability}, that in the case that this annulus  bounds a disk, the $\omega$-energy of the resulting disk will be at least $\frac{4\epsilon} 5$ because   Lemma \ref{strict derivative bound}  implies that each spherical component of the smooth part of $\ex C$ with at most one puncture must have $\omega$-energy at least $\epsilon$. (If a spherical component with at most one puncture had $\omega$-energy less than $\epsilon$, we could use an injective map of a disk to our sphere with arbitrarily large derivative at $0$ to contradict Lemma \ref{strict derivative bound}.)

 We can choose these $\ex A_i$ to be  mutually disjoint. Note that the complement of these $\ex A_i$ is a smooth Riemann surface with boundary. 

\

{\noindent\large{Stage 2}:}

 If  $\ex C$ is a smooth sphere, we shall now remove an annulus. As argued above, our sphere must have energy at least $\epsilon$ in order to be stable. We can therefore put some round metric in the conformal class determined by the complex structure so that there exist $3$ mutually perpendicular geodesics which divide our sphere into $8$ regions so that two antipodal regions each have  $\omega$-energy at least $\frac \epsilon8$.  In particular, for any metric, we may choose one geodesic that divides the sphere into two hemispheres with  $\omega$-energy at least $\frac \epsilon 2$, then choose another perpendicular geodesic that divides one of these hemispheres into two regions of $\omega$-energy at least $\frac \epsilon 4$. At least one of these regions will have an antipodal partner region with $\omega$-energy at least $\frac \epsilon 4$. We may then choose a circle perpendicular to our two geodesics which divides this region into two regions with $\omega$-energy at least $\frac \epsilon 8$. Again, one of the antipodal partner regions will have $\omega$-energy at least $\frac \epsilon 8$. We may then choose our metric so that our circle becomes a geodesic and our first two original geodesics remain geodesics. 
 
 By choosing $\epsilon_0>0$ small enough, we can then get that there exist two disks with $\omega$-energy at least $\epsilon_0$ which intersect each of these two antipodal regions, and which have radius as small as we like. Choose $\epsilon_0>0$ small enough that the complement of these disks  is an annulus of conformal radius at least ${k (2R+1)}$ for some integer $k>\frac {5E}\epsilon$. How small $\epsilon_0$  is required to achieve this depends only on $R$, $\epsilon$ and $E$. Then we can divide this annulus into $k$ annuli with conformal modulus at least $(2R+1)$; at least one of these $k$ annuli must have $\omega$-energy at most  $\frac \epsilon 5$. Add this annulus to our collection. Note that it bounds disks which have energy at least $\epsilon_0$.   
 
 \
 
 \includegraphics{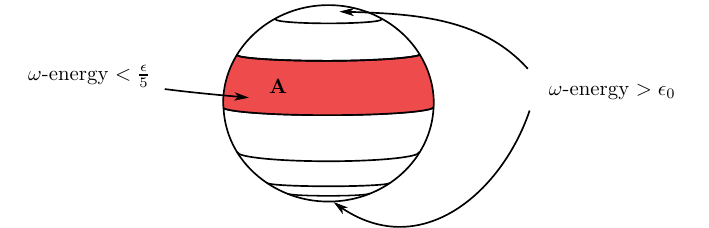}
 
 \
 
 {\noindent\large{Stage 3}:}

Consider a holomorphic injection of the unit disk $\iota:D\longrightarrow\ex C$ into the complement of the union of our annuli $\ex A_{R,i}$ constructed in stages $1$ and $2$. Suppose that $\iota(0)$ is in the complement of $\ex A_{\frac R4,i}$. Lemma \ref{annulus lemma} part \ref{conformal modulus 2} implies that $\{e^{-3R/4}<z<1\}$ is not contained within $\ex A_{\frac R4,i}$. Therefore, each circle homotopic to the boundary within the annulus $\{e^{-\frac {16\pi^2}R-3R/4}<z<e^{-3R/4}\}$ must not be contained entirely within  $\ex A_{\frac R 4,i}$. If such a circle intersects $\ex A_{\frac R2, i}$, it must intersect both boundaries of a component of $\ex A_{\frac R4,i}- \ex A_{\frac R2,i}$. Lemma \ref{annulus lemma} item \ref{conformal modulus 4} tells us that at least one of these circles must not intersect both boundaries of a component of $\ex A_{\frac R4,i}- \ex A_{\frac R2,i}$, so at least one of these circles is in the complement of $\ex A_{\frac R2,i}$. Therefore the restriction of $\iota$ to $\abs z<e^{-\frac {16\pi^2}R-3R/4}$ is in the complement of $\ex A_{\frac R2,i}$.

 Now suppose that at $z=0$, $\abs {d(f\circ \iota)}>c_0e^{\frac {16\pi^2}R+(k+1)(2R+1)}$ where $k$ is some integer greater than $\frac {5E}\epsilon$. Then the restriction of $f$ to the disk $\abs z< e^{-\frac {16\pi^2}R-(k+1)(2R+1)}$ must have energy greater than $\epsilon_{0}$ due to Lemma \ref{strict derivative bound}. This disk is surrounded by $k$ disjoint annuli of conformal modulus $(2R+1)$ which are contained in the complement of $\ex A_{\frac R2,i}$. At least one of these must have energy less than $\frac \epsilon 5$. Add this annulus to our collection.  Continue adding annuli in this manner.

\includegraphics{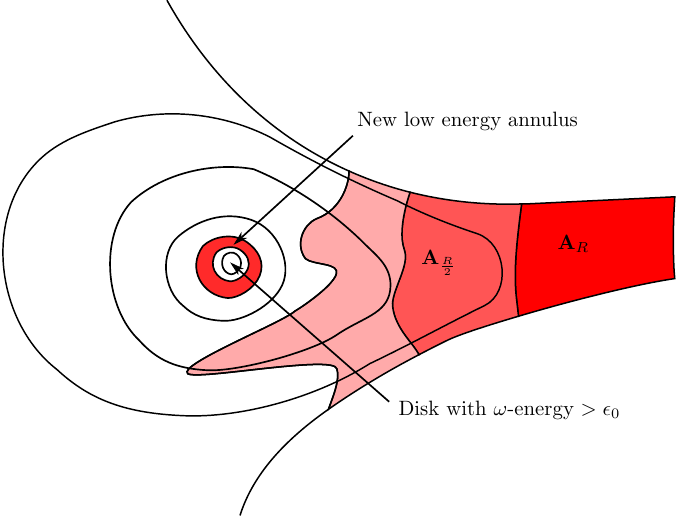}

 After a finite (but not universally bounded) number of times, no more annuli can be added in this way. In particular, each new annulus added this way bounds a new disk with energy at least $\epsilon_{0}$. At most $E/\epsilon_{0}$ such disks exist, so if there are an infinite number of such new annuli, an infinite number of such annuli must be homotopic to the boundary of a single disk. If there are $n$ such annuli homotopic to the boundary of such a disk,  this disk must contain an annulus of  conformal modulus $nR$ bounding an inclusion $\iota:D\longrightarrow \ex C$ so that the derivative of $\iota\circ f$ is universally bounded below.  Then,  $\iota$ may be rescaled by a factor of $e^{nR}$ and still give an inclusion of the unit disk into $C$  with image disjoint from all annuli constructed in stages 1 and 2. If there were an infinite number of our new annuli homotopic to the boundary of a single disk,  the derivative of $f$ would therefore be forced to be unbounded, which is impossible.

\

{\noindent\large{Stage 4}:}

  We now want to add annuli so that condition \ref{conformal bound} is satisfied. Suppose  some smooth component of $\ex C-\bigcup  \ex A_{R,i}$ contains some annulus $\ex A$ of conformal modulus greater than $2R(3+\frac{5E}{\epsilon})$ that does not bound a disk. Note that as we have chosen $R>4\pi$, the size of such an annulus is greater than $2R(1+\frac{5E}{\epsilon})+2R+\frac{32\pi^2}R$. Because we have chosen $\ex A$ in the complement of $\ex A_{R,i}$, Lemma \ref{annulus lemma} part \ref{conformal modulus 2} implies the following: giving $\ex A$ the standard flat metric, the annulus consisting of points of distance at most $R$ from its end can't be contained entirely within $\ex A_i$. Therefore, no circle homotopic to the boundary in $\ex A_{i}$ is  contained within the annulus of points a distance at least $R$ into the interior of $\ex A$.  Applying Lemma \ref{annulus lemma} part \ref{conformal modulus 4}  with the annulus consisting of points in $\ex A$ of distance between $R$ and   $R+\frac {16\pi^2} R$ from this end in the place of $A_{1}$, and a component of $\ex A_{i}-\ex A_{\frac R2,i}$ in place of $A_{2}$, we see that at least one circle homotopic to the boundary within $\ex A_{i}-\ex A_{\frac R2,i}$ must not intersect both boundaries of the other annulus. Therefore our annulus $\ex A$ minus annuli at each end of conformal modulus $R+\frac {16\pi^2} R$ must be contained entirely in the complement of $\ex A_{\frac R2,i}$. We then obtain more than $\frac {5E}{\epsilon}$ disjoint annuli in the compliment of $\ex A_{\frac R2,i}$ with conformal modulus greater than $2R$. The restriction of $f$ to at least one of these annuli must have $\omega$-energy less than $\frac\epsilon 5$.   Add this annulus to our collection.
  
  \includegraphics{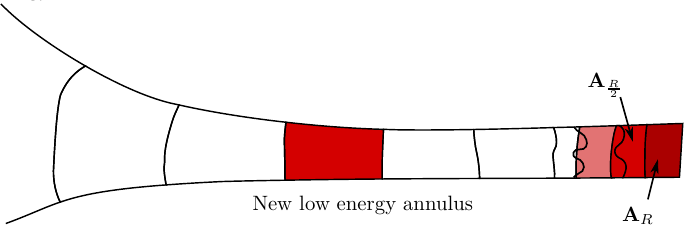}

 \

   We shall now argue that we can only add a finite (but not universally bounded) number of annuli in the above manner:
%
  There are a  finite number of exploded annuli with infinite conformal modulus which we chose in stage 1 to cover the edges of our curve. The number of these is bounded by our topological bounds, and the bound on the number of  spherical components of the smooth part of $\ex C$ with fewer than $3$ punctures from Lemma \ref{omega monotonicity}. If $\ex C$ was equal to a smooth sphere, we then added a single annulus in stage 2. Then in stage 3 we added a finite number of annuli to bound the derivative of $f$. 
  
   Call the complement of all the annuli from stages 1, 2 and 3, $\ex C_0$. Each annulus added to $\ex C_{0}$ in stage 4 is in a nontrivial homotopy class of annuli in $\ex C_{0}$. There is  a bound on the number of homotopy classes in $\ex C_0$ of our remaining annuli, given by our bounds on the topology of $\ex C$, the fact that we have removed a finite number of annuli, and the observation that for any Riemann surface with boundary, there are only a finite number of homotopy classes which contain annuli of conformal modulus greater than $4\pi$. (Lemma \ref{annulus lemma} part \ref{conformal modulus 4} implies that if two annuli have conformal modulus greater than $4\pi$, then the intersection number of a boundary circle of one with a boundary circle of the other must be $0$. This bounds the number of homotopy classes containing an annulus by the largest number of simple closed curves so that no two curves are in the same homotopy class and so that the  intersection number of any two curves is $0$.)
   
    There is also a finite number of annuli from Stage 4 in each homotopy class. If there were an infinite number of annuli $\ex A_{\frac R2,i}$ in our collection in the same homotopy class, then there would exist an annulus in $\ex C_0$ in that homotopy class which contains all of them. As the $\ex A_{\frac R2,i}$ are disjoint and have conformal modulus at least $R$, this annulus would have to have infinite conformal modulus, which is impossible because $\ex C_{0}$ is a compact, smooth Riemann surface with boundary, and this annulus is not in the trivial homotopy class of annuli in $\ex C_{0}$.

\

{\noindent\large{Stage 5}:}

Now we shall merge some annuli so that there exists no annular component of $\ex C-\bigcup \ex A_{R,i}$ with $\omega$-energy less than $\frac\epsilon 5$. Then we will have a bound on the number of annuli which is independent of $f$. Suppose that we have some collection $\{\ex A_1,\dotsc,\ex A_n\}$ of our annuli so that for each $i$, $\ex A_{R,i}$ and $\ex A_{R,i+1}$ bound an annulus which has $\omega$-energy less than $\frac \epsilon 5$. Use the notation $\ex A_{[m, n]}$ to denote the annulus that consists of $\ex A_m$, $\ex A_n$ and everything in between. Then, as we've chosen each of our $\ex A_i$ to have $\omega$-energy less than $\frac\epsilon 5$, $\ex A_{[i,i+2]}$ has $\omega$-energy less than $\epsilon$, so we can apply Corollary \ref{cylinder energy decay} to show that far enough into the interior of this cylinder, there is very little $\omega$-energy. (Note that  the edges of exploded curves  always have no $\omega$-energy, so there is no difficulty in applying this corollary in the seemingly more general setting of exploded annuli.) Apply Lemma \ref{annulus lemma} part \ref{conformal modulus 3}, with $\ex A_{i}$ or $\ex A_{i+2}$ in the place of $A_{1}$ and the appropriate component of $\ex A_{[i,i+2]}-\ex A_{R,[i,i+2]}$ in the place of $A_{2}$. As we have chosen $R>4\pi$, we have $R+\frac{16\pi^2}R<2R$ which is less than the conformal modulus of $\ex A_{i}$ or $\ex A_{i+2}$. It follows that $\ex A_{[i,i+2]}-\ex A_i-\ex A_{i+2}\subset \ex A_{\frac R2, [i,i+2]}$. We have chosen $R$ large enough that Corollary \ref{cylinder energy decay} tells us that the energy of $f$ restricted to $\ex A_{\frac R2,[i,i+2]}$ is less than $\frac \epsilon 5$.
It follows that the energy of $A_{[i,i+3]}$ is less than $\epsilon$, so we can now repeat this argument inductively to show that the energy of $f$ restricted to $\ex A_{[1,n]}$ is less than $\epsilon$ and $\ex A_{\frac R2, [1,n]}$ is less than $\frac\epsilon 5$. 

\includegraphics{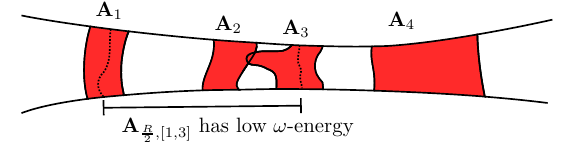}

\

 Now, replace all of our sets of annuli $\{\ex A_1,\dotsc, \ex A_n\}$ of maximal size obeying the above conditions with the annulus given by $\ex A_{\frac R2,1}$, $\ex A_{\frac R2,n}$ and everything in between. As this annulus is contained inside $\ex A_{[1,n]}$, it has energy less than $\epsilon$. (Note that if we had some collection of annuli which bounded small energy annuli in a cyclic fashion, we would get that  $\ex C$ was a torus with energy less than $\epsilon$, and was therefore not stable.)

 \

The construction of our annuli is now complete. We must now check that all the conditions we require are satisfied. First, the number of annuli is bounded independent of $f$. In particular, there are at most $E/\epsilon_{0}$ components of $\ex C$ minus our annuli which are disks, and at most $5E/\epsilon$ components which are annuli. The number of other components is therefore bounded for topological reasons, as the genus and number of punctures of $\ex C$ is bounded.
 
  The size of each annulus is greater than $2R$, so condition \ref{first inductive condition} holds.
 The set of  annuli constructed in stages 1 to 4 obeyed the non-intersection condition. Therefore, the merged annuli constructed in stage 5  also obey the non-intersection condition \ref{non intersection}. In each stage above, we always constructed our annuli with $\omega$-energy less that $\epsilon$, so condition \ref{annulus energy bound} holds. 
Our  set of annuli constructed by the end of Stage 4 obeyed condition \ref{conformal bound} with a bound of $2R(3+\frac{5E}\epsilon)$. After stage 5, each component of the complement of $\ex A_{R,i}$ is given by adding an annulus of conformal modulus  $R$ to each boundary of some component of the complement of $\ex A_{\frac R2,i}$ at stage $4$.  Lemma \ref{annulus lemma} parts \ref{conformal modulus 2} and \ref{conformal modulus 4}   then imply that by increasing our conformal bound by $2R +16\pi^{2}/R$,  we can achieve condition \ref{conformal bound} after Stage 5. A lower bound  of $\epsilon_{0}$ on the $\omega$-energy of unstable components is also satisfied by  construction, so condition \ref{decomposition stability} holds.

\

All that remains is condition \ref{proposition derivative bound}. Note first that the  complement of our annuli $\ex A_{R,i}$ admits metrics as described in condition \ref{proposition derivative bound}. To show the derivative bounds from condition \ref{proposition derivative bound}, we shall show below that if $\iota$ is an injective map of a holomorphic disk into either $\ex A_{i}$ or the complement of $\ex A_{R,i}$, then $d\iota$ and  $d(f\circ\iota)$ are bounded at $0$, and we shall construct such maps $\iota$ with derivatives which are bounded below  at $0$.

\begin{claim}\label{c_{1}}There exists a constant $c_{1}$ so that any holomorphic injection $\iota$ of a unit disk into either $\ex A_{i}$ or the complement of $\ex A_{R,i}$ has $d\iota$ bounded by $c_{1}$ at $0$ when the metrics from condition \ref{proposition derivative bound} are used. 
\end{claim}

We shall prove Claim \ref{c_{1}} using a bubbling argument. If this claim failed for maps to $\ex A_{i}$,  a standard bubbling argument, used as in the proof of Lemma \ref{strict derivative bound}, would produce a holomorphic map of $\mathbb C$ to $\mathbb C^{*}$ with nonzero derivative at $0$ which extends to a map  $\mathbb CP^{1}\longrightarrow \mathbb CP^{1}$ of degree at most $1$.  No such map $\mathbb C\longrightarrow \mathbb C^{*}$ exists, therefore Claim \ref{c_{1}} holds for maps to $\ex A_{i}$. 

Now to prove Claim \ref{c_{1}} for maps to the complement of $\ex A_{R,i}$. The topology of each component of the complement of $\ex A_{R,i}$ is bounded by the bound on the genus and number of punctures of $\ex C$ and the bound on the number of annuli $\ex A_{i}$. The conformal structure of each component is also bounded because of condition \ref{conformal bound}. It follows that there exists a smooth family $\hat{\ex C}$ of Riemann surfaces with boundary (and with metrics), so that some compactly contained subfamily of $\hat{\ex C}$ contains  all possible Riemann surfaces with boundary that may occur in the complement of $\ex A_{R,i}$. Therefore, if Claim \ref{c_{1}} failed, a standard bubbling argument gives that there exists a non-constant, bounded area, holomorphic map  of $\mathbb C$ into one of the surfaces in our family $\hat{\ex C}$. As this fiber must be contained in a compact Riemann surface, Lemma \ref{cylinder bound} and the usual removable singularity theorem for holomorphic maps may be used to extend our map to a non-constant holomorphic map of $\mathbb CP^{1}$  to our fiber of $\hat {\ex C}$. We may choose $\hat{\ex C}$ to not contain any fibers equal to $\mathbb CP^{1}$, so no such map exists. This completes the proof of Claim \ref{c_{1}}.

\

As $R>4\pi>0$, there exists a disk of unit size centered on any point inside $\ex A_{\frac {6R}{10},i}-\ex A_{\frac {9R}{10},i}$. The derivative at the center of this disk is $1$ in the metric on $\ex A_{i}$, 
and is bounded by  $c_{1}$ in the metric on the complement of $\ex A_{R,i}$. Therefore, inside $\ex A_{\frac {6R}{10},i}-\ex A_{\frac {9R}{10},i}$,  the ratio of the metric on the complement of $\ex A_{R,i}$ divided by the metric on  $\ex A_{i}$ is  bounded by $c_{1}$.

\begin{claim}\label{c_{2}}So long as $R$ is chosen large enough, there exists a constant $c_{2}>0$ so that for  any point $p$ in the complement of $\ex A_{\frac{9R}{10},i}$, there exists an injective map $\iota$ of a disk into the complement of $\ex A_{R,i}$ so that $\iota(0)=p$ and so that  $\abs{d\iota}$ is larger than $c_{2}$ at $0$. If $p$ is in $\ex A_{\frac {6R}{10},i}$, then then $\iota$ may be chosen to have image contained in $\ex A_{i}$.
\end{claim}

As in the proof of Claim \ref{c_{1}}, we shall use that there exists a family  $\hat{\ex C}$ of Riemann surfaces with boundary compactly containing all the possible surfaces which may occur in the complement of $\ex A_{R,i}$ with the metrics from condition \ref{proposition derivative bound}. The intersection of $\ex A_{\frac {9R}{10},i}$ with a component of the complement of $\ex A_{R,i}$ is one or two annuli of conformal modulus $\frac R{10}$ surrounding boundary components. In order to prove Claim \ref{c_{2}} with a compactness argument involving $\hat{\ex C}$, we must compare these annuli to standard collar neighborhoods of boundary components defined using the metric from condition \ref{proposition derivative bound}.

Property \ref{non intersection} implies that each boundary of the  complement of $\ex A_{R,i}$ is surrounded by an annulus of conformal modulus $R/2$. We are free to choose $R$ large enough that this implies that there exists a `collar' annulus of conformal modulus $4\pi$ consisting of the set of points within some fixed distance of the boundary. Using Lemma \ref{annulus lemma} part \ref{conformal modulus 3} with this collar in the place of $A_{2}$ gives that any annulus surrounding a boundary component with conformal modulus greater than $8\pi$ must contain the collar of conformal modulus $2\pi$ consisting of points within some fixed distance of the boundary. We are free to choose $R>80\pi$, so $\ex A_{\frac {9R}{10},i}$ contains this collar. 

We have now proved that there is some neighborhood $U$ of the boundary in $\hat{\ex C}$ which is contained in the image of $\ex A_{\frac {9R}{10},i}$ intersected with the complement of $\ex A_{R,i}$. Let $s:\hat{\ex C}\longrightarrow[0,\infty]$ indicate the map so that $s(p)$ is the supremum of the possible derivatives of injective holomorphic maps of disks into a fiber of $\hat{\ex C}$ centered at $p$. The function $s(p)$ is continuous, and positive on the interior of $\hat{\ex C}$. It follows that $s(p)$ is bounded below when restricted to the complement of $U$ and the compact subfamily  of $\hat{\ex C}$ which contains all Riemann surfaces with boundary that may occur in the complement of $\ex A_{R,i}$. 

Therefore, there exists some $M>0$ so that given any point $p$ in the complement of $\ex A_{\frac {9R}{10},i}$, there exists some $\iota$  so that $\iota(0)=p$ and so that  $\abs{d\iota}$ is larger than $M$ at $0$. If $p$ is contained in $\ex A_{\frac {6R}{10},i}$, then Lemma \ref{annulus lemma} part \ref{conformal modulus 4} implies that the image of $\iota$ restricted to $\{\abs z<e^{-8\pi^{2}/R}\}$ is contained entirely within $\ex A_{i}$. (The two annuli used for this application of Lemma \ref{annulus lemma} consist of a component of $\ex A_{i}-\ex A_{\frac R2,i}$ and $\{1>\abs z>e^{-8\pi^{2}/R}\}$.) Claim \ref{c_{2}} therefore holds with $c_{2}=Me^{-8\pi^{2}/R}$.

\

Claims \ref{c_{2}} and \ref{c_{1}} imply that the ratio of the metric on $\ex A_{i}$ divided by the metric on the complement of $\ex A_{R,i}$ is bounded by $c_{1}/c_{2}$ on $\ex A_{\frac{6R}{10},i}-\ex A_{\frac{9R}{10},i}$. 

It remains to prove that the derivative of $f$ is bounded on the complement of $\ex A_{\frac {9R}{10},i}$. To achieve this, we shall use Claim \ref{c_{2}} with Claim \ref{c'} below which bounds the derivative of $f\circ \iota$ at $0$.

\begin{claim}\label{c'} There exists a constant $c'$ so that any holomorphic injection $\iota$ of a unit disk into either $\ex A_{i}$ or the complement of $\ex A_{R,i}$ has 
\[\abs{d(f\circ \iota)}_{g}\]
bounded by $c'$ at $0$.

\end{claim}

As the $\omega$-energy of $\ex A_{i}$ is less than $\epsilon$, Lemma \ref{strict derivative bound} implies Claim \ref{c'} in the case  that the target of $\iota$ is $\ex A_{i}$. Similarly, suppose that $\iota(0)$ is within $\ex A_{\frac R4,i'}$ for some annulus $\ex A_{i'}$ constructed by the end of Stage 4.  Then Lemma \ref{annulus lemma} part \ref{conformal modulus 4} tells us that the restriction of $\iota$ to the disk of radius $e^{-16\frac{\pi^2}R}$ must be contained inside $\ex A_{i'}$, and therefore Lemma \ref{strict derivative bound} applied to this disk gives Claim \ref{c'} for such an $\iota$.

It remains to prove Claim \ref{c'} for an $\iota$ with target in the complement of the  $\ex A_{R,i}$ constructed by the end of Stage 5  so that $\iota(0)$ is in the complement of $\ex A_{\frac R4,i'}$ for all annuli $\ex A_{i'}$ constructed by the end of Stage 4. If $\ex A_{i'}$ is an annulus which survives in Stage 5, then $\iota$ has image in the complement of $\ex A_{R,i'}$. On the other hand, if $\ex A_{i'}$ is an annulus which is absorbed into a larger annulus $\ex A_{i}$ during Stage 5, we must shrink the domain of $\iota$ to ensure that it has image in the complement of $\ex A_{R,i'}$. The only case of concern is when one boundary of $\ex A_{i}$ coincides with a boundary of $\ex A_{\frac R2,i'}$, (recall that in stage 5, $\ex A_{i}$ is constructed by taking the union of some number of the $\ex A_{\frac R2,i'}$ with everything in between them.) In this case, $\iota(0)$ is in the complement of such an $\ex A_{i}$, therefore Lemma \ref{annulus lemma} part \ref{conformal modulus 2} implies that the image of $\{e^{-R}<\abs z<1\}$ is not contained entirely within $\ex A_{i}$, and is therefore not contained entirely within $\ex A_{\frac R2,i'}$. Lemma \ref{annulus lemma} part \ref{conformal modulus 4} then implies that the image of $\{\abs z<e^{-R-8\pi^{2}/R}\}$ is in the complement of $\ex A_{R,i'}$. Stage 3 terminated when there were no more disks in the complement of $\ex A_{R,i'}$ centered in the complement of $\ex A_{\frac R4,i'}$ and derivative at $0$ greater than $c_0e^{\frac {16\pi^2}R+(k+1)(2R+1)}$, so $\abs{d(f\circ \iota)}$ at $0$ is bounded by $c_0e^{\frac {24\pi^2}R+(k+1)(2R+1)+R}$. This completes the proof of Claim \ref{c'}

\

Now suppose that $p$ is any point in the complement of $\ex A_{\frac{9R}{10},i}$. Claim \ref{c_{2}} implies that there exists a injective holomorphic map $\iota$ of a unit disk into the complement of $\ex A_{R,i}$ so that $\abs{d\iota}$ at $0$ is at least $c_{2}$. Claim \ref{c'} tells us that  $\abs{df\circ \iota}_{g}$ at $0$ is at most $c'$. It follows that $\abs{df}_{g}$ at $p$ is at most $c'/c_{2}$. This completes the proof of the derivative bounds from condition \ref{proposition derivative bound}, where $c$ is taken as the maximum of $c'/c_{2}$, $c_{1}/c_{2}$, and $c_{1}$.

   \stop

 \section{Compactness}\label{compactness section}
 
 The following theorem is a version of Gromov compactness for curves in exploded manifolds. It requires the extra assumption that curves have bounded local area in the sense of definition \ref{lab}. This local area bound follows from topological bounds under extra assumptions on $\ex B$ described in section \ref{local area bound section}.

 \begin{thm}\label{completeness theorem}
 Let $\hat{\ex B}\longrightarrow \ex G$ be a family of basic complete exploded manifolds with a metric $g$, and a family of  $\dbar\log$ compatible almost complex structures $J$ tamed by a family of taming  forms $\omega$.

 \
 
  \noindent Given:
  \begin{description}
  \item any convergent sequence of points in $\ex G$ \item and a lift of this sequence to a sequence of $J$-holomorphic curves in $\hat{\ex B}\longrightarrow \ex G$ in the fibers over these points
  \begin{description} \item with bounded $\omega$-energy, local area, genus, and number of punctures,
   \end{description}
   \end{description} 
   
  \noindent  there exists a subsequence of curves $f^i$ which converges to a holomorphic curve $f$ as follows:  
  
  \
  
  There exists a
  sequence of  families of $\C\infty1$ curves,
  
  \[\begin{array}{ccc} (\ex {\hat C},j_i) & \xrightarrow{\hat f^i} &\hat{\ex B}
  \\ \downarrow & &\downarrow
  \\ \ex F& \longrightarrow & \ex G\end{array}\]
 
  so that this sequence of families converges in $\C\infty1$ to the family
  
  \[\begin{array}{ccc} (\ex {\hat C},j) & \xrightarrow{\hat f} &\hat{\ex B}
  \\ \downarrow & &\downarrow
  \\ \ex F& \longrightarrow & \ex G\end{array}\]

  and a sequence of points $p^i$ in $\ex F$ so that \begin{itemize}\item $p^{i}\rightarrow{ p}$,\item $f^i$ is the map given by the restriction of $\hat f^i$ to the fiber over $p^i$,\item and $f$ is given by the restriction of $f$ to the fiber over $p$.\end{itemize} 
  
  \end{thm}
  
  \
  
  Although it is not obvious, the notion of convergence used in Theorem \ref{completeness theorem} is the natural notion of convergence  on the moduli stack of (not necessarily holomorphic) $\C\infty1$ curves. On the moduli stack of $\C\infty1$ curves, we can define an open substack to be a substack $\mathcal U$ with the property that given any $\C\infty1$ family of curves $\hat f$, the subset of $\hat f$ consisting of curves that are in $\mathcal U$ is open. It is proved in \cite{evc} that the notion of convergence used above is equivalent to convergence in this natural topology on the moduli stack of $\C\infty1$ curves.
  
 \   
  
  The proof of Theorem \ref{completeness theorem} uses  Lemma \ref{strong cylinder convergence} on page \pageref{strong cylinder convergence} and  Proposition \ref{decomposition proposition} on page \pageref{decomposition proposition}. Together, these allow us to decompose holomorphic curves into pieces with bounded behavior. A standard Arzela-Ascoli type argument gets a type of convergence of these pieces. As the topology on exploded manifolds is in general not Hausdorff, this limit will not be unique, so we must work a little more to make sure that our limiting pieces match up.

  One way to think of the problem at hand is that there will be a unique limit of our sequence of curves in the smooth part of $\ex B$, but there will be a family of valid tropical choices for our limit. We must construct a  family of curves which contains  these tropical possibilities, and extend this family a little so that we can capture the behavior of our sequence of curves. There is some choice involved in this extension of a family, but we shall construct the family locally in coordinate charts. To make our choices match up, we shall use equivariant coordinate charts, which are constructed in Lemma A3.
 
 \
 
 The following lemma constructs compatible $\ex T^{n}$ actions on neighborhoods of strata in $\ex B$. These $\ex T^{n}$ actions are defined on page \pageref{T action}. Locally, each $\ex T^{n}$ action is given by  multiplying the $\et nP$ coordinates of $\mathbb R^{m}\times \et nP$ by constants, so 
 \[(\tilde c_{1},\dotsc,\tilde c_{n})*(x,\tilde z_{1},\dotsc,\tilde z_{n})=(x,\tilde c_{1}\tilde z_{1},\dotsc,\tilde c_{n}\tilde z_{n})\]
We shall need these $\ex T^{n}$ actions to move parts of our curves around in a consistent way.

\begin{lemma}\label{strata charts}
Given any basic \exploded manifold $\hat {\ex B}$, and a family $\hat{\ex B}\longrightarrow\ex G$,  for each stratum $\hat{\ex B}_i\subset\hat{\ex B}$,  there exists some open neighborhood $\ex U_i\subset\hat{ \ex B}$ containing $\hat{\ex B}_i$ and a free action of $\ex T^{n}$ on $\ex U_{i}$ in the sense of definition A2 (where $n$ is the dimension of the tropical part of $\hat{\ex B}_{i}$). This $\ex T^{n}$ action on $\ex U_{i}$ satisfies the following properties:
\begin{enumerate}
\item The $\ex T^{n}$ action preserves fibers of the family $\hat{\ex B}\longrightarrow\ex G$, in the sense that if $p_{1}$ and $p_{2}$ have the same image in $\ex G$, then $\tilde z*p_{1}$ and $\tilde z*p_{2}$ also have the same image in $\ex G$. 
\item If  the closure of $\hat{\ex B}_i$ contains $ {\hat{\ex B}_j}$, then the restriction of the action of $\ex T^n$ on $\ex U_{i}$ to   $\ex U_i\cap\ex U_j$ is equal to the restriction of the  $\ex T^m$ action on $\ex U_j$ to $\ex U_{i}\cap \ex U_{j}$ and to some subgroup. 

\end{enumerate}

\end{lemma}

\pf

This follows immediately from Lemma A3. Choose equivariant coordinates on $\hat{\ex B}\longrightarrow \ex G$, and let $\ex U_{i}$ be some open  neighborhood of $\ex B_{i}$ contained inside the union of all our coordinate charts that intersect $\ex B_{i}$. As each of these coordinate charts intersect $\ex B_{i}$, and $\ex B$ is basic, there is a canonical inclusion of $\totb{\ex B_{i}}$ into the tropical part of each of these coordinate charts, and a corresponding action of $\ex T^{n}$ on these coordinate charts given by multiplying the corresponding coordinates by constants. The fact that our coordinates are equivariant implies that this gives a $\ex T^{n}$ action on $\ex U_{i}$. The compatibility of the actions on $\ex U_{i}\cap \ex U_{j}$ and the compatibility of the action with the projection map $\hat {\ex B}\longrightarrow \ex G$ all follow from the fact that these coordinate charts are equivariant.  

\stop

\

We are going to be translating pieces of our curve around by locally defined $\ex T^{m}$ actions. The assumption that  $\ex B$ is basic simplifies the possibilities of how these locally defined $\ex T^{m}$ actions are related to each other. So long as a local area bound for curves can be achieved, this basic assumption can be weakened, but the assumption of an immersion $\totb{\ex B}\longrightarrow \mathbb R^{n}$ used in section \ref{local area bound section} to obtain local area bounds implies that $\ex B$ is basic anyway. In the family case, by restricting to a subset of $\ex G$ which contains the image of a subsequence, we can also assume that $\hat{\ex B}$ is basic.

 \
 
 We shall now start working towards proving Theorem \ref{completeness theorem}. The proof shall rely on Lemma \ref{strong cylinder convergence} on page \pageref{strong cylinder convergence} and  Proposition \ref{decomposition proposition} on page \pageref{decomposition proposition}. We shall be using the notation from Proposition \ref{decomposition proposition}.
 First, choose the following:
 
 \begin{enumerate}
 \item Choose a finite collection of coordinate charts on ${\ex B}$, each contained in some $\ex U_i$ from Lemma \ref{strata charts} and small enough to apply Lemma \ref{strong cylinder convergence}. For the case of a family, restrict $\hat{\ex B}\longrightarrow\ex G$ to some small coordinate chart on $\ex G$  in which the image of our sequence of curves converges (after passing to a suitable subsequence), and choose our finite collection of coordinate charts  on (the smaller, renamed) $\hat {\ex B}$ satisfying the above. 
 \item 
  Choose  exploded annuli $\ex A_n^i\subset \ex C^i$ satisfying the conditions in Proposition \ref{decomposition proposition} with $\omega$-energy bound $\epsilon$ small enough and $R$ large enough that each smaller annulus $\ex A^i_{\frac {6R}{10},n}$ is contained well inside some $\ex U_j$ in the sense that it is of distance greater than $2$ to the boundary of $\ex U_{j}$. (The fact that we can achieve this follows from Lemma \ref{omega cylinder bound} on page \pageref{omega cylinder bound}.) Also choose $\epsilon$ small enough and $R$ large enough so that each connected smooth component of $\ex A^i_{\frac {6R}{10},n}$  can have Lemma \ref{strong cylinder convergence} applied to it.
\end{enumerate}

\
 
 Use the notation $C^i_m$ to indicate connected components of $\ex C^i-\bigcup_k \ex A^i_{\frac {8R}{10}, k}$. If there are no $C^{i}_{m}$, then $\ex C^{i}$ is isomorphic to  $\ex T$. The proof of Theorem \ref{completeness theorem} is trivial if there exists a subsequence of curves with domains equal to $\ex T$, so we shall assume that $\ex C^{i}-\bigcup_{k}\ex A^i_{\frac {8R}{10}, k}$ has at least one component. We can choose a subsequence so that\begin{itemize}\item the number of such components is the same for each $\ex C^i$,\item $C^i_m$ has topology that is independent of $i$,
 \item the combinatorics of which boundaries of $C^{i}_{m}$ are connected by annuli is independent of $i$
 \item and there are diffeomorphisms identifying $C^i_m$ for different $i$ so that the sequence of components $C^{i}_{m}$ converges as $i\rightarrow \infty$ to some $C_m$ in the sense that the metric from Proposition \ref{decomposition proposition} and the complex structure converges to one on $C_m$. That this may be achieved follows from the bounded geometry from items \ref{proposition derivative bound} and \ref{conformal bound} of Proposition \ref{decomposition proposition}.
 \end{itemize}
 
   \includegraphics{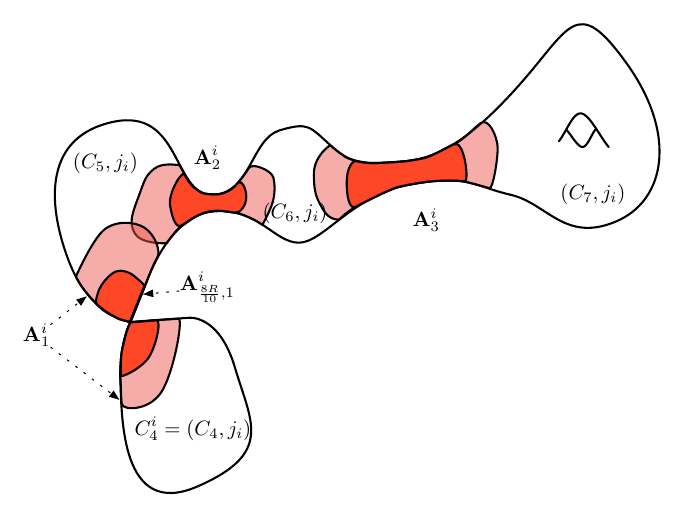}

 We shall view any  annulus $\ex A^{i}_{n}$ with finite conformal modulus as glued out of two semi infinite annuli, so $\ex A^{i}_{n}$ has two holomorphic coordinates $\tilde z^{+}$ and $\tilde z^{-}$ related by
 
 \[\tilde z^+\tilde z^-= \w^{i}_{n}\in(0,\infty)\e {[0,\infty)}\]
 where the conformal modulus of $\ex A_{n}^{i}$ is $-\log \w^{i}_{n}$.  Regard $\ex A^{i}_{n}\subset (\et 11)^{2}$ as given by 
 \[\ex A_{n}^{i}=\{\tilde z^+\tilde z^-=\w^{i}_{n}, \ \abs{\tilde z^{+}}<1,\ \abs{\tilde z^{-}}<1 \}\subset (\et 11)^{2}\]

 Choose a subsequence so that   $\{\w^{i}_{n}\}$ converges. (In other words, the positive real numbers $\totl{\w^{i}_{n}}$ converge--- note that $\totl{\w^{i}_{n}}$ is $0$ whenever $\w^{i}_{n}$ is not a finite real number.) If $\totl{\w^{i}_{n}}\rightarrow 0$, then a gluing parameter $\w_{n}\in\et 1{[0,\infty)}$  will be a function defined on our limiting family $\ex F$, as it is part of the choices involved in describing a limiting curve.
  
  \
  
  Define 
 \[\ex A^{\pm}_{n}:=\left\{\abs {\tilde z^{\pm}}\leq e^{-\frac{7R}{10}}\right\}\subset\et 11\] 
	 \[\ex A^{i}_{\frac {7R}{10},n}:=\{(\tilde z^{+},\tilde z^{-})\text{ so that }\tilde z^{+}\tilde z^{-}=\w^{i}_{n}\}\subset\ex A^+_{n}\times\ex A^{-}_{n}\]
 We shall need to keep track of what $\ex A_{n}$ is attached to, so use the notation $n^{\pm}$ to define $C^{i}_{n^{\pm}}$ as the component attached to the end of $\ex A_{n}^{i}$ with the boundary $\abs{\tilde z^{\pm}_{n}}=1$. (Assume, by passing to a subsequence, that $n^{\pm}$ is well defined independent of $i$.) As some annuli $\ex A^{i}_{k}$ have semi-infinite conformal modulus, we shall use the convention that $\ex A^{+}_{k}=\ex A^{i}_{\frac{7R}{10},k}$ on the rare occasions that we must mention these semi-infinite annuli.

Consider transition maps between $C^{i}_{n^{\pm}}$ and $\ex A^{i}_{\frac {7R}{10},n}$ to give transition maps between $C_{n^{\pm}}$ and $\ex A^{\pm}_{n}$. 
 The bound on the derivative of transition maps from item \ref{proposition derivative bound} of Proposition \ref{decomposition proposition}  on the region $\ex A^{i}_{\frac {6R}{10},n}-\ex A^{i}_{\frac {9R}{10},n}$ (and standard elliptic bootstrapping to get bounds on higher derivatives on the smaller region $\ex A^{i}_{\frac {7R}{10},n}-\ex A^{i}_{\frac {8R}{10},n}$) tells us that we can choose a subsequence so that the transition maps between $C^i_{n^{\pm}}$ and $\ex A^\pm_{n}$ converge to some smooth transition map between $C_{n^{\pm}}$ and $\ex A^{\pm}_{n}$. Similarly, by passing to a subsequence, we may require that the transition maps between any of our  annuli $\ex A^{+}_{k}$ with semi infinite conformal modulus and $C_{k^{+}}$ will converge.   We can modify our diffeomorphisms identifying $C_m^i$ with $C_m$ so that these transition maps all give exactly the same map. Denote the complex structure on $C_m$ induced by this identification by $j_i$. (We shall do something similar for the annuli $\ex A^{i}_{n}$ later on.) We can do this so that $j_i$ and the metrics given by this identification still converge in $C^\infty$ to those on $C_m$.
 
 \includegraphics{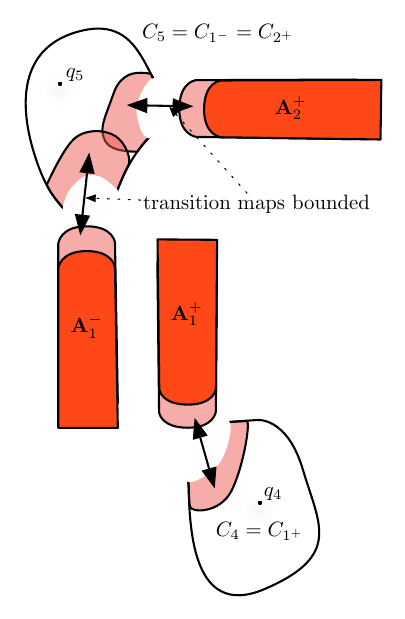}
 
 We shall now start to construct our family $\ex {\hat C}\longrightarrow \ex F$. The first step shall be to get some kind of convergence for each piece of our holomorphic curve.
 
  \

 For each $C_m$ choose some point $q_m\in  C_m$, and consider \[\Q^i_m:=f^i(q_m)\in \hat{\ex B}\] This sequence $\Q^{i}_{m}$ has a subsequence that converges in $\hat{\ex  B}$ to some point $\Q^{\infty}_m\in\hat{ \ex B}$.  Label the stratum of $\hat{\ex B}$ that contains $\Q^{\infty}_m$ by $\hat{\ex B}_m$, and consider the chart $\ex U_m$ from Lemma \ref{strata charts} containing $\hat{\ex B}_m$. As the derivative of $f^i$ restricted to $C_m$ is bounded, we can choose a subsequence so that $f^i(C_m)$ is contained well inside $\ex U_m$. (Note that given any smooth metric on $\hat {\ex B}$,  any point in the stratum $\hat{\ex B}_{m}$ is infinite metric distance to the boundary of $\ex U_{m}$. We can therefore choose a subsequence so that the image of $f^{i}(C_{m})$ is contained in the subset of $\ex U_{m}$ which consists of all points some arbitrary distance from the boundary of $\ex U_{m}$.)

 Recall that $C^{i}_{m}$ is a connected component of $\ex C^i-\bigcup \ex A^i_{\frac {8R} {10},k}$, and is contained some fixed distance in the interior of $\ex C^i-\bigcup \ex A^i_{\frac {9R} {10},k}$.
 As the derivative of $f^i$ is uniformly bounded on $\ex C^i-\bigcup \ex A^i_{\frac {9R} {10},k}$, we can use lemmas \ref{coordinate bound}, \ref{dbar of derivative}, and \ref{elliptic regularity} to get bounds on the higher derivatives of $f^i$ restricted to $C^i_m$.  Lemma \ref{bounded geometry} then tells us that if $f^i(q_m)$ converges  to $\Q^{\infty}_m\in \hat{\ex B}$, the geometry around $f^i(q_m)$ converges to that around $\Q^{\infty}_m$, so we can choose a subsequence so that $f^i$ restricted to $ C^i_m$ converges in some sense to a map  
 \[f_{m,\Q^{\infty}_m}: C_m\longrightarrow\hat{\ex B}\text{ so that }f_{m,\Q^{\infty}_m}(q_m)=\Q^{\infty}_m\subset\ex U_{m}\]  
 More specifically, remembering that  the image of $C^{i}_{m}$ is contained inside $\ex U_m$, we can use our $\ex T^{n}$ action on $\ex U_{m}$ to say this more precisely.  There exists some sequence $\tilde c_m^i\in\ex T^n$ so that
 \[\tilde c_m^i*f^i:C_m\longrightarrow \ex U_m\]
 converges in $C^\infty$ to $f_{m,\Q^{\infty}_m}:C_m\longrightarrow\ex U_m$. 
 
 Of course, there was a choice involved here. We could also have chosen a different  limit $\Q'_{m}$  of $f^{i}(q_{m})$, which would give a translate of $f_{m,\Q^{\infty}_{m}}$ by our $\ex T^{n}$ action. These choices will turn up as  parameters on our family. They will need to be `compatible' in the sense that these pieces will need to fit together.

Note that this ambiguity in the limit does not appear if we look only at the smooth part of these maps. Then $\totl {f^{i}}$ converges in $C^{\infty}$  to $\totl{f_{m,\Q^{\infty}_{m}}}$, which is the unique limit in $\totl{\hat{\ex B}}$. Another way to see the convergence above is that $\totl{Tf^{i}}:TC_{m}\longrightarrow \totl{T\hat{\ex B}}$ has a subsequence that converges in $C^{\infty}$  to a map to $\totl{T\hat{\ex B}}$. This specifies uniquely what the derivative of $f_{m,\Q^{\infty}_{m}}$ should be, but there is still a choice of which point $f_{m,\Q^{\infty}_{m}}(q_{m})$ should be  in $\hat{\ex B}$; this choice is the choice of the limit $y^{\infty}_{m}$ of $y^{i}_{m}$. Multiplying $f^{i}$ by $\tilde c_{m}^{i}$ allows us to translate $y^{i}_{m}$ so that the distance in any metric between $y^{i}_{m}$ and $y^{\infty}_{m}$ converges to zero. Then $f^{i}$ will converge in $C^{\infty}$ to $f_{m,\Q^{\infty}_{m}}$ because its derivative converges  and it converges at $q_{m}$.

 \
  
We now consider the analogous convergence on annular regions. Because the $\omega$-energy of $f^{i}$ restricted to our annular regions $\ex A^i_n$ is small, we can apply Lemma \ref{strong cylinder convergence} on page \pageref{strong cylinder convergence} to tell us that if the limit of the conformal modulus of $\ex A^i_n$ is not a finite real number, then $f^i$ restricted to the smooth parts of $\ex A^i_n$ converges   to some unique pair of holomorphic maps
 \[f^{\pm}_{n,\Q^{\infty}_{n^{\pm}}}:\{e^{-\frac {7R}{10}}>\abs {\tilde z_{n}^\pm}>0\}\longrightarrow\hat{\ex B}\] compatible with $f_{n^{\pm},\Q^{\infty}_{n^{\pm}}}:C_{n^{\pm}}\longrightarrow \hat{\ex B}$ and the transition maps. (Recall that $C_{n^{\pm}}$ is the component attached to $\ex A^{\pm}_{n}$.) Choose our subsequence so that $\totl{Tf^{i}}$ converges in $C^{\infty}$ on bounded subsets of $T\ex A^{i}_{n}$ close to the appropriate boundaries to $\totl{Tf^{\pm}_{n,\Q^{\infty}_{n^{\pm}}}}$.
  
Lemma \ref{strong cylinder convergence} gives that $\totl{f^{\pm}_{n,\Q^{\infty}_{n^{\pm}}}(\tilde z^{\pm})}$ converges as $\totl{\tilde z^{\pm}}\to 0$.    As we have made the assumption that $J$ is $\dbar\log$ compatible, we can use the  removable singularity theorem for finite energy holomorphic curves in $\totl{\hat{\ex B}}$ to see that these limit maps extend uniquely to  smooth maps on $\{e^{-\frac {7R}{10}}>\abs {\tilde z_{n}^\pm}>\e x\}\subset\ex A^{\pm}_{n}$ for some $x>0$. (Recall that the notation $\abs{ct^{a}}>t^{x}$ means that $a<x$ or $a=x$ and $\abs c>1$.)

  \includegraphics{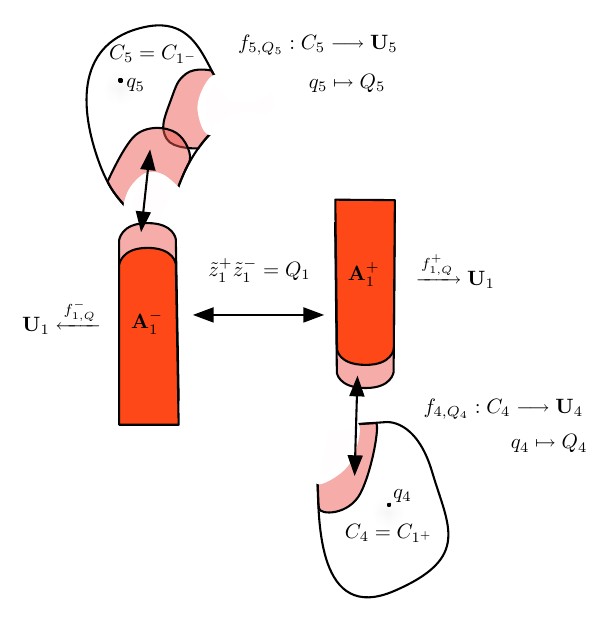}

   We may assume, after passing to a subsequence, that $f^i(\ex A^i_{\frac {6R}{10},n})$ are all contained inside a single $\ex U_n$. Recalling our convention that $C_{n^{\pm}}$ intersects $\ex A^i_n$, note that this means  that $\ex U_{n^{\pm}}$ intersects $\ex U_n$ (and $\totb{\ex U_{n^{\pm}}}\subset\totb{\ex U_n}$), so our sequence $\tilde c^i_{n^{\pm}}$ defined earlier also has an action on $\ex U_{n}$, and  
   \[\tilde c^i_{n^{\pm}}* f^i: \ex A_{\frac {7R}{10},n}^i\longrightarrow \ex U_n\] converges in $C^\infty$ on compact subsets adjacent to $C_{n^{\pm}}$ to  $f^\pm_{n,\Q^{\infty}_{n^{\pm}}}$. As we shall see in the proofs of Lemmas \ref{matching} and \ref{cinfty} below, Lemma \ref{strong cylinder convergence} allows us to strengthen the above  convergence.
   
   \
   
 We now have a type of convergence of the individual pieces we have cut our holomorphic curve into. The limit is unique in $\totl{\hat{\ex B}}$ and even $\totl{T\hat{\ex B}}$, but there is a choice of a translation in defining the limit of each piece in $\hat{\ex B}$. The translation chosen for all the pieces must be compatible or we will not be able to glue together the resulting limit pieces to get a curve.
 
  We shall now define a model for the exploded structure on our family that is too large, as it ignores the requirement that these pieces must fit together. Later, we shall cut down this model to reflect the requirement of the pieces gluing together.
  
\

 Define $\ex V_m \subset \ex U_m$ to be the exploded manifold consisting of all points $\Q_m\in\ex U_m$ so that
 
 \begin{enumerate}
 \item There exists some $\tilde  c$ so that $\Q_m=\tilde c*\Q^{\infty}_m$ 
 \item  Defining 
 \begin{equation}\label{fmydef} f_{m,\Q_m}:=\tilde c*f_{m,\Q^{\infty}_m}\ ,\end{equation}
 the image of $C_{m}$, $f_{m,\Q_{m}}(C_{m})$ is contained well inside $\ex U_{m}$ in the sense that the distance to the boundary of $\ex U_{m}$ is greater than $1$.
 \item If $C_{m}$ is attached to $\ex A^{\pm}_{n}$, (so $m=n^{\pm}$) then define 

 \begin{equation}\label{fpmdef}f^{\pm}_{n,\Q_{n^{\pm}}}:=\tilde c*f^\pm_{n,\Q^{\infty}_{n^{\pm}}}\end{equation}

The image of the smooth part of $\ex A^{\pm}_{n}$, $\{f^{\pm}_{n,\Q_{n^{\pm}}}(\tilde z),\ 0<\abs{\tilde z}<e^{-\frac {7R}{10}}\}$ is also contained well inside $\ex U_{n}$.
  \end{enumerate}

 Note that $\ex V_{m}$ is a smooth exploded manifold which includes all $\Q_m$ which are  limits for $\Q^i_m$. For such $\Q_m$, the map $f_{m,\Q_m}$ will be holomorphic, but for other $\Q_m$, this may not be the case. 

\

We shall consider the family $\ex F$ as a sub exploded manifold of  

\[\prod_m\ex V_m\times \prod_n\ex W_n\subset \prod_m\ex U_m\times \prod_n\ex W_n\]

The $\ex W_n$ above stands for the `gluing parameter' for identifying coordinates $\tilde z_{n}^\pm$ on $\ex A^\pm_n$ via 
\[\tilde z_{n}^+\tilde z_{n}^-=\w_n\in\ex W_n:=\{\abs{\tilde z}<e^{-2R}\}\subset\et 11\] 
   
   \
   
 We shall consider the following sequence of points \[\Q^i:=(f^i(q_1),f^i(q_2),\dotsc,\w^i_{n_1},\dotsc)\in\prod_m \ex U_m\times \prod_n\ex W_n\]
  where the conformal modulus of $\ex A^i_{n}$ is equal to $-\log \w^i_{n}$.
 
 \
 
 Note that there is a transitive action of $\ex T^k$ on $
 \prod \ex V_m\times \prod \ex W_n$ which is the action from Lemma \ref{strata charts} on $\ex V_m\subset \ex U_m$, and multiplication by some coordinate of $\ex T^k$ on $\ex W_n$. Our family will be given by a complete inclusion
 \[\ex F\longrightarrow \prod_m\ex V_m\times \prod_n\ex W_n\subset \prod_m\ex U_m\times \prod_n\ex W_n\]
  satisfying the following conditions:
  \begin{enumerate}
  \item \label{F condition 1}The image of $\ex F$ contains some point $\Q\in \prod_m\ex V_m\times \prod_n\ex W_n $ which is a limit of $\Q^i$, and the image of $\ex F$ is given by the orbit of $\Q$ under the action of some subgroup of $\ex T^k$.
  \item\label{F condition 2} The distance in any smooth metric from some subsequence of $\Q^i$ to the image of $\ex  F$ converges to $0$.
  \item There is no other inclusion satisfying the above conditions which has smaller dimension than $\ex F$.
  \end{enumerate}
  
  It is clear that such an $\ex F$ must exist. The entire space satisfies the first two properties, so we may simply choose an $\ex F$ satisfying the first two properties which has minimal dimension. It follows from the definition of an action of $\ex T^{k}$  (page \pageref{T action}) that $\ex F$ is smooth as it is the orbit of some $\ex T^{n}$ action. 
The minimality of the dimension of $\ex F$ is a trick which will ensure that our limiting maps  $f^{\pm}$ defined on $\ex A^{\pm}$ glue together well. (See Lemma \ref{matching} below.)

\

For each component $C_{m}$, there is a corresponding map 
\[y_{m}:\ex F\longrightarrow \ex U_{m} \subset \hat {\ex B}\]
which should be thought of as specifying where the point $q_{m}\in C_{m}$ is sent. With reference to equation (\ref{fmydef}), for any point $y\in \ex F$, define
\[ f_{m,y}:C_{m}\longrightarrow \hat{\ex B}\]
by
\begin{equation}\label{fmydef2}f_{m,y}:=f_{m,y_{m}}=\frac {y_{m}}{y^{\infty}_{m}}*f_{m,y^{\infty}_{m}}\end{equation}
Also, recalling the notational convention that  $m=n^{\pm}$ if $\ex A_{n}^{\pm}$ is attached to $C_{m}$, define
\[f^{\pm}_{n,y}:\ex A^{\pm}_{n}\longrightarrow \ex B\]
by
\begin{equation}\label{fpmdef2}f^{\pm}_{n,y}:=f^{\pm}_{n,y_{n^{\pm}}}=\frac{y_{n^{\pm}}}{y_{n^{\pm}}^{\infty}}*f_{n,y^{\infty}_{n^{\pm}}}\end{equation}
as in equation (\ref{fpmdef}).

\

Now let us construct our family of curves $\hat{\ex C}\longrightarrow\ex F$. We shall have charts on $\hat{\ex C}$ given by $\hat C_m:=C_m\times\ex F$, $\hat{\ex A}_{k}:=\ex A_{k}\times \ex F$ for annuli $\ex A_{k}$ with semi infinite conformal modulus,  and $\hat{\ex A}_n$ for other annuli $\ex A_{n}$. $\hat {\ex A}_n$ has coordinates $(\tilde z_{n}^+,\tilde z_{n}^-,\Q)$ where $\tilde z^{\pm}_{n}\in\ex A^{\pm}_{n}$, $\Q\in \ex F$ and $\w_n(\Q)=\tilde z_{n}^+\tilde z_{n}^-$.
\[\begin{array}{ccl}\hat{\ex A}_{k}&\longrightarrow &\ex F
\\ \downarrow && \downarrow \w_{n}
\\ \ex A_{n}^{+}\times \ex A_{n}^{-}& \xrightarrow{ \tilde z^{+}\tilde z^{-} }& \et 11\end{array}\]

 Transition maps between $\hat C_m$ and $\hat{\ex A}_n$ are simply given by the transition maps between $C_m$ and $\ex A^\pm_{n}$ times the identity on the $\ex F$ component. For example, if $\phi$ is a transition map between $C_{m}$ and $\ex A^{+}_{n}$, the corresponding transition map between $\hat C_{m}$ and $\hat{\ex A}_{n}$ is given by 
\[(z,\Q)\mapsto\left(\phi(z),\frac{\w_{n}(\Q)}{\phi(z)},\Q\right) \]
  This describes the exploded manifold $\hat{\ex C}$. The map down to $\ex F$ is simply given by the obvious projection to the last component of each of the above charts. This gives a smooth family of exploded curves.

\begin{lemma}\label{matching}
 If $\Q\in \ex F$, and ${\w^{i}_{n}}\to {\w_{n}(y)}$, then $f^{\pm}_{n,\Q}$ can be glued  by the identification $\tilde z_{n}^+\tilde z_{n}^-=\w_n(\Q)$. The condition ${\w^{i}_{n}}\to {\w_{n}(y)}$  is always satisfied for $y\in \ex F$ if the limit of the conformal modulus of $\ex A^{i}_{n}$ is finite.
\end{lemma}

\pf

   In the case that the limit of the conformal modulus of $\ex A^i_n$ is finite, $\lim \w^{i}_{n}=\w_{n}(\Q)$ for any $\Q\in \ex F$. If not, we could simply restrict $\ex F$ to the subset  $\ex F_{0}\subset \ex F$ consisting  of all points that satisfy this condition. This subset $\ex F_{0}$ would include any point in $\ex F$ that is the limit of $\Q^{i}$. $\ex F_{0}$ would also satisfy the other conditions required of $\ex F$, but have smaller dimension contradicting the minimality of $\ex F$.
    Therefore in the case that $\w_{n}(\Q)$ is finite, we can glue the limit, because we have $f^{i}$ converging on $C_{n^{\pm}}$ and $\ex A^{i}_{n}$, and transition maps between these also converging.
    
    Similarly, note that the limit of $\Q^{i}$ in $\totl{\ex F}$ is in a $0$-dimensional stratum of $\totl{\ex F}$, or else $\ex F$ will not have minimal dimension.

  In the case that the limit of the conformal modulus of $\ex A^i_n$ is not finite, we want to  glue $f^+_{n,\Q}(\tilde z_{n}^+)$ to $f^{-}_{n,\Q}(\tilde z_{n}^-)$ over the region where the smooth coordinates $ \totl{\tilde z^+_{n}}=\totl{\tilde z^-_{n}}=0$ via the identification $\tilde z_{n}^+\tilde z_{n}^-=\w_n(\Q)$. Define the following map
   
    \[\phi_n:\ex F\longrightarrow \ex T^l \text{ (the group acting on $\ex U_n$)}\] which detects the failure for $f^\pm_n$ to glue for any $\Q\in \ex F$: If $\totl{\w_n(\Q)}=0$,  and $\tilde z_{n}^+\tilde z_{n}^-=\w_n(\Q)$ so that $\totl{\tilde z^+}=\totl{\tilde z^-}=0$, let $\phi_{n}(\Q)$ be the element of $\ex T^l$   so that $\phi_n(\Q)*f_{n,\Q}^+(\tilde z_{n}^+)=f_{n,\Q}^-(\tilde z_{n}^-)$. This is defined for any $\Q$ which is the limit of $\Q^{i}$, because for such $\Q$,
\[ \totl{   f_{n,\Q}^+(\tilde z_{n}^+)}=\totl{f_{n,\Q}^-(\tilde z_{n}^-)}\]
    
    There exists some homomorphism from the subgroup of $\ex T^k$ acting transitively on $\ex F$ to $\ex T^l$ so that  $\phi_{n}$ is equivariant. We can therefore uniquely extend $\phi_n$ to the rest of $\ex F$ in an equivariant way. 
    
    We want to show that $\phi_n$ is identically $1$ on $\ex F$. By the construction of $\ex F$, for each $\Q^i$, there exists a close point $\check{\Q}^i$ inside $\ex F$ so that the distance in any smooth metric between $\check \Q^i$ and $\Q^i$ converges to $0$ as $i\rightarrow\infty$. We shall show below that $\phi_n(\check{\Q}^i)$ converges to $1$:
    
    $\totl{Tf^{i}}$ converges to $\totl{Tf_{m,y_{m}^{\infty}}}$ on $C_{m}$, and converges to $\totl{Tf^{\pm}_{n,y^{\infty}_{n^{\pm}}}}$ on compact subsets of $\ex A^{\pm}_{n}$ adjacent to $C_{ n^{\pm}}$. It follows that given any fixed choice of point $x$ in $\ex A^{\pm}_{n}$, the metric distance between $f^{\pm}_{n,\check\Q^{i}}(x)$ and $f^{i}(x)$ converges to $0$ as $i\to\infty$. On the other hand, in the interior of $\ex A^{i}_{n}$, both $f^{\pm}_{n,\check \Q^{i}}$ and  $f^{i}$ are approximated by a cylinder, as specified in Lemma \ref{strong cylinder convergence}. Therefore, for each $i$, we may choose an $x_{i}$ close to the center of $\ex A^{i}$, identify this with an $x_{i^{\pm}}$ in $\ex A^{\pm}_{n}$, and then the metric distance between $f^{\pm}_{\check\Q^{i},n}(x_{i}^{\pm})$ and $f^{i}(x_{i})$ will converge to $0$.  As the metric distance between $
    \w_{n}(\check\Q^{i})$ and $\w_{n}(\Q^{i})$ converges to $0$, the metric distance between $x_{i}^{+}x_{i}^{-}$ and $\w_{n}(\check \Q^{i})$ converges to $0$. It follows that the distance between $f^{+}_{\check\Q^{i},n}(x_{i}^{+})$ and $f^{-}_{\check\Q^{i},n}(\w_{n}(\check \Q^{i})/x_{i}^{+})$ converges to $0$, so   $\phi_{n}(\check\Q^{i})$ converges to $1$.

    As the limit of $\Q^{i}$ is in a zero dimensional stratum of $\totl{\ex F}$, either $\phi_{n}$ is constant and equal to $1$ on this stratum, or $\phi_{n}$ is non constant on this stratum, and must take all values in $\mathbb C^{*}$. In either case, $\phi^{-1}(1)\subset \ex F$ is nonempty and satisfies the conditions \ref{F condition 1} and \ref{F condition 2} above, so by the minimality of the dimension of $\ex F$, it must be all of $\ex F$.  This tells us that if $\totl{\w_{n}(\Q)}=0$, then we can glue together $f^{\pm}_n$ without any modifications. 
 
 \stop
 
 \

 Note that choosing any point $\Q\in \ex F$ so that $\Q$ is a limit of $\Q^i$ (such a point must exist by the definition of $\ex F$), the above lemma allows us to glue together $f_{m,\Q}$ and $f^{\pm}_{n,\Q}$ to obtain our limiting holomorphic curve $f$ for Theorem \ref{completeness theorem}. Note that condition \ref{decomposition stability} from Proposition \ref{decomposition proposition} ensures that $f$ is a stable holomorphic curve. 
 
 \
 
 It remains to prove that $f^{i}$ converges to $f$ in the sense of Theorem \ref{completeness theorem}. This involves constructing a family of curves containing $f$ large enough to capture the behavior of $f^{i}$ for $i$ large enough--- for example, if the smooth part of the domain of $f$ has more nodes than the smooth part of the domains of $f^{i}$, we need to construct a family of curves containing $f$ which smooths these extra nodes. Unfortunately, whenever nodes must be smoothed, we need to make gluing choices, however as our smoothed curves need not be holomorphic themselves, we are not required to go through the full details of a traditional holomorphic curve `gluing theorem'.

 To understand the gluing procedure that will follow, roughly what we shall need is a `gluing' map which takes a matching pair of maps from $\ex A^{\pm}_{n}$, and produces a family of maps with domain equal to  the family of annuli obtained by gluing $\ex A^{\pm}_{n}$. We shall also need a cutting map, which produces a matching pair of maps from $\ex A^{\pm}_{n}$ from a given map from $\ex A^{i}_{n}$. These gluing and cutting maps must be compatible in the sense that if we cut a given map, then apply the gluing map, then restrict the domain to $\ex A^{i}_{n}$, we get back the original map. Roughly speaking, we shall use the gluing map to glue together the maps $f^{\pm}_{n}$ to create a family of curves $\hat f$ containing our limit of $f^{i}$.  We  shall also use the cutting map to cut apart each $f^{i}$, then glue it back together to obtain a family of curves $\hat f^{i}$ containing $f^{i}$ so that $\hat f^{i}$ converges to $\hat f$.

 \
 
 Make the following gluing choices.

\begin{enumerate}

\item 
Choose some sequence $\check \Q^i$ of points in  $\ex F$ so that the distance between $\check \Q^i$ and $\Q^i$ converges to $0$. 

\

\item\label{cc}
We need to take care of the different complex structures obtained by gluing $\ex A^{\pm}_{n}$ by $\check \w^{i}_{n}:=\w_{n}(\check \Q^{i})$ and $\w^{i}_{n}$. For each $i$, choose an almost complex structure $j_{i}$ on $\ex A^{+}_{n}$ as follows: \\ Choose smooth isomorphisms
\[  \Phi^{i}_{n}:\ex A^{+}_{n}\longrightarrow \ex A^{+}_{n}\]
so that 
\begin{enumerate}
\item \[\Phi^{i}_{n}\lrb{\tilde z_{n}^{+}}=\tilde z_{n}^{+}\text{ for }\abs {\tilde z_{n}^{+}}\geq e^{-\frac {8R}{10}}\]
\item \[\Phi^{i}_{n}\lrb{\tilde z_{n}^{+}}=\frac{\check \w^{i}_{n}}{\w^{i}_{n}}\tilde z_{n}^{+}\text{ for }\abs{\tilde z_{n}^{+}}\leq e^{-R}\]
\item On the region $ \totl{\tilde z_{n}^{+}}\neq 0$, the sequence of maps $\{\Phi^{i}_{n}\}$ converges in $C^{\infty}$ to the identity.
\end{enumerate}

Now define $j_{i}$ on $\ex A^{+}_{n}$ to be the pullback under $\Phi^{i}_{n}$ of the standard complex structure. Using the standard complex structure on $\ex A^{-}_{n}$, we then get our $j_{i}$ defined on $\hat{\ex A}_{n}$. This $j_{i}$ is compatible with  the  $j_{i}$ already defined on  $\hat C_{m}$, so we get $j_{i}$ defined on $\hat{\ex C}$. Note that $j_{i}$ restricted to the fiber over $\check \Q^{i}$ is the complex structure on $\ex C^{i}$.  

\

From now on, we shall use these new coordinates on $\ex A^{i}_{n}$, so $\ex A^{i}_{n}$ is obtained by gluing $\ex A^{\pm}_{n}$ using the gluing parameter $\check \w_{n}^{i}:=\w_{n}(\check \Q^{i})$ instead of $\w^{i}_{n}$.

\

\item Define a linear gluing map as follows:
\begin{enumerate}
\item Chose some smooth cutoff function 
\[\rho:\mathbb R^{*}\longrightarrow [0,1]\]
so that \[\rho(x)=0\text{ for all }x\geq e^{-\frac {8R}{10} }\]
\[\rho(x)=1\text{ for all }x\leq e^{-\frac{9R}{10}}\] 
Extend this to a map
\[\rho:\mathbb R^{*}\e{\mathbb R}\longrightarrow [0,1]\]
satisfying both the above conditions--- so $\rho(x\e 0)$ is equal to the old $\rho(x)$,  $\rho(x\e a)=0$ for all $a<0$, and   $\rho(x\e a)=1$ for all $a>0$. ($\rho$ was defined first on $\mathbb R^{*}$ so that it was clear what `smooth' meant.)

\item\label{cutnglue}

Given  maps $\phi^+,\phi^-:\et 11\longrightarrow \mathbb C^k$ which vanish at $z=0$ define the gluing map 

\[G_{(\phi^+,\phi^-)}(\tilde z^+,\tilde z^-):=\rho \left( \abs {\tilde z^-}\right)\phi^+(\tilde z^+)+\rho \left( \abs {\tilde z^+}\right)\phi^-(\tilde z^-)\]

\end{enumerate}

Note that if $\phi^+$ and $\phi^-$ are smooth, $G_{(\phi^+,\phi^-)}:(\et 11)^{2}\longrightarrow \mathbb C^k$ is smooth, so $G$ can be regarded as a smooth gluing map.

\begin{claim}\label{cG}If $\phi^{+}$ and $\phi^{-}$ are small in $C^{\infty,\delta}$, then $G_{(\phi^{+},\phi^{-})}$ is small too. \end{claim}

The definition of $C^{\infty,\delta}$ is found in section $7$ of  \cite{iec}. The cutoff function $\rho(\abs{ \tilde z^{\pm}})$ is smooth, and therefore $C^{\infty,\delta}$ for all $\delta<1$. Multiplication and addition are continuous operations on $C^{\infty,\delta}$, so \[G_{(\cdot,\cdot)}:(C^{\infty,\delta}(\et 11))^{2}\longrightarrow C^{\infty,\delta}\lrb{(\et 11)^{2}}\]  is continuous, and Claim \ref{cG} holds.

\

An important aspect of our gluing construction shall be that it is localized within the region of $\hat {\ex A}_{n}$ which does not intersect $\hat {\ex C}_{n^{\pm}}$. The particular property we shall need is as follows:

\begin{claim}\label{Gloc}
\[G_{\phi^{+},\phi^{-}}(\tilde z^{+},\tilde z^{-})=\phi^{\pm}(\tilde z^{\pm})\ \text{ when } \abs{\tilde z^{\pm}}\geq e^{-\frac {8R}{10}}\ \text{ and }\abs{\tilde z^{+}\tilde z^{-}}<e^{-2R}\]
\end{claim}

 Indeed, when $\abs {\tilde z^{\pm}}\geq e^{\frac {-8R}{10}}$, the cut off function $\rho(\abs{\tilde z^{\pm}})=0$. When $\abs{\tilde z^{+}\tilde z^{-}}<e^{-2R}$ too, (as is the case in $\hat{\ex A}_{n}$), $\abs{\tilde z^{\mp}}<e^{-R}$, so $\rho(\abs{\tilde z^{\mp}})=1$.

\

\item Define a linear `cutting' map: as follows:
\begin{enumerate}
\item Choose a smooth cutoff function $\beta:\mathbb R^{*}\longrightarrow[0,1]$ satisfying the following:
\[\beta(x)+\beta\lrb{x^{-1}}=1\]
\[\beta(x)=1\text{ for all } x>e^{\frac R{10}}\]

\item \label{cut}
Given a map $\phi:\ex A^{i}_{n}\longrightarrow \mathbb C^k$, and  $\check \w_{n}^{i}:=\w_{n}(\check \Q^{i}) \in\mathbb C^*$, where 
\[\ex A^{i}_{n}:=\{\tilde z^{+}\tilde z^{-}=\check \w^{i}_{n}\}\subset(\et 11)^{2}\]
we can define the cutting of $\phi$ as follows:
\[\phi^+(\tilde z^+):=\beta \left(\frac{\abs{\tilde z^+}^2}{\abs {\check \w^{i}_{n}}}\right)\phi(\tilde z^+)\]
\[\phi^-(\tilde z^-):=\beta \left(\frac{\abs{\tilde z^-}^2}{\abs {\check \w^{i}_{n}}}\right)\phi(\tilde z^-):=
\beta \left(\frac{\abs{\tilde z^-}^2}{\abs {\check \w^{i}_{n}}}\right)\phi\left(\tilde z^{+}=\frac {\check \w^{i}_{n}}{\tilde z^-}\right)\]
As $\phi^{\pm}(\tilde z^{\pm})$ vanishes for $\totl{\tilde z^{\pm}}$ small, we may extend $\phi^{\pm}$ to be zero for $\tilde z^{\pm}$ small.

\end{enumerate}

\

\begin{claim}\label{c6.4} $G_{(\phi^+,\phi^-)}$ restricted to $\tilde z^+\tilde z^{-}=\check \w^{i}_{n}$ is equal to $\phi$. \end{claim} To prove Claim \ref{c6.4}, note that as required in the definition of $\ex F$,  $\abs{\check \w^{i}_{n}}< e^{-2R}$. The functions $\phi^{\pm}(\tilde z^{\pm})$ are supported in the set where $\abs{\tilde z^{\pm}/\tilde z^{\mp}}>e^{-\frac{R}{10}}$, so $\abs{\tilde z^{\mp}}<e^{-\frac{19R}{20}}$, but $\rho(\abs {\tilde z^{\mp}})$ is equal to $1$ where $\abs{\tilde z^{\mp}}\leq e^{-\frac{9R}{10}}$. Therefore, $\rho$ may be dropped from our expression for $G_{\phi^{+},\phi^{-}}$ and 
\[G_{\phi^{+},\phi^{-}}(\tilde z^{+},\tilde z^{-})=\beta \left(\frac{\abs{\tilde z^+}^2}{\abs {\check \w^{i}_{n}}}\right)\phi(\tilde z^{+})+\beta \left(\frac{\abs{\tilde z^-}^2}{\abs {\check \w^{i}_{n}}}\right)\phi(\check \w^{i}_{n}/\tilde z^{-})=\phi(\tilde z^{+})\]
Therefore Claim \ref{c6.4} holds.

\

We chose our cutoff function $\beta$ so that our cutting construction takes place within the region of $\ex A^{i}_{n}$ which does not intersect $C_{n^{\pm}}$:

\begin{claim}\label{cloc} \[\phi^{\pm}_{n}(\tilde z^{\pm})=\phi(\tilde z^{\pm}) \ \text{ where }\tilde z^{\pm}>e^{-\frac {8R}{10}}\]
\end{claim} 

Indeed, as $\abs{\check \w^{i}_{n}}<e^{-2R}$, and $\abs{\tilde z^{\pm}}>e^{-\frac {8R}{10}}$, 
\[\beta\left(\frac{\abs{\tilde z^\pm}^2}{\abs {\check \w^{i}_{n}}}\right)= 1 \ \ \ \text{ and }\beta\left(\frac{\abs{\tilde z^\mp}^2}{\abs {\check \w^{i}_{n}}}\right)=0\]

\item In this step, we use the above cutting map to construct $\phi^{i\pm}_{n}$ from $f^{i}_{n}$.

\

Use Lemma \ref{cylinder bound} to ensure that the image of each of our annuli of finite conformal modulus, $f^i(\ex A^i_{\frac {7R}{10},n})$ is contained in some coordinate chart appropriate for Lemma \ref{strong cylinder convergence}. Choose this coordinate chart to be contained inside $\ex U_n$ so that the $\ex T^k$ action on $\ex U_n$ just consists of multiplying coordinate functions by a constant. 

In the case that $\ex A^{i}_{n}$ has finite conformal modulus,  Lemma \ref{strong cylinder convergence} together with the conditions on our coordinate change for $\ex A^{i}_{n}$ from item \ref{cc} from page \pageref{cc} tells us that 
\begin{equation}\label{phi1}\begin{split} f^i_n(z)&=(e^{\phi^i_{n,1}(z)}c_{n,1}^i z^{\alpha_1},\dotsc,e^{\phi^i_{n,k}(z)}c_{n,k}^i z^{\alpha_k},c_{n,k+1}^{i}+\phi^{i}_{n,k+1}, \dotsc)
\\&:=e^{\phi^i_n} c^i_nz^\alpha\end{split}\end{equation}

where $\phi^i_n$ is exponentially small on the interior of $\ex A^i_{R,n}$ as required by Lemma \ref{strong cylinder convergence}. In the case of a family of targets $\hat{\ex B}$, we may ensure that $\phi^{i}_{n}$ is a vertical vector field. Use the linear cutting map from item \ref{cut} above on $\phi^{i}_{n}$ to produce $\phi^{i\pm}_{n}$. In the case of a family of targets,  $\phi^{i\pm}_{n}$ is also vertical. 
\[\phi^{i}_{n}\xrightarrow{\ \ \ \text{cut}\ \ \ } \phi^{i+}_{n},\phi^{i-}_{n}\]
 We shall use these $\phi^{i\pm}_{n}$ in our gluing map to produce a family of maps on different domains.
\[\phi^{i}_{n}\xrightarrow{\ \ \ \text{cut}\ \ \ } \phi^{i+}_{n},\phi^{i-}_{n}\xrightarrow{\ \ \ \text{glue}\ \ \ } G_{(\phi_{n}^{i+},\phi^{i-}_{n})}\]

As we didn't define our cutting map in the case when the conformal modulus of $\ex A^i_n$ is not a finite number,  we may define $\phi^{i\pm}_{n}$ separately in that case so that

\begin{equation}\label{phi2}f^i_n(\tilde z^+):=e^{\phi^{i+}_n(\tilde z^{+})} c^i_n(\tilde z^+)^\alpha\ \ \  \text {when } \totl{\tilde z^{+}}\neq 0\end{equation}
\begin{equation}\label{phi3}f^i_{n}(\tilde z^-):=e^{\phi^{i-}_n(\tilde z^-)}c^i_n\left(\frac {\check \w^i_n}{\tilde z^-}\right)^\alpha\ \ \  \text {when } \totl{\tilde z^{-}}\neq 0\end{equation}
 Lemma \ref{strong cylinder convergence}, Lemma \ref{J embedding}, and the removable singularity theorem for holomorphic curves gives that $\phi^{i\pm}_{n}$ is smooth and vanishes when $\totl{ \tilde z^{\pm}}=0$.

\

If $\Q^\infty\in\ex F$ is some  limit of $\Q^i$, then $f^{\pm}_{n,\Q^\infty}$ as defined in equation (\ref{fpmdef2}) on page \pageref{fpmdef2} is given in these coordinates  by 
\begin{equation}\label{phi4}f^+_{n,\Q^\infty}(\tilde z^+):=e^{\phi^{+}_n(\tilde z^+)} c_n(\tilde z^+)^\alpha\end{equation}
\[:=(e^{\phi^{+}_{n,1}(\tilde z^+)}c_{n,1} (\tilde z^+)^{\alpha_1},\dotsc,e^{\phi^+_{n,k}(\tilde z^+)}c_{n,k} (\tilde z^+)^{\alpha_k},c_{n,k+1}+\phi^{+}_{n,k+1}(\tilde z^{+}),\dotsc )
\]
\begin{equation}\label{phi5}f^-_{n,\Q^\infty}(\tilde z^-):=e^{\phi^-_n(\tilde z^-)}c_n\left(\frac {\Q^\infty_n}{\tilde z^-}\right)^\alpha\end{equation}
where $\phi^\pm$ are smooth and vanish when $\totl{z^\pm}=0$.

\

\item In this final step of our gluing construction, we construct our family of maps $\hat f^{i}$  on $\hat{\ex A}_{n}$.

\

 If $p_{1}=c*p_{2}$ using one of our $\ex T^{n}$ actions, use the notation $c=\frac {p_1}{p_2}$, so $p_1=\frac {p_1}{p_2}*p_2$. Recall that $\Q_{n^+}$ records the position of $C_{n^{+}}$, which is attached to  $\ex A^{+}_n$, so $\frac{\Q_{n^{+}}}{\check \Q^i_{n^{+}}}$ makes sense.
 
Define $F^i_n:\hat{\ex A}_n\longrightarrow \hat{\ex B}$
by
\[F^i_n(\tilde z^+,\tilde z^-,\Q):=\frac{\Q_{n^{+}}}{\check \Q^i_{n^{+}}}*c^i_n(\tilde z^+)^\alpha\]
 where $c^{i}_{n}(\tilde z^{+})^{\alpha}$ has the same meaning as in equation (\ref{phi1}). Similarly, define
\[F_n(\tilde z^+,\tilde z^-,\Q):=\frac{\Q_{n^{+}}}{ \Q^\infty_{n^{+}}}*c_n(\tilde z^+)^\alpha\]
where $\Q^{\infty}$, $c_{n}$ and $\alpha$ are as in equation (\ref{phi4}).

\

Define $\hat f^i:\hat{\ex A}_n\longrightarrow \hat{\ex B}$ using the cutting and gluing maps above by

\[\hat f^i(\tilde z^+,\tilde z^-,\Q):=e^{G_{(\phi^{i+}_n,\phi^{i-}_n)}(\tilde z^+,\tilde z^-)}F^i_n(\tilde z^+,\tilde z^-,\Q)\]
Similarly, define $\hat f:\hat{\ex A}_n\longrightarrow \hat{\ex B}$ by
\[\hat f(\tilde z^+,\tilde z^-,\Q):=e^{G_{(\phi^{+}_n,\phi^{-}_n)}(\tilde z^+,\tilde z^-)}F_n(\tilde z^+,\tilde z^-,\Q)\]

\end{enumerate}

 \begin{remark}\label{r1}Claim \ref{c6.4} implies that $\hat f^{i}$ restricted to $\Q=\check \Q^{i}$, is equal to $f^{i}$ restricted to $\ex A^{i}_{n}$.\end{remark}
 \begin{remark}\label{r2} Claims \ref{Gloc} and \ref{cloc} imply that $\hat f^{i}$ restricted to the region which intersects $\hat C_{n^{+}}$ is just equal to $ f^{i}$ translated by $\frac{\Q_{n^{+}}}{\check \Q^i_{n^{+}}}$. Lemma \ref{matching} implies that $\hat f^{i}$ restricted to the region which intersects $\hat C_{n^{-}}$ is just equal to $f^{i}$ translated by $\frac{\Q_{n^{-}}}{\check \Q^i_{n^{-}}}$.\end{remark}
 
 \begin{remark}\label{r3} Claim \ref{Gloc} and Lemma \ref{matching} imply that $\hat f$ restricted the regions which intersect $\hat C_{n^{\pm}}$ is equal to $ f^{\pm}_{n,\Q}(\tilde z^{\pm})$.\end{remark}

\

To define $\hat f^{i}:\hat{\ex C}\longrightarrow \hat{\ex B}$, we still need to define $\hat f^{i}$ on our other coordinate charts $\hat C_{m}$ and $\hat {\ex A}_{k}$ for semi-infinite annuli $\ex A_{k}$.

\

Define $\hat f^i$ on $\hat C_m$  by translating $f^i$ appropriately depending on the difference between $\check \Q^i_{m}$ and $\Q_{m}$, so
\[\hat f^{i}(z,\Q)=\frac {\Q_{m}}{\check \Q^{i}_{m}}*f^{i}(z)\]
Similarly, if $\ex A_{k}$ is an annulus with semi-infinite conformal modulus attached to $\hat C_{m}$, define $\hat f^{i}$ on $\hat{\ex A}_{k}$ using the same formula above. Remark \ref{r2} implies that on the overlap of any two coordinate charts, $\hat f^{i}$ is just the appropriate translate of $f^{i}$. Therefore these maps $\hat f^{i}$ piece together to give a well defined smooth map
\[\hat f^i:\hat {\ex C}\longrightarrow \hat{\ex B}\]
Remark \ref{r1} implies that   $f^i$ is given by the restriction of $\hat f^i$ to the fiber of $\hat{\ex C}\longrightarrow \ex F$ over $\check \Q^i\in\ex F$.

\

Define
\[\hat f:\hat{ C}_m\longrightarrow \hat{\ex B}\]
using the definition of $f_{m,\Q}(z)$ found in equation (\ref{fmydef2}) on page \pageref{fmydef2}
\[\hat f(z,\Q):=f_{m,\Q}(z)=\frac{y_{m}}{y_{m}^{\infty}}*f_{m,y^{\infty}_{m}}(z)\]
and similarly use equation (\ref{fpmdef2}) to define $\hat f$ on cylinders $\hat{\ex A}_{k}$ with infinite conformal modulus attached to $C_{m}$ by translating appropriately:
\[\hat f(z,\Q)=f^{+}_{k,\Q}(z)=\frac{y_{k^{+}}}{y_{k^{+}}^{\infty}}*f^{+}_{k,\Q_{k^{+}}^{\infty}}(z)\]
As $ f^{\pm}_{n,\Q}$ coincides with $ f_{n^{\pm},\Q}$ on their common domain of definition, Remark \ref{r3} implies that this defines a smooth map
\[\hat f:\hat{\ex C}\longrightarrow \hat{\ex B}\]
Note that $\hat f$ restricted to the fiber over  $\Q^{\infty}$ is a holomorphic curve. 

\

 As the vector fields $\phi^{\pm}$ are vertical with respect to $\hat{\ex B}\longrightarrow \ex G $, and translation in our coordinates preserves fibers of $\hat{\ex B}\longrightarrow \ex G$, our maps $\hat f^{i}$ and $\hat f$ are fiber preserving maps. We have now constructed our required families of curves.

  \[\begin{array}{ccccccc} (\ex {\hat C},j_i) & \xrightarrow{\hat f^i} &\hat{\ex B} &\ \ \ \ \ \ \ \ \ &
  (\ex {\hat C},j) & \xrightarrow{\hat f} &\hat{\ex B}
  \\ \downarrow & &\downarrow& &\downarrow & &\downarrow
  \\ \ex F& \longrightarrow & \ex G & & \ex F &\longrightarrow &\ex G\end{array}\]

Note that as required in Theorem \ref{completeness theorem}, $\hat f^{i}$ restricted to the curve over $\check y^{i}$ is equal to $f^{i}$. It remains to prove that $\hat f^{i}$ converges in $\C\infty1$ to $\hat f$.

\begin{lemma}\label{cinfty}
 \[\hat f^i:(\hat{\ex C},j_i)\longrightarrow \hat{\ex B}\] converges in $\C\infty1$ to 
 \[\hat f:(\hat{\ex C},j)\longrightarrow \hat{\ex B}\] 
\end{lemma}

\pf

We work in coordinates. First note that restricted to $\hat C_{m}$,  $\hat f^i$ converges in $C^\infty$. This also holds for $\hat {\ex A}_n$ if the size of $\ex A^i_n$ stays bounded. Note also that $j_i$ converges in $C^\infty$ to $j$.  

It remains only to show that for any $\delta<1$, $\hat f^i$ converges in $C^{\infty,\delta}$ to $\hat f$ on $\hat{\ex A}_n$. We shall deal with the interesting case that $\ex A^{i}_{n}$ is not semi-infinite first. Recall that $\hat {\ex A}_n$ has coordinates $(\tilde z^+,\tilde z^-,\Q)$ where $\Q\in \ex F$ and $\w_n(\Q)=\tilde z^+\tilde z^-$. Note that the maps $F^i_n:\hat{\ex A}_n\longrightarrow\hat{\ex B}$ defined above converge in $C^\infty$ to $F_n:\hat{\ex A}_n\longrightarrow\hat{\ex B}$. We have that

\[\hat f^i(\tilde z^+,\tilde z^-,\Q):=e^{G_{(\phi^{i+}_n,\phi^{i-}_n)}(\tilde z^+,\tilde z^-)}F^i_n(\tilde z^+,\tilde z^-,\Q)\]
 Similarly, we write
\[\hat f(\tilde z^+,\tilde z^-,\Q):= e^{G_{(\phi^+_n,\phi^-_n)}(\tilde z^+,\tilde z^-)}F_n(\tilde z^+,\tilde z^-,\Q)\]

We have already  chosen a subsequence so that  $\phi^{i\pm}_n( z^\pm)$ converges in $C^\infty$ on compact subsets of $\{0<\abs{z^{\pm}}\leq e^{-\frac {7R}{10}}\} $ to $\phi^\pm_n(z^\pm)$. Lemma \ref{strong cylinder convergence} and our cutting construction above tells us that for any $\delta<1$, there exists some constant $c$ independent of $i$  so that 
\[\abs{\phi_n^{i\pm}( z)}\leq c\abs { z}^\delta\]  
 Standard elliptic bootstrapping implies that the same inequality holds for  $\phi^\pm_n$ and all derivatives (with a different constant $c$).
 The above inequality and $C^{0}$ convergence we already have implies that for all $\delta'<\delta$, \[\abs{\phi_{n}^{i\pm}(z)-\phi_{n}^{\pm}(z)}\abs {z}^{\delta'}\] converges to $0$ as $i\rightarrow 0$. The same inequality for derivatives and the  $C^{\infty}$ convergence on complex subsets implies that the same expression for derivatives converges to $0$ as $i\rightarrow 0$. As this convergence holds for any $\delta'<1$, it is what is needed to show that $\phi^{i\pm}_n$ converges to $\phi^\pm_n$ in $\C\infty1$.

 Pre-composition with smooth maps preserves  $\C\infty1$ convergence. Therefore $\phi^{i\pm}_{n}(\tilde z^{\pm})$ converges in $\C\infty 1$ to $\phi^{\pm}_{n}(\tilde z^{\pm})$ on $\hat {\ex A}_{n}$.
 Then $G_{(\phi^{i+}_n,\phi^{i-}_n)}$ is obtained by adding these two functions multiplied by fixed smooth function cutoff functions. Therefore, $G_{(\phi^{i+}_n,\phi^{i-}_n)}$ converges in $\C\infty1$ to $G_{(\phi^+_n,\phi^-_n)}$ on $\hat{\ex A}_{n}$.  This implies that $\hat f^i$ converges to $\hat f$ in $\C\infty1$ as required. 

In the remaining case that $\ex A_{n}$ is semi-infinite, we may analogously define $F^{i}_{n}$, $F_{n}$, $\phi^{i+}_{n}$, and $\phi^{+}_{n}$ so that now $\hat f^{i}=e^{\phi^{i+}_{n}}F^{i}_{n}$. Again, $F^{i}_{n}$ converges to $F_{n}$ in $C^{\infty}$, and $\phi^{i+}_{n}$ converges to $\phi^{+}_{n}$ in $\C\infty1$, so $\hat f^{i}$ converges to $\hat f$ in $\C\infty1$.

 \stop

\

This completes  the proof of Theorem \ref{completeness theorem}.

\appendix\section{}
 \label{equivariant section}
In this appendix, we prove that any exploded manifold has coordinate charts with transition functions in a special form. 

Suppose that $\mathbb R^{n}\times \et mP$ and $\mathbb R^{n'}\times \et {m'}Q$ are two coordinate charts on $\ex B$ which intersect in a connected open set $U$ so  that $P$ is a face of the polytope $Q$. By making a monomial change of coordinates, we may also assume that the coordinates on $P$ are the first $m$ coordinates on $Q$, and that the remaining $m'-m$ coordinates on $Q$ vanish on $P$. Then the transition map between $\mathbb R^{n}\times \et mP$ and $\mathbb R^{n'}\times \et {m'}Q$ must be in the form
\[\phi(\tilde x,\tilde z)=(f,g_{1}\tilde z_{1},\dotsc,g_{m}\tilde z_{m}, g_{m+1},\dotsc,g_{m'})\]
where 
\[\mathbb R^{n}\times \totl{\et mP}\supset \totl U\xrightarrow{\ \ \ f\ \ \ }\mathbb R^{n'}\]
\[\mathbb R^{n}\times\totl{\et mP}\supset \totl U\xrightarrow{\ \ \ g_{i}\ \ \ } \mathbb C^{*}\]
are smooth functions. A special case is when $f$ and $g_{i}$ depend only on $\mathbb R^{n}$ and not on $\totl{\et mP}$. We shall prove that any exploded manifold $\ex  B$ can be covered  with coordinate charts modeled on open subsets of $\mathbb R^{n}\times \et mP$ with transition maps independent of $\totl{\et mP}$ in this sense. 

\begin{defna1} 
Define a $(\mathbb C^{*})^{m}$ action on $\mathbb R^{n}\times\et mP$ by
\[(c_{1},\dotsc,c_{m})*(x,\tilde z)=(x,c_{1}\tilde z_{1},\dotsc,c_{m}\tilde z_{m})\]
Similarly, if $(\tilde c_{1},\dotsc,\tilde c_{m})\in \ex T^{m}$, 
\[\text{let }(\tilde c_{1},\dotsc,\tilde c_{m})*(x,\tilde z)=(x,\tilde c_{1}\tilde z_{1},\dotsc,\tilde c_{m}\tilde z_{m})\text{  where defined. }\]

Call a map $f$ from a subset of $\mathbb R^{n}\times \et mP$ to $\mathbb R^{n'}\times \et {m'}Q$ equivariant if there exists a map \[\alpha:\ex T^{m}\longrightarrow\ex T^{m'}\] 
so that  
\[f(c*(x,\tilde z))=\alpha(c)*f(x,\tilde z)\]
 whenever the left hand side is defined.
\end{defna1}

Note that every map from $\mathbb R^{n}\times \et mP$ is equivariant when restricted to the interior of $P$.

\begin{defna2}\label{T action} 
A free action of $\ex T^{k}$ on an exploded manifold $\ex B$ is a map from a subset $\Dom (*)$ of $\ex T^{k}\times \ex B$ to $\ex B$
\[\begin{array}{ccc} \Dom(*)&\longrightarrow &\ex B
\\ (\tilde z,p)&\mapsto& \tilde z*p\end{array}\]
so that
\begin{enumerate}
\item There exists a covering of $\ex B$ by coordinates modeled on open subsets of $\mathbb R^{n}\times\et mP$ so that the $\ex T^{k}$ action is given by multiplication on the last $k$ coordinates of $\et mP$, in the sense that
\begin{enumerate}
\item if $p$ and $\tilde z*p$ are in the same coordinate chart, then $\tilde z*p$ is determined by multiplying the last $k$ coordinates of $p$ by $\tilde z$. 
\item If $q$ is a point in the same coordinate chart as $p$, and $q$ is obtained from $p$ by multiplying the last $k$ coordinates by $\tilde z$, then $(\tilde z,p)\in \Dom(*)$ and $\tilde z*p=q$.
\end{enumerate} 
\item If $(\tilde z,p)\in \Dom(*)$, then $(\tilde w,\tilde z*p)\in \Dom(*)$ if and only if $(\tilde w\tilde z,p)\in \Dom(*)$, and in this case, 
\[\tilde w*(\tilde z*p)=(\tilde w\tilde z)*p\]
\end{enumerate}

\end{defna2}

Note that an action of $\ex T^{n}$ on $\ex B$ is not quite an action of $\ex T^{n}$ considered as a group, simply  because the action may not be defined everywhere. (This may be thought of as similar to a not always defined $\mathbb R$ action given by the flow of a vectorfield.) For a given $p\in\ex B$, the action is always locally defined in the sense that $(\tilde z,p)\in\Dom (*)$ for $\tilde z$ close enough to $1$. 

\begin{lemmaa}\label{equivariant coordinates} Every (second countable) exploded manifold $\ex B$ has a cover by coordinate charts modeled on open subsets of $\mathbb R^{n}\times \et mP$ so that every transition map is either  equivariant or its inverse is equivariant.

Every family $\pi: \ex B\longrightarrow\ex G$ of exploded manifolds has a cover by coordinate charts so that transition maps on $\ex G$ and $\ex B$ are equivariant in the above sense,  and $\pi$ is equivariant in the sense that given any of our coordinate charts $U$ on $\ex B$, one of our coordinate charts $V$ on $\ex G$ contains $\pi(U)$, and $\pi:U\longrightarrow V$ is equivariant. 
\end{lemmaa}

\pf

As we are assuming $\ex B$ is second countable, we may choose an exhaustion of $\ex B$ by a sequence of compact subsets $K_{i}$ so that $K_{i-1}$ is contained in the interior of $K_{i}$. As any exploded manifold is locally basic, we may cover $\ex B$ by open subsets $W'$ with closures contained in basic open subsets $W$
which are themselves contained inside some $K_{i+1}\setminus K_{i-1}$. Then as each $K_{i}$ is compact, we may choose a countable sub cover by $W'_{l}$ so that only finitely many of the corresponding $W_{l}$ are contained inside any $K_{i}$, so each $W_{l}$ intersects only finitely many $W_{j}$. We shall choose all coordinate charts used in our construction  to have the following property:
\begin{enumerate}
\item\label{property}Each coordinate chart $V$ shall be  contained in some $W'_{l}$, and also be contained inside $W_{j}$ if  $V\cap W'_{j}\neq \emptyset$.\end{enumerate}

 Note that every point in $\ex B$ is contained in a coordinate chart satisfying property \ref{property}, and the union  of all coordinate charts satisfying property \ref{property} that intersect a given coordinate chart satisfying property \ref{property} is contained in some $W_{i}$, and is therefore basic. This simplifies the construction of equivariant transition maps because if the union of two coordinate charts is not basic, the tropical part of their transition maps may be different on different connected components.
 
\

Suppose that we have already constructed some coordinate charts satisfying property \ref{property} with equivariant transition maps which cover an open subset $O\subset \ex B$. 
%


 Let $\psi:V\longrightarrow \mathbb R^{n}\times\et m{P}$ be a coordinate chart satisfying property \ref{property} which restricts to a coordinate chart $\psi^{\circ}:U\longrightarrow \mathbb R^{n}\times\et m{P^{\circ}}$ on a stratum $\ex B_{i}$ of $\ex B$, where $P^{\circ}\subset P$ is the interior of $P$, and $U$ is $V\cap \ex B_{i}$.  This coordinate chart $\psi^{\circ}$ is automatically compatible with all previous coordinate charts because $P^{\circ}$ is an open polytope.   We shall now try to extend $\psi^{\circ}$ to a coordinate chart $\psi':V'\longrightarrow \mathbb R^{n}\times \et mP$ on an open neighborhood $V'$ of $U$ in $\ex B$. 
 
 \begin{claima}\label{claimA4} Given any open subset $O'$ of $\ex B$ so that the closure of $O'$ is contained in the open set $O$ already covered by equivariant coordinate charts,  there exists an extension of the coordinate chart $\psi^{\circ}$ on $\ex B_{i}$ to a coordinate chart $\psi'$ on $\ex B$ satisfying property \ref{property} with equivariant transition maps to the old coordinate charts restricted to $O'$.
 \end{claima}
 
 To prove Claim A4, consider one of our previously defined coordinate charts $V_{0}$  which intersects $U$.
 As $\totb U$ is a $m$ dimensional open polytope, and the union of $V_{0}$  with $V\supset U$ is basic, the restriction of the action on $V_{0}$ to some  $\ex T^{m}$ subgroup is compatible with the action of $\ex T^{m}$ on $U$. (The basic property is used to rule out the possibility that $\totb{U}$ intersects $\totb{V_{0}}$ in more than one way, giving rise to incompatible $\ex T^{m}$ actions.) The set of points $p$ in $V_{0}$  so that there is some $\tilde z\in\ex T^{m}$ so that $\tilde z*p\in U$ is open. 
 Denote by $O_{U}\subset V$ the intersection with $V$ of the union of all such subsets of all our previously defined coordinate charts. Note that as our previously defined coordinate charts have equivariant transition maps, and the union of the coordinate charts which intersect $U$ is basic, there is a well defined action of $\ex T^{m}$ on $O_{U}$.

  Because every point in $O_{U}$ is related to a point in $U$ by the $\ex T^{m}$ action, there is a unique map $\psi': O_{U}\longrightarrow\mathbb R^{n}\times \et mP$ so that $\psi'=\psi^{\circ}$ restricted to $U$ and $\psi'(\tilde z*p)=\tilde z*\psi'(p)$ for all $(\tilde z,p)\in \Dom (*)$.  The new map $\psi'$ is an isomorphism onto its image because $\psi$ was. The chart $\psi'$ is appropriately compatible with all our previous coordinate charts, but it does not cover all of $U$, only $U\cap O$.  The goal is now to patch together $\psi$ and $\psi'$ to create a map defined on a neighborhood of all $U$, so that this patched together map agrees with $\psi'$ on $O'\subset O$. Patching together two such maps is a standard construction; we give the details in the following paragraph.
  
  There exists a smooth vectorfield $v$ on $\psi(O_{U})$ so that $v$ vanishes on $U$, and restricted to some open neighborhood of $U\cap O_{U}$ in $O_{U}\cap V$, $\psi'$ is given by $\psi$ composed with the geodesic flow of $v$ using the standard flat metric on $\mathbb R^{n}\times \et mP$.  To extend $\psi'$, we shall extend this vectorfield. As the closure of $O'$ intersected with $U$   is contained in $O\cap U=O_{U}\cap U$, there exists a neighborhood $N$ of $U$ so that the closure of $O'\cap N$ in $N$ is contained in $O_{U}\cap N$, and so that $N$ is small enough that on $N$, $\psi'$ is $\psi$ composed with the geodesic flow of $v$. Therefore, we may extend $v$ to a vector field on $\psi(N)$ by multiplying by a smooth cut off function which is $1$ in $\psi(O'\cap N)$ and $0$ outside of $\psi(O_{U}\cap N)$. Restricted to some neighborhood $V'$ of $U$ in $N$, the resulting map $\psi'$ given by composing $\psi$ restricted to $V'$ with the geodesic flow of $v$ will be an isomorphism on to its image. 
  
  This map $\psi':V'\longrightarrow \mathbb R^{n}\times\et mP$  is our required coordinate chart  which is compatible with all previous coordinate charts restricted to $O'$ and satisfies property \ref{property} so long as $V'$ was chosen small enough. Therefore Claim A4 holds.
  
  \
 
The remainder of the proof of Lemma A3 shall now follow by a double induction.  Our first inductive assumption is  that every stratum of $\ex B$ with dimension strictly less than $ k$ in $\totl{\ex B}$ can be covered by coordinate charts satisfying property \ref{property} with equivariant transition maps (this is trivially true for k=0). We shall prove that we can cover all strata of dimension less than or equal to $k$ in $\totl{\ex B}$ by coordinate charts satisfying property \ref{property} with equivariant transition maps. Our main tool shall be Claim A4, which requires we shrink the domain of definition of formerly constructed coordinate charts a little each time we define a new coordinate chart--- this requirement that we must shrink the domain of definition of previously constructed coordinate charts shall be the main complication that our inductive argument must overcome.

Assume that we already have our equivariant cover of the strata of lower dimension. Let $O$ denote the open set already covered, and let $O'$ be an open subset of $O$ which contains all strata of lower dimension so that the closure of $O'$ is contained in $O$. Then choose a nested sequence of open sets $O_{i}$ containing $O'$ so that $O_{0}=O$ and the closure of $O_{i+1}$ is contained inside $O_{i}$. 

 Denote by  $U_{0}$ the intersection of $O'$ with the strata of dimension $k$ in $\totl{\ex B}$. Choose a countable open cover  $\{U_{i}\}$ of the strata of dimension $k$ in $\totl{\ex B}$ so that for all $i>0$, $\bar U_{i}$ is contained inside some coordinate chart satisfying property \ref{property}.  Our second inductive assumption is  that for all $i<n$, there exist  coordinate charts $V_{i}$ with open subsets $V_{i,j}\subset V_{i}$ for all $j\geq i$ so that
 \begin{itemize} \item the coordinate charts $V_{i}$ satisfy property \ref{property},
 \item each of our open subsets $V_{i,j}\subset V_{i}$ contains the closure of some open neighborhood  $V'_{i}$ of $U_{i}$,
 \item   $\bar V_{i,j+1}\subset V_{i,j}$
 \item and all transition maps between  all $V_{i,n-1}$ for $i< n$ and all our previous coordinate charts restricted to $O_{n-1}$ are equivariant.
 \end{itemize}
 
 As the closure of   $O_{n}\cup\bigcup _{i=1}^{n-1} V_{i,n}$ is contained in $O_{n-1}\cup\bigcup_{i=1}^{n-1}V_{i,n-1}$,  Claim A4 implies that we can construct a coordinate chart $V_{n}$ satisfying property \ref{property} which contains the closure of some open neighborhood $V_{n}'$ of $U_{n}$ so that for $i<n$, all transition maps between $V_{n}$ and $V_{i,n}$ are equivariant, and so that all transition maps between $V_{n}$ and the restriction of all previous coordinate charts to $O_{n}$ are equivariant. We then may choose a sequence of open subsets $V_{n,j}\subset V_{n}$ for $j\geq n$, so that  $V_{n,j}$ contains $V'_{n}$ and so that $\bar V_{n,j+1}\subset V_{n,j}$. Then  our second inductive statement holds for $n$. Therefore, our second inductive statement holds for all $n$. Now consider our previously constructed coordinate charts restricted to $O'$ together with the new coordinate charts $V_{i}'$. The transition maps between these coordinate charts are all equivariant, and together these coordinate charts cover all strata of dimension $k$ in $\ex B$. By induction on $k$, we may cover all of $\ex B$ by coordinate charts satisfying property \ref{property} with equivariant transition maps. 

\

We can prove the family case similarly. First choose equivariant coordinate charts on $\ex G$. We can choose  coordinate charts on $\ex B$ so that $\pi$ is equivariant as follows: 

Around any point  $p\in\ex B$, we can choose a coordinate chart $V$ small enough that its image under $\pi$ is contained in one of our coordinate charts on $\ex G$, and $p$ is contained in the stratum of $V$ with maximal dimensional tropical part. On this interior stratum, our map is automatically equivariant, and we may extend this to an equivariant map $\pi_{1}$ from a neighborhood of $p$  in $V$ to our coordinate chart on $\ex G$. Suppose that our coordinate chart on $\ex G$ is an open subset of $\mathbb R^{n}\times\et mP$, so we may think of adding a vector in $\mathbb R^{n}$ to the first coordinates, and multiplying the last coordinates by an element of $(\mathbb C^{*})^{m}$. On a neighborhood of $p$, $\pi$ and $\pi_{1}$ are related by 

\[\pi_{1}(q)=h(q)*(\pi(q)+f(q))\]
where $f$ is a smooth $\mathbb R^{n}$ valued function which vanishes on the interior stratum of $V$, and $h$ is a smooth $(\mathbb C^{*})^{m}$ valued function which is $1$ on the interior stratum of $V$. There is therefore a smooth family $\pi_{t}$ of maps from some neighborhood of $p$ so that $\pi_{0}=\pi$ and $\pi_{1}$ is the equivariant map defined earlier. Lift each  standard basis vector field on $\mathbb R^{n}\times \et mP$ to a vector field on $V$. Then $\frac \partial {\partial t}\pi_{t}(q)$ may be considered as some linear combination of standard basis vectors on $\mathbb R^{n}\times \et mP$, which corresponds to the same linear combination of lifted vectorfields which is a vectorfield $X_{q,t}$ on $V$. We can define a family of maps $\psi_{t}$ from a neighborhood of $p$ to $V$ so that $\pi_{t}=\pi\circ\psi_{t}$ using the differential equation
\[\frac {\partial\psi_{t}}{\partial t} (q)=X_{q,t}(\psi_{t}(q))\]
with the initial condition $\psi_{0}(q)=q$. As this differential equation can be reformulated as the flow of a smooth  vector field on $[0,1]$ times the product of two neighborhoods of $p$ (one neighborhood for each appearance of $q$ in the right hand side of the above equation), its solution is smooth, and as $X_{q,t}$ always vanishes on the stratum of $V$ containing $p$, $\psi_{t}$ exists for all $t\in[0,1]$ restricted to some small enough neighborhood of $p$. As $\pi\circ \psi_{t}$ satisfies the same differential equation as $\pi_{t}$, $\pi\circ \psi_{t}=\pi_{t}$.
  As $\psi_{t}$ is the identity restricted to the stratum containing $p$, $\psi_{t}$ is the flow of a smooth time dependent vectorfield restricted to some neighborhood of $p$, so $\psi_{t}$ is an isomorphism onto its image restricted to this neighborhood. Therefore, $\psi_{1}$ restricted to this neighborhood gives a coordinate chart with an equivariant projection to our coordinate chart on $\ex G$.

At this stage we have shown that given a choice of equivariant coordiate charts on $\ex G$, we may choose coordinate charts on $\ex B$ so that $\pi$ is equivariant. Now we may repeat the above argument for $\ex B$ choosing all coordinate charts so that the map $\pi$ is equivariant.  Claim A4 also holds for coordinate charts with equivariant maps to some coordinate chart on $\ex G$, because the step where a coordinate chart needs to be modified using the flow of a vector field can  be achieved using vertical vectorfields so that this procedure preserves the property that $\pi$ is equivariant.
The remainder of the proof is identical.

\stop

\bibliographystyle{plain}
\bibliography{ref.bib}

\end{document}